

\ifx\shlhetal\undefinedcontrolsequence\let\shlhetal\relax\fi

\input amstex
\expandafter\ifx\csname mathdefs.tex\endcsname\relax
  \expandafter\gdef\csname mathdefs.tex\endcsname{}
\else \message{Hey!  Apparently you were trying to
  \string\input{mathdefs.tex} twice.   This does not make sense.} 
\errmessage{Please edit your file (probably \jobname.tex) and remove
any duplicate ``\string\input'' lines}\endinput\fi




\catcode`\X=12\catcode`\@=11

\def\n@wcount{\alloc@0\count\countdef\insc@unt}
\def\n@wwrite{\alloc@7\write\chardef\sixt@@n}
\def\n@wread{\alloc@6\read\chardef\sixt@@n}
\def\r@s@t{\relax}\def\v@idline{\par}\def\@mputate#1/{#1}
\def\l@c@l#1X{\firstpart.#1}\def\gl@b@l#1X{#1}\def\t@d@l#1X{{}}

\def\crossrefs#1{\ifx\all#1\let\tr@ce=\all\else\def\tr@ce{#1,}\fi
   \n@wwrite\cit@tionsout\openout\cit@tionsout=\jobname.cit 
   \write\cit@tionsout{\tr@ce}\expandafter\setfl@gs\tr@ce,}
\def\setfl@gs#1,{\def\@{#1}\ifx\@\empty\let\next=\relax
   \else\let\next=\setfl@gs\expandafter\xdef
   \csname#1tr@cetrue\endcsname{}\fi\next}
\def\m@ketag#1#2{\expandafter\n@wcount\csname#2tagno\endcsname
     \csname#2tagno\endcsname=0\let\tail=\all\xdef\all{\tail#2,}
   \ifx#1\l@c@l\let\tail=\r@s@t\xdef\r@s@t{\csname#2tagno\endcsname=0\tail}\fi
   \expandafter\gdef\csname#2cite\endcsname##1{\expandafter
     \ifx\csname#2tag##1\endcsname\relax?\else\csname#2tag##1\endcsname\fi
     \expandafter\ifx\csname#2tr@cetrue\endcsname\relax\else
     \write\cit@tionsout{#2tag ##1 cited on page \folio.}\fi}
   \expandafter\gdef\csname#2page\endcsname##1{\expandafter
     \ifx\csname#2page##1\endcsname\relax?\else\csname#2page##1\endcsname\fi
     \expandafter\ifx\csname#2tr@cetrue\endcsname\relax\else
     \write\cit@tionsout{#2tag ##1 cited on page \folio.}\fi}
   \expandafter\gdef\csname#2tag\endcsname##1{\expandafter
      \ifx\csname#2check##1\endcsname\relax
      \expandafter\xdef\csname#2check##1\endcsname{}%
      \else\immediate\write16{Warning: #2tag ##1 used more than once.}\fi
      \multit@g{#1}{#2}##1/X%
      \write\t@gsout{#2tag ##1 assigned number \csname#2tag##1\endcsname\space
      on page \number\count0.}%
   \csname#2tag##1\endcsname}}

\def\multit@g#1#2#3/#4X{\def\t@mp{#4}\ifx\t@mp\empty%
      \global\advance\csname#2tagno\endcsname by 1 
      \expandafter\xdef\csname#2tag#3\endcsname
      {#1\number\csname#2tagno\endcsnameX}%
   \else\expandafter\ifx\csname#2last#3\endcsname\relax
      \expandafter\n@wcount\csname#2last#3\endcsname
      \global\advance\csname#2tagno\endcsname by 1 
      \expandafter\xdef\csname#2tag#3\endcsname
      {#1\number\csname#2tagno\endcsnameX}
      \write\t@gsout{#2tag #3 assigned number \csname#2tag#3\endcsname\space
      on page \number\count0.}\fi
   \global\advance\csname#2last#3\endcsname by 1
   \def\t@mp{\expandafter\xdef\csname#2tag#3/}%
   \expandafter\t@mp\@mputate#4\endcsname
   {\csname#2tag#3\endcsname\lastpart{\csname#2last#3\endcsname}}\fi}
\def\t@gs#1{\def\all{}\m@ketag#1e\m@ketag#1s\m@ketag\t@d@l p
\let\realscite\scite
\let\realstag\stag
   \m@ketag\gl@b@l r \n@wread\t@gsin
   \openin\t@gsin=\jobname.tgs \re@der \closein\t@gsin
   \n@wwrite\t@gsout\openout\t@gsout=\jobname.tgs }
\outer\def\localtags{\t@gs\l@c@l}
\outer\def\globaltags{\t@gs\gl@b@l}
\outer\def\newlocaltag#1{\m@ketag\l@c@l{#1}}
\outer\def\newglobaltag#1{\m@ketag\gl@b@l{#1}}

\newif\ifpr@ 
\def\m@kecs #1tag #2 assigned number #3 on page #4.%
   {\expandafter\gdef\csname#1tag#2\endcsname{#3}
   \expandafter\gdef\csname#1page#2\endcsname{#4}
   \ifpr@\expandafter\xdef\csname#1check#2\endcsname{}\fi}
\def\re@der{\ifeof\t@gsin\let\next=\relax\else
   \read\t@gsin to\t@gline\ifx\t@gline\v@idline\else
   \expandafter\m@kecs \t@gline\fi\let \next=\re@der\fi\next}
\def\pretags#1{\pr@true\pret@gs#1,,}
\def\pret@gs#1,{\def\@{#1}\ifx\@\empty\let\n@xtfile=\relax
   \else\let\n@xtfile=\pret@gs \openin\t@gsin=#1.tgs \message{#1} \re@der 
   \closein\t@gsin\fi \n@xtfile}

\newcount\sectno\sectno=0\newcount\subsectno\subsectno=0
\newif\ifultr@local \def\ultralocal{\ultr@localtrue}
\def\firstpart{\number\sectno}
\def\lastpart#1{\ifcase#1 \or a\or b\or c\or d\or e\or f\or g\or h\or 
   i\or k\or l\or m\or n\or o\or p\or q\or r\or s\or t\or u\or v\or w\or 
   x\or y\or z \fi}

\def\resetall{\global\advance\sectno by 1\subsectno=0
   \gdef\firstpart{\number\sectno}\r@s@t}
\def\resetsub{\global\advance\subsectno by 1
   \gdef\firstpart{\number\sectno.\number\subsectno}\r@s@t}
\def\newsection#1\par{\resetall\vskip0pt plus.3\vsize\penalty-250
   \vskip0pt plus-.3\vsize\bigskip\bigskip
   \message{#1}\leftline{\bf#1}\nobreak\bigskip}
\def\subsection#1\par{\ifultr@local\resetsub\fi
   \vskip0pt plus.2\vsize\penalty-250\vskip0pt plus-.2\vsize
   \bigskip\smallskip\message{#1}\leftline{\bf#1}\nobreak\medskip}


\newdimen\marginshift

\newdimen\margindelta
\newdimen\marginmax
\newdimen\marginmin

\def\margininit{       
\marginmax=3 true cm                  
				      
\margindelta=0.1 true cm              
\marginmin=0.1true cm                 
\marginshift=\marginmin
}    

\def\t@gsjj#1,{\def\@{#1}\ifx\@\empty\let\next=\relax\else\let\next=\t@gsjj
   \def\@@{p}\ifx\@\@@\else
   \expandafter\gdef\csname#1cite\endcsname##1{\citejj{##1}}
   \expandafter\gdef\csname#1page\endcsname##1{?}
   \expandafter\gdef\csname#1tag\endcsname##1{\tagjj{##1}}\fi\fi\next}
\newif\ifshowstuffinmargin
\showstuffinmarginfalse
\def\jjtags{\showstuffinmargintrue
\ifx\all\relax\else\expandafter\t@gsjj\all,\fi}

\def\tagjj#1{\realstag{#1}\mginpar{\zeigen{#1}}}
\def\citejj#1{\zeigen{#1}\mginpar{\rechnen{#1}}}

\def\rechnen#1{\expandafter\ifx\csname stag#1\endcsname\relax ??\else
                           \csname stag#1\endcsname\fi}

\newdimen\theight

\def\marginfont{\sevenrm}

\def\trymarginbox#1{\setbox0=\hbox{\marginfont\hskip\marginshift #1}%
		\global\marginshift\wd0 
		\global\advance\marginshift\margindelta}

\def \mginpar#1{%
\ifvmode\setbox0\hbox to \hsize{\hfill\rlap{\marginfont\quad#1}}%
\ht0 0cm
\dp0 0cm
\box0\vskip-\baselineskip
\else 
             \vadjust{\trymarginbox{#1}%
		\ifdim\marginshift>\marginmax \global\marginshift\marginmin
			\trymarginbox{#1}%
                \fi
             \theight=\ht0
             \advance\theight by \dp0    \advance\theight by \lineskip
             \kern -\theight \vbox to \theight{\rightline{\rlap{\box0}}%
\vss}}\fi}


\def\t@gsoff#1,{\def\@{#1}\ifx\@\empty\let\next=\relax\else\let\next=\t@gsoff
   \def\@@{p}\ifx\@\@@\else
   \expandafter\gdef\csname#1cite\endcsname##1{\zeigen{##1}}
   \expandafter\gdef\csname#1page\endcsname##1{?}
   \expandafter\gdef\csname#1tag\endcsname##1{\zeigen{##1}}\fi\fi\next}
\def\verbatimtags{\showstuffinmarginfalse
\ifx\all\relax\else\expandafter\t@gsoff\all,\fi}
\def\zeigen#1{\hbox{$\langle$}#1\hbox{$\rangle$}}
\def\margincite#1{\ifshowstuffinmargin\mginpar{\rechnen{#1}}\fi}

\def\(#1){\edef\dot@g{\ifmmode\ifinner(\hbox{\noexpand\etag{#1}})
   \else\noexpand\eqno(\hbox{\noexpand\etag{#1}})\fi
   \else(\noexpand\ecite{#1})\fi}\dot@g}

\newif\ifbr@ck
\def\eat#1{}
\def\[#1]{\br@cktrue[\br@cket#1'X]}
\def\br@cket#1'#2X{\def\temp{#2}\ifx\temp\empty\let\next\eat
   \else\let\next\br@cket\fi
   \ifbr@ck\br@ckfalse\br@ck@t#1,X\else\br@cktrue#1\fi\next#2X}
\def\br@ck@t#1,#2X{\def\temp{#2}\ifx\temp\empty\let\neext\eat
   \else\let\neext\br@ck@t\def\temp{,}\fi
   \def\teemp{#1}\ifx\teemp\empty\else\rcite{#1}\fi\temp\neext#2X}
\def\resetbr@cket{\gdef\[##1]{[\rtag{##1}]}}
\def\references{\resetbr@cket\newsection References\par}

\newtoks\symb@ls\newtoks\s@mb@ls\newtoks\p@gelist\n@wcount\ftn@mber
    \ftn@mber=1\newif\ifftn@mbers\ftn@mbersfalse\newif\ifbyp@ge\byp@gefalse
\def\defm@rk{\ifftn@mbers\n@mberm@rk\else\symb@lm@rk\fi}
\def\n@mberm@rk{\xdef\m@rk{{\the\ftn@mber}}%
    \global\advance\ftn@mber by 1 }
\def\rot@te#1{\let\temp=#1\global#1=\expandafter\r@t@te\the\temp,X}
\def\r@t@te#1,#2X{{#2#1}\xdef\m@rk{{#1}}}
\def\b@@st#1{{$^{#1}$}}\def\str@p#1{#1}
\def\symb@lm@rk{\ifbyp@ge\rot@te\p@gelist\ifnum\expandafter\str@p\m@rk=1 
    \s@mb@ls=\symb@ls\fi\write\f@nsout{\number\count0}\fi \rot@te\s@mb@ls}
\def\byp@ge{\byp@getrue\n@wwrite\f@nsin\openin\f@nsin=\jobname.fns 
    \n@wcount\currentp@ge\currentp@ge=0\p@gelist={0}
    \re@dfns\closein\f@nsin\rot@te\p@gelist
    \n@wread\f@nsout\openout\f@nsout=\jobname.fns }
\def\m@kelist#1X#2{{#1,#2}}
\def\re@dfns{\ifeof\f@nsin\let\next=\relax\else\read\f@nsin to \f@nline
    \ifx\f@nline\v@idline\else\let\t@mplist=\p@gelist
    \ifnum\currentp@ge=\f@nline
    \global\p@gelist=\expandafter\m@kelist\the\t@mplistX0
    \else\currentp@ge=\f@nline
    \global\p@gelist=\expandafter\m@kelist\the\t@mplistX1\fi\fi
    \let\next=\re@dfns\fi\next}
\def\symbols#1{\symb@ls={#1}\s@mb@ls=\symb@ls} 
\def\bigsymbol{\textstyle}
\symbols{\bigsymbol\ast,\dagger,\ddagger,\sharp,\flat,\natural,\star}
\def\ftnumbers{\ftn@mberstrue} \def\ftsymbols{\ftn@mbersfalse}
\def\paginal{\byp@ge} \def\resetftnumbers{\ftn@mber=1}
\def\ftnote#1{\defm@rk\expandafter\expandafter\expandafter\footnote
    \expandafter\b@@st\m@rk{#1}}

\long\def\jump#1\endjump{}
\def\ssum{\mathop{\lower .1em\hbox{$\textstyle\Sigma$}}\nolimits}

\def\qed{\nobreak\kern 1em \vrule height .5em width .5em depth 0em}
\def\newneq{\hbox{\rlap{\hbox to 1\wd9{\hss$=$\hss}}\raise .1em 
   \hbox to 1\wd9{\hss$\scriptscriptstyle/$\hss}}}
\def\subsetne{\setbox9 = \hbox{$\subset$}\mathrel{\hbox{\rlap
   {\lower .4em \newneq}\raise .13em \hbox{$\subset$}}}}
\def\supsetne{\setbox9 = \hbox{$\subset$}\mathrel{\hbox{\rlap
   {\lower .4em \newneq}\raise .13em \hbox{$\supset$}}}}

\def\vbar{\mathchoice{\vrule height6.3ptdepth-.5ptwidth.8pt\kern-.8pt}
   {\vrule height6.3ptdepth-.5ptwidth.8pt\kern-.8pt}
   {\vrule height4.1ptdepth-.35ptwidth.6pt\kern-.6pt}
   {\vrule height3.1ptdepth-.25ptwidth.5pt\kern-.5pt}}
\def\f@dge{\mathchoice{}{}{\mkern.5mu}{\mkern.8mu}}
\def\b@c#1#2{{\rm \mkern#2mu\vbar\mkern-#2mu#1}}
\def\b@b#1{{\rm I\mkern-3.5mu #1}}
\def\b@a#1#2{{\rm #1\mkern-#2mu\f@dge #1}}
\def\bb#1{{\count4=`#1 \advance\count4by-64 \ifcase\count4\or\b@a A{11.5}\or
   \b@b B\or\b@c C{5}\or\b@b D\or\b@b E\or\b@b F \or\b@c G{5}\or\b@b H\or
   \b@b I\or\b@c J{3}\or\b@b K\or\b@b L \or\b@b M\or\b@b N\or\b@c O{5} \or
   \b@b P\or\b@c Q{5}\or\b@b R\or\b@a S{8}\or\b@a T{10.5}\or\b@c U{5}\or
   \b@a V{12}\or\b@a W{16.5}\or\b@a X{11}\or\b@a Y{11.7}\or\b@a Z{7.5}\fi}}

\catcode`\X=11 \catcode`\@=12


\expandafter\ifx\csname citeadd.tex\endcsname\relax
\expandafter\gdef\csname citeadd.tex\endcsname{}
\else \message{Hey!  Apparently you were trying to
\string\input{citeadd.tex} twice.   This does not make sense.} 
\errmessage{Please edit your file (probably \jobname.tex) and remove
any duplicate ``\string\input'' lines}\endinput\fi

\sectno=-1   
\localtags
\NoBlackBoxes
\define\mr{\medskip\roster}
\define\sn{\smallskip\noindent}
\define\mn{\medskip\noindent}
\define\bn{\bigskip\noindent}
\define\ub{\underbar}
\define\wilog{\text{without loss of generality}}
\define\ermn{\endroster\medskip\noindent}
\define\dbca{\dsize\bigcap}
\define\dbcu{\dsize\bigcup}
\define \nl{\newline}
\magnification=\magstep 1
\documentstyle{amsppt}

{    
\catcode`@11

\ifx\alicetwothousandloaded@\relax
  \endinput\else\global\let\alicetwothousandloaded@\relax\fi

\gdef\subjclass{\let\savedef@\subjclass
 \def\subjclass##1\endsubjclass{\let\subjclass\savedef@
   \toks@{\def\usualspace{{\rm\enspace}}\eightpoint}%
   \toks@@{##1\unskip.}%
   \edef\thesubjclass@{\the\toks@
     \frills@{{\noexpand\rm2000 {\noexpand\it Mathematics Subject
       Classification}.\noexpand\enspace}}%
     \the\toks@@}}%
  \nofrillscheck\subjclass}
} 

\pageheight{8.5truein}
\topmatter
\title{The null ideal restricted to some non-null set may be
$\aleph_1$-saturated} \endtitle
\rightheadtext{The null ideal restricted, etc.}
\author {Saharon Shelah \thanks {\null\newline 
I would like to thank Alice Leonhardt for the beautiful typing. \null\newline
This research was supported by The Israel Science Foundation founded
by the Israel Academy of Sciences and Humanities. Publication 619 \null\newline
} \endthanks} \endauthor 
\affil{Institute of Mathematics\\
 The Hebrew University\\
 Jerusalem, Israel
 \medskip
 Rutgers University\\
 Mathematics Department\\
 New Brunswick, NJ  USA} \endaffil

\abstract  Our main result is that possibly some non-null set of reals cannot
be divided to uncountably many non-null sets.  We deal also with a non-null
set of reals, the graph of any function from it is null and deal with our
iterations somewhat more generally. \endabstract
\endtopmatter
\document  

\expandafter\ifx\csname alice2jlem.tex\endcsname\relax
  \expandafter\xdef\csname alice2jlem.tex\endcsname{\the\catcode`@}
\else \message{Hey!  Apparently you were trying to
\string\input{alice2jlem.tex}  twice.   This does not make sense.}
\errmessage{Please edit your file (probably \jobname.tex) and remove
any duplicate ``\string\input'' lines}\endinput\fi

\expandafter\ifx\csname bib4plain.tex\endcsname\relax
  \expandafter\gdef\csname bib4plain.tex\endcsname{}
\else \message{Hey!  Apparently you were trying to \string\input
  bib4plain.tex twice.   This does not make sense.}
\errmessage{Please edit your file (probably \jobname.tex) and remove
any duplicate ``\string\input'' lines}\endinput\fi

\def\renewcommand{\newcommand}	       
\edef\cite{\the\catcode`@}%
\catcode`@ = 11
\let\@oldatcatcode = \cite
\chardef\@letter = 11
\chardef\@other = 12
%
%
%
%
\def\@innerdef#1#2{\edef#1{\expandafter\noexpand\csname #2\endcsname}}%
%
%
\@innerdef\@innernewcount{newcount}%
\@innerdef\@innernewdimen{newdimen}%
\@innerdef\@innernewif{newif}%
\@innerdef\@innernewwrite{newwrite}%
%
%
%
\def\@gobble#1{}%
%
%
%
\ifx\inputlineno\@undefined
   \let\@linenumber = \empty 
\else
   \def\@linenumber{\the\inputlineno:\space}%
\fi
%
%
%
\def\@futurenonspacelet#1{\def\cs{#1}%
   \afterassignment\@stepone\let\@nexttoken=
}%
\begingroup 
\def\\{\global\let\@stoken= }%
\\ 
\endgroup
\def\@stepone{\expandafter\futurelet\cs\@steptwo}%
\def\@steptwo{\expandafter\ifx\cs\@stoken\let\@@next=\@stepthree
   \else\let\@@next=\@nexttoken\fi \@@next}%
\def\@stepthree{\afterassignment\@stepone\let\@@next= }%
%
%
%
\def\@getoptionalarg#1{%
   \let\@optionaltemp = #1%
   \let\@optionalnext = \relax
   \@futurenonspacelet\@optionalnext\@bracketcheck
}%
%
%
\def\@bracketcheck{%
   \ifx [\@optionalnext
      \expandafter\@@getoptionalarg
   \else
      \let\@optionalarg = \empty
      \expandafter\@optionaltemp
   \fi
}%
\def\@@getoptionalarg[#1]{%
   \def\@optionalarg{#1}%
   \@optionaltemp
}%
%
%
%
\def\@nnil{\@nil}%
\def\@fornoop#1\@@#2#3{}%
\def\@for#1:=#2\do#3{%
   \edef\@fortmp{#2}%
   \ifx\@fortmp\empty \else
      \expandafter\@forloop#2,\@nil,\@nil\@@#1{#3}%
   \fi
}%
\def\@forloop#1,#2,#3\@@#4#5{\def#4{#1}\ifx #4\@nnil \else
       #5\def#4{#2}\ifx #4\@nnil \else#5\@iforloop #3\@@#4{#5}\fi\fi
}%
\def\@iforloop#1,#2\@@#3#4{\def#3{#1}\ifx #3\@nnil
       \let\@nextwhile=\@fornoop \else
      #4\relax\let\@nextwhile=\@iforloop\fi\@nextwhile#2\@@#3{#4}%
}%
%
%
%
\@innernewif\if@fileexists
\def\@testfileexistence{\@getoptionalarg\@finishtestfileexistence}%
\def\@finishtestfileexistence#1{%
   \begingroup
      \def\extension{#1}%
      \immediate\openin0 =
         \ifx\@optionalarg\empty\jobname\else\@optionalarg\fi
         \ifx\extension\empty \else .#1\fi
         \space
      \ifeof 0
         \global\@fileexistsfalse
      \else
         \global\@fileexiststrue
      \fi
      \immediate\closein0
   \endgroup
}%
%
%
%
%
\def\bibliographystyle#1{%
   \@readauxfile
   \@writeaux{\string\bibstyle{#1}}%
}%
\let\bibstyle = \@gobble
%
%
\let\bblfilebasename = \jobname
\def\bibliography#1{%
   \@readauxfile
   \@writeaux{\string\bibdata{#1}}%
   \@testfileexistence[\bblfilebasename]{bbl}%
   \if@fileexists
      \nobreak
      \@readbblfile
   \fi
}%
\let\bibdata = \@gobble
%
%
\def\nocite#1{%
   \@readauxfile
   \@writeaux{\string\citation{#1}}%
}%
\@innernewif\if@notfirstcitation
%
%
\def\cite{\@getoptionalarg\@cite}%
%
%
\def\@cite#1{%
   \let\@citenotetext = \@optionalarg
   \printcitestart
   \nocite{#1}%
   \@notfirstcitationfalse
   \@for \@citation :=#1\do
   {%
      \expandafter\@onecitation\@citation\@@
   }%
   \ifx\empty\@citenotetext\else
      \printcitenote{\@citenotetext}%
   \fi
   \printcitefinish
}%
\def\@onecitation#1\@@{%
   \if@notfirstcitation
      \printbetweencitations
   \fi
   \expandafter \ifx \csname\@citelabel{#1}\endcsname \relax
      \if@citewarning
         \message{\@linenumber Undefined citation `#1'.}%
      \fi
      \expandafter\gdef\csname\@citelabel{#1}\endcsname{%
\strut
\vadjust{\vskip-\dp\strutbox
\vbox to 0pt{\vss\parindent0cm \leftskip=\hsize 
\advance\leftskip3mm
\advance\hsize 4cm\strut\openup-4pt 
\rightskip 0cm plus 1cm minus 0.5cm ?  #1 ?\strut}}
         {\tt
            \escapechar = -1
            \nobreak\hskip0pt
            \expandafter\string\csname#1\endcsname
            \nobreak\hskip0pt
         }%
      }%
   \fi
   \csname\@citelabel{#1}\endcsname
   \@notfirstcitationtrue
}%
%
%
\def\@citelabel#1{b@#1}%
%
%
\def\@citedef#1#2{\expandafter\gdef\csname\@citelabel{#1}\endcsname{#2}}%
%
%
%
\def\@readbblfile{%
   \ifx\@itemnum\@undefined
      \@innernewcount\@itemnum
   \fi
   \begingroup
      \def\begin##1##2{%
         \setbox0 = \hbox{\biblabelcontents{##2}}%
         \biblabelwidth = \wd0
      }%
      \def\end##1{}
      %
      %
      \@itemnum = 0
      \def\bibitem{\@getoptionalarg\@bibitem}%
      \def\@bibitem{%
         \ifx\@optionalarg\empty
            \expandafter\@numberedbibitem
         \else
            \expandafter\@alphabibitem
         \fi
      }%
      \def\@alphabibitem##1{%
         \expandafter \xdef\csname\@citelabel{##1}\endcsname {\@optionalarg}%
         \ifx\biblabelprecontents\@undefined
            \let\biblabelprecontents = \relax
         \fi
         \ifx\biblabelpostcontents\@undefined
            \let\biblabelpostcontents = \hss
         \fi
         \@finishbibitem{##1}%
      }%
      \def\@numberedbibitem##1{%
         \advance\@itemnum by 1
         \expandafter \xdef\csname\@citelabel{##1}\endcsname{\number\@itemnum}%
         \ifx\biblabelprecontents\@undefined
            \let\biblabelprecontents = \hss
         \fi
         \ifx\biblabelpostcontents\@undefined
            \let\biblabelpostcontents = \relax
         \fi
         \@finishbibitem{##1}%
      }%
      \def\@finishbibitem##1{%
         \biblabelprint{\csname\@citelabel{##1}\endcsname}%
         \@writeaux{\string\@citedef{##1}{\csname\@citelabel{##1}\endcsname}}%
         \ignorespaces
      }%
      %
      %
      \let\em = \bblem
      \let\newblock = \bblnewblock
      \let\sc = \bblsc
      \frenchspacing
      \clubpenalty = 4000 \widowpenalty = 4000
      \tolerance = 10000 \hfuzz = .5pt
      \everypar = {\hangindent = \biblabelwidth
                      \advance\hangindent by \biblabelextraspace}%
      \bblrm
      \parskip = 1.5ex plus .5ex minus .5ex
      \biblabelextraspace = .5em
      \bblhook
      \input \bblfilebasename.bbl
   \endgroup
}%
%
%
\@innernewdimen\biblabelwidth
\@innernewdimen\biblabelextraspace
%
%
%
\def\biblabelprint#1{%
   \noindent
   \hbox to \biblabelwidth{%
      \biblabelprecontents
      \biblabelcontents{#1}%
      \biblabelpostcontents
   }%
   \kern\biblabelextraspace
}%
%
%
%
\def\biblabelcontents#1{{\bblrm [#1]}}%
%
%
\def\bblrm{\rm}%
%
%
\def\bblem{\it}%
%
%
\def\bblsc{\ifx\@scfont\@undefined
              \font\@scfont = cmcsc10
           \fi
           \@scfont
}%
%
%
\def\bblnewblock{\hskip .11em plus .33em minus .07em }%
%
%
\let\bblhook = \empty
%
%
%
\def\printcitestart{[}
\def\printcitefinish{]}
\def\printbetweencitations{, }
\def\printcitenote#1{, #1}
%
%
%
\let\citation = \@gobble
%
%
%
\@innernewcount\@numparams
%
%
\def\newcommand#1{%
   \def\@commandname{#1}%
   \@getoptionalarg\@continuenewcommand
}%
%
%
\def\@continuenewcommand{%
   \@numparams = \ifx\@optionalarg\empty 0\else\@optionalarg \fi \relax
   \@newcommand
}%
%
%
\def\@newcommand#1{%
   \def\@startdef{\expandafter\edef\@commandname}%
   \ifnum\@numparams=0
      \let\@paramdef = \empty
   \else
      \ifnum\@numparams>9
         \errmessage{\the\@numparams\space is too many parameters}%
      \else
         \ifnum\@numparams<0
            \errmessage{\the\@numparams\space is too few parameters}%
         \else
            \edef\@paramdef{%
               \ifcase\@numparams
                  \empty  No arguments.
               \or ####1%
               \or ####1####2%
               \or ####1####2####3%
               \or ####1####2####3####4%
               \or ####1####2####3####4####5%
               \or ####1####2####3####4####5####6%
               \or ####1####2####3####4####5####6####7%
               \or ####1####2####3####4####5####6####7####8%
               \or ####1####2####3####4####5####6####7####8####9%
               \fi
            }%
         \fi
      \fi
   \fi
   \expandafter\@startdef\@paramdef{#1}%
}%
%
%
%
%
\def\@readauxfile{%
   \if@auxfiledone \else 
      \global\@auxfiledonetrue
      \@testfileexistence{aux}%
      \if@fileexists
         \begingroup
            \endlinechar = -1
            \catcode`@ = 11
            \input \jobname.aux
         \endgroup
      \else
         \message{\@undefinedmessage}%
         \global\@citewarningfalse
      \fi
      \immediate\openout\@auxfile = \jobname.aux
   \fi
}%
%
%
\newif\if@auxfiledone
\ifx\noauxfile\@undefined \else \@auxfiledonetrue\fi
%
%
%
%
\@innernewwrite\@auxfile
\def\@writeaux#1{\ifx\noauxfile\@undefined \write\@auxfile{#1}\fi}%
%
%
%
\ifx\@undefinedmessage\@undefined
   \def\@undefinedmessage{No .aux file; I won't give you warnings about
                          undefined citations.}%
\fi
%
%
\@innernewif\if@citewarning
\ifx\noauxfile\@undefined \@citewarningtrue\fi
%
%
%
\catcode`@ = \@oldatcatcode


\def\widestnumber#1#2{}

\def\rm{\fam0 \tenrm}

\def\fakesubhead#1\endsubhead{\bigskip\noindent{\bf#1}\par}



%
%
%

%

\font\textrsfs=rsfs10
\font\scriptrsfs=rsfs7
\font\scriptscriptrsfs=rsfs5

\newfam\rsfsfam
\textfont\rsfsfam=\textrsfs
\scriptfont\rsfsfam=\scriptrsfs
\scriptscriptfont\rsfsfam=\scriptscriptrsfs

\edef\oldcatcodeofat{\the\catcode`\@}
\catcode`\@11

\def\Cal@@#1{\noaccents@ \fam \rsfsfam #1}

\catcode`\@\oldcatcodeofat


\expandafter\ifx \csname margininit\endcsname \relax\else\margininit\fi

\newpage

\head {Anotated Content} \endhead  \resetall 
\bn
\S0 Introduction
\mr
\item "{{}}"  [We review results and background, and give notation.]
\endroster
\bn
\S1 The null ideal resricted to a non-null set may be $\aleph_1$-saturated 
\mr
\item "{{}}"  [We explain the difficulty for the null case, solved by adding
$\kappa$ random reals $\eta_{\lambda + \alpha}$ in the end, $\kappa$
measurable, but they are random over some subuniverses, so we add in the
beginning $\lambda$ Cohens ${\underset\tilde {}\to r_i}(i < \lambda)$.  The
memory is devised such that $\{A:\{\eta_{\lambda + \alpha}:\alpha \in A\}$ null$\}$
will include a ${\Cal P}(\kappa)^{\bold V} \backslash D,D$ a normal ultrafilter on
$\kappa$.  Whereas in \cite{Sh:592}, the $A_\alpha$ (memory) were chosen closed
enough, here we use automorphism of the memory structure.]
\endroster
\bn
\S2 Non-null set with no non-null function
\mr
\item "{{}}"  [We show that consistently for some non-null set $A$ of reals
for every function from $A$ to the reals, its graph is a null subset of the
plain.]
\endroster
\bn
\S3 The $L_{\aleph_1,\aleph_1}$-elmentary submodels and the forcing 
\mr
\item "{{}}"  [We deal with general FS iterations of c.c.c. definable in
subuniverses, and give sufficient condition for $P'_A \lessdot P_\alpha$.
As application we show consistency of some values of $2^{\aleph_0}$, 
add(meagre), cov(meagre), unif(meagre) ${\frak b}, {\frak d}$.]
\endroster
\newpage

\head {\S0 Introduction} \endhead  \resetall \sectno=0
\bn
Note that the result stated in the abstract tells us that the positive
result explained below cannot be improved (to $\aleph_1$ sets).  It is (Gitik, Shelah 
\cite{GiSh:582})
\mr
\item "{$(*)$}"  given sets $A_n$ of reals for $n < \omega$, we can find
$B_n \subseteq A_n$, pairwise disjoint such that $A_n,B_n$ have the same outer
Lebesgue measure.
\endroster
\bn
Lately, we have proved (\cite{Sh:592}): 
\proclaim{\stag{0.1} Theorem}  Con(cov(null) = $\aleph_\omega + MA_{\aleph_n}$)
for each $n < \omega$.
\endproclaim
\bigskip

The idea of the proof was to use finite support $\langle P_i,Q_i:i < \alpha \rangle$, 
where say ${\underset\tilde {}\to Q_i}$ has generic real 
${\underset\tilde {}\to r_i}$ and $Q_i$ is random forcing in $\bold V
[\langle {\underset\tilde {}\to r_j}:j \in a_i \rangle],a_i \subseteq i,a_i$
closed enough or ${\underset\tilde {}\to r_i}$ is Cohen real but the
``memory" is not transitive, i.e. $j \in a_i \nRightarrow a_j
\subseteq a_i$.  

As \scite{0.1}
was hard for me for long it seems reasonable to hope the solution will open
my eyes on other problems as well.

In this paper we deal mainly with ``can every non-null set be partitioned to
uncountably many non-null sets?", equivalently: ``can the ideal of null sets
which are subsets of a fixed non-null subset of $\Bbb R$ be 
$\aleph_1$-saturated?".  P. Komj\'ath \cite{Ko}, proved that it is consistent
that there is a non-meager set $A$ such that the ideal of meager subsets of
$A$ is $\aleph_1$-saturated.  The question whether a similar fact may hold for
measure dates back to Ulam, see also Prikry's thesis, Fremlin had asked both
versions since the seventies.  
\bn
So we prove the following:
\proclaim{\stag{0.2} Theorem}  It is consistent that there is a non-null set
$A \subseteq \Bbb R$ such that the ideal of null subsets of $A$ is
$\aleph_1$-saturated (of course, provided that ``ZFC $+ \exists$ measurable"
is consistent).
\endproclaim
\bigskip

The proof of \scite{0.1} was not directly applicable, but ``turning the
tables" make it relevant as explained in \S1.

The question appears on the current Fremlin's list of problems, \cite{Fe94}
as problem EL(a) respectively.

We also have some further remarks, e.g. the exact cardinal assumption for
\scite{0.1}.  We try to make the paper self-contained for readers with basic
knowledge of forcing and of \cite{Sh:592}.
\bn
Also we answer the following problem which Komj\'ath draws our attention to:
\proclaim{\stag{0.3} Theorem}  It is consistent that:
\mr
\item "{$\oplus$}"  there is a non-null $A \subseteq \Bbb R$ such that: for
every $f:A \rightarrow \Bbb R$, the function $f$ as a subset of the plane
$\Bbb R \times \Bbb R$ is null
\ermn
provided that ``ZFC + there is a measurable cardinal" is consistent.
\endproclaim
\bigskip

Lastly, in \S3 we investigate when partial memory iteration behaves as in
\cite{Sh:592} and gives an example how to apply the method.  This originally
was a part of \cite{Sh:592}, but as the main part of \cite{Sh:592} was in 
final form and this was not and we like to add \S3, we separate.
\bigskip

\demo{\stag{0.4} Notation}  We denote:
\mr
\item  natural numbers by $k,l,m,n$ and also $i,j$
\sn
\item  ordinals by $\alpha,\beta,\gamma,\delta,\zeta,\xi$ ($\delta$ always
limit)
\sn
\item  cardinals by $\lambda,\kappa,\chi,\mu$
\sn
\item  reals by $a,b$ and positive reals (normally small) by $\varepsilon$
\sn
\item  subsets of $\omega$ or ${}^{\omega \ge}2$ or Ord by $A,B,C,X,Y,Z$
\sn
\item  ${\Cal B}$ is a Borel function
\sn
\item  finitely additive measures by $\Xi$
\sn
\item  sequences of natural numbers or ordinals by $\eta,\nu,\rho$
\sn
\item  $s$ is used for various things.
\ermn
${\Cal T}$ is as in Definition 2.9 of \cite{Sh:592}, $t$ is a member
of ${\Cal T}$. \nl
We denote
\mr
\item  forcing notions by $P,Q$
\sn
\item  forcing conditions by $p,q$
\sn
\item  $p \ge q$ means $p$ stronger than $q$
\ermn
and use $r$ to denote members of Random (see below)
\mr
\item  Leb is Lebesgue measure (on $\{A:A \subseteq {}^\omega 2\}$)
\sn
\item  Random will be the family
$$
\align
\{r \subseteq {}^{\omega \ge}2:&r \text{ is a subtree of }
({}^{\omega >} 2,\triangleleft) \text{ (i.e. closed under initial
segments},<> \in r \\
  &\text{with no } \triangleleft \text{-maximal element (so lim}(r)
\subseteq {}^\omega 2 \\
  &\text{is closed)) and Leb(lim}(r)) > 0 \text{ and moreover } \eta \in r
\Rightarrow \text{ lim}(r^{[\eta]}) \\
  &\text{is not null (on }r^{[\eta]} \text{ see below)} \}
\endalign
$$
ordered by inverse inclusion.  We may sometimes use instead
$$
\{B:B \text{ is a Borel non-null subset of } {}^\omega 2\}.
$$
\ermn
For $\eta \in {}^{\omega >}2,A \subseteq {}^{\omega \ge}2$ let

$$
A^{[\eta]} = \{\nu \in A:\nu \triangleleft \eta \vee \eta \trianglelefteq
\nu\}.
$$

We thank Tomek Bartoszy\'nski and Mariusz Rabus and Heike Mildenberger for
reading and commenting and correcting.
\enddemo
\newpage

\head {\S1 The null ideal restricted to a non-null set may be
$\aleph_1$-saturated} \endhead  \resetall 
\bn
Let us first describe an outline of Komj\'ath's solution to the problem for
the meagre ideal.  Note that by a theorem of Solovay, the conclusion implies
that there is a measurable cardinal in some inner model.  So Komj\'ath starts
with

$$
\bold V \models ``\kappa \text{ is measurable and } D \text{ is a normal
ultrafilter on } \kappa".
$$
\mn
He uses finite support iteration $\langle P_\alpha,{\underset\tilde {}\to 
Q_\beta}:\alpha \le 2^\kappa,\beta < 2^\kappa \rangle$ of c.c.c.
forcing notions such that for $\beta < \kappa$ the forcing notion $Q_\beta$
adds a Cohen real ${\underset\tilde {}\to \eta_\beta} \in {}^\omega 2$ (so
${\underset\tilde {}\to Q_\beta} = ({}^{\omega >}2,\triangleleft))$ and for
each $\beta \in [\kappa,2^\kappa)$ for some $A_\beta \in D$, the forcing
notion ${\underset\tilde {}\to Q_\beta}$ makes the set
$\{{\underset\tilde {}\to \eta_\gamma}:\gamma \in \kappa \backslash A_\beta\}$
meagre (and every $A \in D$ appears).  The point is that finite support
iterations tend to preserve non-meagreness, so in $\bold V^{P_{2^\kappa}}$ 
the set $\{{\underset\tilde {}\to \eta_\alpha}:\alpha < \kappa\}$ remains 
non-meagre (and $D$ is (i.e. generates) in $\bold V^{P_{2^\kappa}}$ an
$\aleph_1$-saturated filter).

For our aims this per se is doomed to failure: finite support iterations add
Cohen reals which make the set of old reals null, whereas countable support
iterations tend to collapse $\aleph_2$.  This seems to indicate that the
solution should be delayed till we have better other support iterations (the
$2^{\aleph_0} = \aleph_3$ problem; see \cite[Ch.VII,Ch.VIII]{Sh:b} and
\cite[Ch.VII,Ch.VIII]{Sh:f}).  But we start with a simpler remedy: we use
finite support iteration $\langle P_\alpha,{\underset\tilde {}\to Q_\beta}:
\alpha \le \delta^*,\beta < \delta^* \rangle$ of c.c.c. forcing notions with
cf$(\delta^*) = \kappa,\langle \beta_\xi:\xi < \kappa \rangle$ increasing with
limit $\delta^*,{\underset\tilde {}\to Q_{\beta_\xi}}$ is a partial random,
say Random$^{\bold V^{P_{A(\beta_\xi)}}}$ where $A(\beta_\xi) =
A_{\beta_\xi} \subseteq \beta_\xi$,
adding the real ${\underset\tilde {}\to \tau_{\beta_\xi}}$. \nl
Clearly we want that

$$
\models_{P_{\delta^*}} ``\{{\underset\tilde {}\to \tau_{\beta_\xi}}:\xi <
\kappa\} \text{ is not null}".
$$
\mn
For this it suffices to have: every countable $A \subseteq \delta^*$ is
included in some $A(\beta_\xi)$.  However, we also want that for $B \in D$
the set $\{{\underset\tilde {}\to \tau_{\beta_\xi}}:\xi \in \kappa \backslash
B\}$ will be null.  For this it is natural to demand that for some $\alpha,
{\underset\tilde {}\to Q_\alpha}$ is Cohen forcing and

$$
(\forall \xi \in \kappa \backslash B)(\alpha \notin A_{\beta_\xi}).
$$
\mn
So we try to force a null set including $\{{\underset\tilde {}\to \tau_{\beta_\xi}}:
\xi \in \kappa \backslash B\}$ before we force the
${\underset\tilde {}\to \tau_\beta}$'s!  If 
${\underset\tilde {}\to \tau_\alpha}$ is similar enough to being Cohen over
$\langle {\underset\tilde {}\to \tau_{\beta_\xi}}:\xi \in \kappa \backslash
B \rangle$ we are done.   So this becomes similar to the problem in proving
\scite{0.1}.  But there we use $2^\kappa = \chi$, so that we need to carry
only $\kappa$ finitely additive measures $\langle
{\underset\tilde {}\to \Xi^t_\alpha}:t \in {\Cal T} \rangle$, i.e. have few
blueprints hence can make $A_\alpha$ (when $|Q_\alpha| \ge \kappa$) closed
enough.  However, the ${\underset\tilde {}\to \Xi^t_\alpha}$'s come only to
prove the existence of $p^\oplus$ which forces that many $p_\ell$'s
($\ell < \omega$) are in the generic set 
(see the proof of \cite[Claim 3.3]{Sh:592}).  
So in fact we can define the name of the finitely
additive measure after we have the sequence $\langle p_\ell:\ell < \omega
\rangle$.  Actually, we need this only for some specific cases and/or can
embed our $P_\alpha$ into another iteration.  So the problem boils down to
having the $A_\beta$'s closed enough to enable us to produce finitely
additive measures like the one in clause (i) of 2.11 or 
clause (d) of 2.16, \cite{Sh:592}.
The way we materialize the idea is by having enough automorphisms of the
structure

$$
\langle P_\alpha,{\underset\tilde {}\to Q_\beta},A_\beta,\mu_\beta,
{\underset\tilde {}\to \tau_\beta}:\alpha \le \alpha^*,\beta < \alpha^*
\rangle.
$$
\mn
In fact we can replace below $\lambda + \kappa$ by $\lambda \cdot \kappa$.
If we would like to have 
e.g. $\{A:\Xi(A) = 1\}$ to be a selective filter then the $\lambda
\times \kappa$ version is better: we can use 
$\{{\underset\tilde {}\to \tau_{\lambda \times \xi + 1 + \gamma}}:\gamma <
\lambda\}$, which are Cohen, for ensuring this.

If you do not like the use of ``automorphisms" of $\bar Q$ and doing it
through higher $\lambda$'s, later we analyze the method (of 
\cite[\S2,\S3]{Sh:592}) more fully (using essentially
$\prec_{L_{\aleph_1,\aleph_1}}$), and then the proof is more direct
(see \scite{3.1} - \scite{3.12}). 
\bigskip

\proclaim{\stag{1.1} Theorem}  Let $D$ be a $\kappa$-complete nonprincipal
ultrafilter on $\kappa$.  \ub{Then} for some c.c.c. forcing notion $P$ of
cardinality $2^\kappa$, in $\bold V^P$ we have
\mr
\item "{$(*)$}"  for some $A \in [{}^\omega 2]^\kappa$ we have
{\roster
\itemitem{ $(a)$ }  $A$ is not null
\sn
\itemitem{ $(b)$ }  the ideal $I = \{B \subseteq A:B \text{ is null}\}$ is
$\aleph_1$-saturated.
\endroster}
\ermn
Moreover,
\mr
\item "{$(**)$}"  we can find pairwise distinct $\eta_\xi \in {}^\omega 2$
(for $\xi < \kappa$) such that
\ermn
$$
A = \{\eta_\xi:\xi < \kappa\} \text{ and } Y \in D \Leftrightarrow
\{\eta_\xi:\xi \in \kappa \backslash Y\} \text{ is null}.
$$
\endproclaim
\bigskip

\remark{\stag{1.2} Remark}  1) We can replace ${}^\omega 2$ by $\Bbb R$. \nl
2) We can use as $D$ any uniform $\aleph_2$-complete filter on $\kappa$
(such that $D$ is $\aleph_1$-saturated in $\bold V$ and so in $\bold V^P$,
see \scite{1.10}).
\nl
3) In $(**)$ of course $D$ stands for the filter that $D$ generated in
$\bold V^P$.
\endremark
\bigskip

\demo{\stag{1.3} Convention}  For $\lambda \ge 2^\kappa$
(but we use only $\lambda < (2^\kappa)^{+ \omega}$) let $g_\lambda:\lambda
\rightarrow D$ be such that

$$
Y \in D \Rightarrow |\{\alpha < \lambda:g_\lambda(\alpha) = Y\}| = \lambda.
$$
\mn
For simplicity $g_\lambda$ is increasing with $\lambda$, and let $h:D
\rightarrow 2^\kappa$ be such that $g_\lambda \circ h  = \text{ id}_D$.  For
$\xi < \kappa$ let

$$
E_\xi = E^\lambda_\xi = \{\alpha < \lambda:\xi \in g_\lambda(\alpha)\}.
$$
\mn
Finally we let

$$
\bar E = \bar E^\lambda = \langle E^\lambda_\xi:\xi < \kappa \rangle.
$$
\mn
If not said otherwise we assume that $\lambda = \lambda^{\aleph_0} \ge
2^\kappa$ (for simplicity).
\sn
Note: our intention is that ${\underset\tilde {}\to \tau_\alpha}$ will
exemplify $\{{\underset\tilde {}\to \tau_{\lambda + \xi}}:\xi \in \kappa
\backslash g_\lambda(\alpha)\}$ is null, for which it is enough that
${\underset\tilde {}\to \tau_\alpha}$ will be (at least somewhat) like Cohen
over $\bold V[\langle {\underset\tilde {}\to \tau_{\lambda + \xi}}:\xi \in
\kappa \backslash g_\lambda(\alpha)\rangle]$, so it is reasonable to
ask that
${\underset\tilde {}\to \tau_\alpha}$ is not in the subuniverse over which
${\underset\tilde {}\to \tau_{\lambda + \xi}}$ is random when $\xi \in \kappa
\backslash g_\lambda(\alpha)$, i.e. $\alpha \notin E_\xi$.
\enddemo
\bigskip

\definition{\stag{1.4} Definition}  ${\Cal K}'_\theta$ is 
the family of $\langle P_\alpha,{\underset\tilde {}\to Q_\alpha},A_\alpha,
\mu_\alpha,{\underset\tilde {}\to \tau_\alpha}:\alpha < \beta \rangle$ such
that
\mr
\item "{$(A)$}"  $\langle P_\alpha,{\underset\tilde {}\to Q_\alpha}:\alpha <
\beta \rangle$ is a FS iteration with direct limit $P_\beta$
\sn
\item "{$(B)$}"  $A_\alpha \subseteq \alpha,\mu_\alpha = \aleph_0$, and either
$|A_\alpha| < \theta,Q_\alpha$ is Cohen forcing,
${\underset\tilde {}\to \tau_\alpha}$ the Cohen real or $|A_\alpha| \ge
\underset\tilde {}\to \theta,{\underset\tilde {}\to Q_\alpha}$ is Random$^{\bold V[\langle
{\underset\tilde {}\to \tau_\gamma}:\gamma \in A_\alpha \rangle]}$ 
(defined as in \cite[2.2]{Sh:592}) and
${\underset\tilde {}\to \tau_\alpha}$ the random real.  This is a particular
case of \cite[2.2]{Sh:592} except replacing $\kappa$ there by $\theta$ here
so we can use \cite{Sh:592}. 
\ermn
Note that Definition \scite{1.5}(1), Claim \scite{1.6}(2) are vacuous
when $\theta = 1$, the case we use here but they are natural for less specific cases.
Let ${\underset\tilde {}\to {\bar \tau}} = \langle
{\underset\tilde {}\to \tau_\alpha}:\alpha < \beta \rangle$, note that
more accurately we should write not ${\underset\tilde {}\to Q_\alpha}
= \text{ Random}^{\bold V[{\underset\tilde {}\to {\bar \tau}}
\restriction A_\alpha]}$ but Random$^{\bold V[{\underset\tilde {}\to
{\bar \tau}} \restriction A_\alpha],{\underset\tilde {}\to {\bar \tau}}
\restriction A_\alpha}$ or Random$^{\bold V,{\underset\tilde {}\to {\bar
\tau}} \restriction A_\alpha}$ as $Q_\alpha$ may be a proper subset of
Random$^{\bold V[{\underset\tilde {}\to {\bar \tau}} \restriction
A_\alpha]}$ depending on ${\underset\tilde {}\to {\bar \tau}}
\restriction A_\alpha$, too.  For the ``good" cases equality holds.
In shortness we shall write Random$^{{\underset\tilde {}\to {\bar \tau}}
\restriction A}$.
\enddefinition
\bigskip

\definition{\stag{1.5} Definition}  For $\bar Q \in {\Cal K}'_\theta$: \nl
1) $A \subseteq \ell g(\bar Q)$ is called $\bar Q$-closed if

$$
(\alpha < \ell g(\bar Q) \and |A_\alpha| < \theta \and \alpha \in A)
\Rightarrow A_\alpha \subseteq A.
$$
\mn
2) Let \footnote{this suffices when we iterate just partial random forcing
(when $|A_\beta| \ge \kappa$) and Cohen forcing (when $|Q_\beta| < \kappa$),
which is enough here.  For a more complicated situation, see \S3.}

$$
\align
\text{PAUT}(\bar Q) = \{f:&f \text{ is a one to one partial function from }
\ell g(\bar Q) \text{ to} \\
  &\ell g(\bar Q) \text{ such that dom}(f) \text{ and rang}(f) \text{ are }
\bar Q \text{-closed and for} \\
  &\beta \in \text{ dom}(f) \text{ we have } |A_\beta| < \theta
\Leftrightarrow |A_{f(\beta)}| < \theta \text{ and if} \\
  &\beta,\gamma \in \text{ dom}(f) \text{ then } \beta \in A_\gamma
\Leftrightarrow f(\beta) \in A_{f(\gamma)}\}.
\endalign
$$
\mn
3) We define the following by induction on $\ell g(\bar Q)$.  For $f \in 
\text{ PAUT}(\bar Q)$ let $\hat f$ be the partial function from
$P'_{\text{dom}(f)}$ to $P'_{\text{rang}(f)}$, 
(see \cite[Definition 2.2(3)]{Sh:592}) defined by:

$$
p_1 = \hat f(p_0) \text{ when } p_0 \in P'_{\text{dom}(f)},\text{dom}(p_1) =
\{f(\beta):\beta \in \text{ dom}(p_0)\}
$$

$$
p_1(f(\beta)) = {\Cal B}_\beta(\ldots, \text{truth value}(\zeta_\gamma \in
{\underset\tilde {}\to \tau_{f(\gamma)}}),\ldots)_{\gamma \in w^{p_0}_\beta}
$$ 

$$
\text{when }p_0(\beta)= {\Cal B}_\beta(\ldots \text{ truth value}(\zeta_\gamma
\in {\underset\tilde {}\to \tau_\gamma}) \ldots)_{\gamma \in w^{p_0}
_\beta};
$$
\mn
(so $w^{p_1}_{f(\beta)} = \{f(\gamma):\gamma \in w^{p_0}_\beta\}$). \nl
For such $\bar Q,f$ for any $P'_{\text{dom}(f)}$-name
${\underset\tilde {}\to \tau}$ we define $\hat f(\underset\tilde {}\to \tau)$,
a $P'_{\text{rang}(f)}$-name, naturally. \nl
4) For $\bar Q \in {\Cal K}'_\theta$ let

$$
\text{EAUT}(\bar Q) = \{f \in \text{PAUT}(\bar Q):\text{ dom}(f) =
\text{ rang}(f) = \ell g(\bar Q)\}.
$$
\enddefinition
\bigskip

\demo{\stag{1.6} Fact}  Let $\bar Q \in {\Cal K}'_\theta$. \nl
1) If $\alpha < \ell g(\bar Q),\bar Q' = \bar Q \restriction \alpha$ then
$\bar Q' \in {\Cal K}'_\theta$ and PAUT$(\bar Q') \subseteq \text{ PAUT}
(\bar Q)$. \nl
2) If $\alpha \le \ell g(\bar Q)$, \ub{then} $\{\beta:\beta < \alpha\}$ is
$\bar Q$-closed; and the family of $\bar Q$-closed sets is closed under intersection
and union (of any family). \nl
3) If $f \in \text{ PAUT}(\bar Q)$ and $A \subseteq \ell g(\bar Q)$ is
$\bar Q$-closed, \ub{then} $f \restriction (A \cap \text{ dom}(f))$ belongs to
$\text{PAUT}(\bar Q)$. \nl 
4) If \footnote{does $f \in \text{ PAUT}(\bar Q)$ imply $\hat f$ is an
isomorphism from $P'_{\text{dom}(f)}$ onto $P'_{\text{rang}(f)}$?  The problem
is that the order is inherited from $P_\alpha$ so is not necessarily the
same.} $f \in \text{ EAUT}(\bar Q)$, \ub{then} $\hat f$ is an automorphism
of $P'_{\ell g(\bar Q)}$. \nl 
5) If $A$ is $\bar Q$-closed and $\otimes_{\bar Q,A}$ below holds, \ub{then}
$P'_A \lessdot P_{\ell g(\bar Q)}$.  Moreover, if $q \in P'_{\ell g(\bar Q)}$
\ub{then}:
\mr
\item "{$(a)$}"  $q \restriction A \in P'_A$
\sn
\item "{$(b)$}"  $P'_{\ell g(\bar Q)} \models ``q \restriction A \le q"$
\sn
\item "{$(c)$}"  if $q \restriction A \le p \in P'_A$ \ub{then} $q' = p \cup
(q \restriction (\ell g(\bar Q) \backslash A))$ belongs to 
$P'_{\ell g(\bar Q)}$ and is a lub of $p,q$
\ermn
where
\mr
\item "{$\otimes_{\bar Q,A}$}"  if $\alpha \in A,|A_\alpha| \ge
\theta$ and $B \subseteq \alpha$ is countable then for some \nl
$f \in \text{ EAUT}(\bar Q \restriction \alpha)$ we have:
{\roster
\itemitem{ $(i)$ }  $f \restriction (B \cap A) = \text{ the identity}$
\sn
\itemitem{ $(ii)$ }  $f''(B) \subseteq A$
\sn
\itemitem{ $(iii)$ }  $f''(B \cap A_\alpha) \subseteq A \cap A_\alpha$.
\endroster}
\endroster
\enddemo
\bigskip

\demo{Proof}  Straight.  For part (5) we prove by induction on $\beta \le
\ell(\bar Q)$ that replacing $P'_{\ell g(\bar Q)},A$ by $P'_\beta,
A' = A \cap \beta$ so $q \in P'_\beta$, clauses
(a), (b), (c) in (5) hold.

In successor stages $\beta = \alpha +1$, if $\alpha \notin A$ or $|A_\alpha|
< \theta$ it is trivial, so assume $\alpha \in A \and |A_\alpha| \ge \theta$.
So $P'_{A \cap \alpha} \lessdot P_\alpha$ etc., and it is enough to show:
\mr
\item "{$(*)$}"  if in $\bold V^{P'_{A \cap \alpha}},{\Cal I}$ is a maximal
antichain in Random$^{\bold V^{P'_{A \cap \alpha \cap A_\alpha}}}$
\ub{then} in $\bold V^{P'_\alpha}$, the set ${\Cal I}$ is a maximal antichain
of Random$^{\bold V^{P'_{A_\alpha}}}$.
\ermn
By the c.c.c. this is equivalent to
\mr
\item "{$(*)'$}"  if $\zeta^* < \omega_1,\{p_\zeta:\zeta < \zeta^*\}
\subseteq P'_{A \cap (\alpha +1)},p \in P'_{A \cap \alpha}$ and $p
\Vdash_{P'_{A \cap \alpha}} ``\{p_\zeta(\alpha):\zeta < \zeta^*$ and
$p_\zeta \restriction \alpha \in G_{P'_{A \cap \alpha}}\}$ is a predense
subset of Random$^{{\underset\tilde {}\to {\bar \tau}} \restriction 
A \cap \alpha \cap A_\alpha}\,"$ then
$p \Vdash_{P'_\alpha} ``\{p_\zeta(\alpha):\zeta < \zeta^*$ and $p_\zeta
\restriction \alpha \in G_{P'_\alpha}\}$ is a predense subset of
Random$^{\bold V^{P'_{A_\alpha}}}"$.
\ermn
Assume $(*)'$ fails, so we can find $q$ such that

$$
p \le q \in P'_\alpha
$$

$$
q \Vdash_{P'_\alpha} ``\{p_\zeta(\alpha):\zeta < \zeta^* \text{ and } p_\zeta \restriction
\alpha \in G_{P'_\alpha}\} 
\text{ is not a predense subset of Random}^{{\underset\tilde {}\to
{\bar \tau}} \restriction A_\alpha}".
$$
\mn
So for some $G_{P'_\alpha}$-name $\underset\tilde {}\to r$ 

$$
q \Vdash_{P'_\alpha} ``\underset\tilde {}\to r \in 
\text{ Random}^{{\underset\tilde {}\to {\bar \tau}} \restriction
A_\alpha} (= {\underset\tilde {}\to Q_\alpha})
\text{ and is incompatible with every } p_\zeta(\alpha) \in
{\underset\tilde {}\to Q_\alpha}".
$$
\mn
Possibly increasing $q$, \wilog \,\, $\underset\tilde {}\to r$ is a ${\Cal B}
(\ldots$, truth value$(\zeta_\gamma \in {\underset\tilde {}\to \tau_\gamma}),
\ldots))_{\gamma \in w}$ where $w \subseteq A_\alpha$ is countable.  Let us define (supp
for support is from \cite[Def.2.2(1)(F)]{Sh:592}, i.e. above
supp$(\underset\tilde {}\to r) = w$)

$$
\align
B = w \cup \text{ dom}(q) &\cup \dbcu_{\zeta < \zeta^*} \text{ dom}(p_\zeta
\restriction \alpha) \cup \bigcup\{\text{supp}(q(\beta)):\beta \in
\text{ dom}(q)\} \\
  &\cup \bigcup\{\text{supp}(p_\zeta(\beta)):\beta \in \text{ dom}(p_\zeta
\restriction \alpha) \text{ and } \zeta < \zeta^*\}.
\endalign
$$
\mn
Clearly $B$ is a countable subset of $\alpha$.  By $\otimes_{\bar Q,A}$
there is $f \in \text{ EAUT}(\bar Q \restriction \alpha)$ such that

$$
f \restriction (B \cap A) = \text{ the identity}
$$

$$
f''(B) \subseteq A
$$

$$
f''(B \cap A_\alpha) \subseteq A_\alpha.
$$
\mn
As $\hat f$ is an automorphism of $P'_\alpha$ and is the identity on
$P'_{A \cap B}$, clearly

$$
\hat f(p) = p,
$$

$$
\hat f(p_\zeta) = p_\zeta,
$$

$$
p \le \hat f(q) \in P'_\alpha,
$$

$$
\hat f(\underset\tilde {}\to r) \text{ is } {\Cal B}(\dotsc,
{\underset\tilde {}\to \tau_{f(\gamma)}},\ldots)_{\gamma \in w}, \text{ and }
f''(w) \subseteq f''(B \cap A_\alpha) \subseteq A \cap A_\alpha
$$

$$
\text{ hence} \Vdash_{P'_\alpha} f(\underset\tilde {}\to r) \in
\text{ Random}^{{\underset\tilde {}\to {\bar \tau}} \restriction (A \cap A_\alpha)}
$$
\mn
and \nl
$\hat f(q) \Vdash_{P'_A}$ ``in ${\underset\tilde {}\to Q_\alpha}$, the
conditions $f(\underset\tilde {}\to r)$ and $p_\zeta(\alpha)$ are incompatible
for $\zeta < \zeta^*$", getting a contradiction.
\sn
Note that we use: in $\bold V^{P'_\alpha}$, two conditions in 
Random$^{{\underset\tilde {}\to {\bar \tau}} \restriction B}$ 
are compatible in Random$^{{\underset\tilde {}\to {\bar \tau}}
\restriction B}$ iff they
are compatible in Random$^{\bold V^{P'_\alpha}}$, for every $B \subseteq \alpha$;
in particular for $B = A_\alpha \cap A$ and $B = f''(A_\alpha \cap A)$.
\hfill$\square_{\scite{1.6}}$\margincite{1.6}
\enddemo
\bigskip

\remark{\stag{1.7} Remark}  1) Instead automorphisms $(f \in \text{ EAUT}
(\bar Q \restriction \alpha))$ we can use \nl
$f \in \text{ PAUT}(\bar Q
\restriction \alpha)$ but having more explicit assumptions and more to carry
by induction, see \S3. \nl
2) The use of Random is not essential for the last claim, we just need enough
absoluteness.
\endremark
\bigskip

\demo{Proof of 1.1}  We shall define

$$
\bar Q^\lambda = \langle P_\alpha,{\underset\tilde {}\to Q_\beta},A_\beta,\mu_\beta,
{\underset\tilde {}\to \tau_\beta}:\alpha \le \lambda + \kappa \text{ and }
\beta < \lambda + \kappa \rangle
$$
\mn
as an iteration from ${\Cal K}'_\theta$.  More accurately, $P_\alpha =
P^\lambda_\alpha = P_{\lambda,\alpha},{\underset\tilde {}\to Q_\alpha}
= {\underset\tilde {}\to Q^\lambda_\alpha} = {\underset\tilde {}\to Q_{\lambda,\alpha}},
A_\alpha = A^\lambda_\alpha = A_{\lambda,\alpha},\mu = \mu_\alpha = 
\mu^\lambda_\alpha,{\underset\tilde {}\to \tau_\alpha} = 
{\underset\tilde {}\to \tau^\lambda_\alpha} = 
{\underset\tilde {}\to \tau_{\lambda,\alpha}}$
are such that
\mn
\block   for $\beta < \lambda$ we have: $A_\beta = \emptyset$ and
${\underset\tilde {}\to Q_\beta}$ is as in \cite[Def.2.5]{Sh:592}, 
i.e. the Cohen forcing notion (and 
${\underset\tilde {}\to \tau_\beta}$ is the Cohen real),
\endblock
\sn
for $\beta \in [\lambda,\lambda + \kappa)$ we let:

$$
A_\beta = E^\lambda_{\beta - \lambda} \cup [\lambda,\beta)
$$

$$
{\underset\tilde {}\to Q_\beta} = 
\text{ Random}^{{\underset\tilde {}\to {\bar \tau}} \restriction A_\beta}
$$

$$
\text{(and } {\underset\tilde {}\to \tau_\beta} \text{ is the (partial)
random real)}
$$
\mn
(as in clause (F) of Definition 2.2 in \cite{Sh:592}. \nl
Let $P = P^\lambda =
P_{\lambda + \kappa}$. We define $P'_\alpha (= P'_{\lambda,\alpha})$
as in \cite[Definition 2.2(3)]{Sh:592}, i.e.

$$
\align
\bigl\{ p \in P_\alpha:&\text{for each } \beta \in \text{ dom}(p)
\text{ the condition } p(\beta) \text{ is either a Cohen} \\
  &\text{condition or has the form } {\Cal B}(\ldots, \text{ truth value}
(\zeta_\gamma \in {\underset\tilde {}\to \tau_\gamma}),\ldots)_{\gamma \in
w^p_\beta}, \\
  &\text{where } w^p_\beta \subseteq A_\beta \text{ is a countable set},
{\Cal B} \text{ is a Borel function with} \\
  &\text{domain and range of the right form (and } 
{\Cal B},w^p_\beta \text{ are not} \\
  &P_\beta \text{-names but actual objects)} \bigr\}.
\endalign
$$
\mn
More generally, $P'_{\lambda,A}$ for $A \subseteq \lambda + \kappa$
denote $P'_A$ for $\bar a^\lambda$ (this to help when we deal with
more than one $\lambda$).
\enddemo
\bigskip

\definition{\stag{1.7A} Definition/Fact}  1) We define ${\Cal E}^\lambda_\xi$
(for $\xi < \kappa,\lambda = \lambda^{\aleph_0} \ge 2^\kappa$), an equivalence
relation on $\lambda$ by: $\alpha {\Cal E}^\lambda_\xi \beta$ iff $g_\lambda
(\alpha) \cap \xi = g_\lambda(\beta) \cap \xi$. \nl
2) ${\Cal E}^\lambda_\xi$ (for $\xi < \kappa,\lambda = \lambda^{\aleph_0} \ge
2^\kappa$) is an equivalence relation on $\lambda$ with $\le 2^{|\xi|}$
equivalence classes, each of cardinality $\lambda$. \nl
3) If $\xi \le \zeta < \kappa$ then ${\Cal E}^\lambda_\zeta$ refines
${\Cal E}^\lambda_\xi$ and $|E^\lambda_\zeta \cap (\alpha/{\Cal E}^\lambda_\xi)| =
\lambda$ for every $\alpha < \lambda$.
\enddefinition
\bigskip

\demo{\stag{1.8} Fact}  1) If $\xi < \kappa$ then (in $\bar Q
^\lambda$)
\mr
\item "{$(a)$}"  $P'_{A_{\lambda + \xi}} = P'_{(\lambda \cap
A_{\lambda + \xi}) \cup [\lambda,\lambda + \xi)} = 
P'_{E^\lambda_\xi \cup [\lambda,\lambda + \xi)} \lessdot P'_{\lambda + \xi}$.
\ermn
Moreover
\mr
\item "{$(b)$}"  if $q \in P'_{\lambda + \xi}$ and $q \restriction
(E^\lambda_\xi \cup [\lambda,\lambda + \xi)) \le p \in P'_{E^\lambda_\xi
\cup [\lambda,\lambda + \xi)}$ then
$$
p \cup (q \restriction (\lambda \backslash E^\lambda_\xi) \cup
[\lambda,\lambda + \xi)) \in P'_{\lambda + \xi}
$$
is least upper bound of $p,q$.
\ermn
2) If $\lambda \le \lambda'$ then $P'_{\lambda',\lambda \cup [\lambda',
\lambda' + \kappa)} \lessdot P'_{\lambda',\lambda' + \kappa}$ and
$P'_{\lambda',\lambda \cup [\lambda',\lambda' + \kappa)}$ is isomorphic to
$P'_{\lambda,\lambda + \kappa}$ say by $h:P'_{\lambda,\lambda + \kappa}
\rightarrow P'_{\lambda',\lambda \cup [\lambda',\lambda' + \kappa)}$
where $h = h^{\lambda,\lambda'}$ is the canonical mapping, i.e. let $h:\lambda + \kappa
\rightarrow \lambda' + \kappa$ be $\alpha < \lambda \Rightarrow h(\alpha) =
\alpha,\xi < \kappa \Rightarrow h(\lambda + \xi) = \lambda' + \xi$, now if
$h(p) = p'$ then dom$(h(p)) = \{h(\alpha):\alpha \in \text{
dom}(p)\}$, etc. \nl
3) If $\alpha < \lambda + \kappa,|A_\alpha| \ge \theta$ (equivalently,
$\alpha \in [\lambda,\lambda + \kappa)$) then ${\underset\tilde {}\to
Q_{\lambda,\alpha}} = \bold V^{{\underset\tilde {}\to {\bar \tau}}
\restriction A_\alpha} = \bold V^{P'_{A_\alpha}}$.
\enddemo
\bigskip

\demo{Proof}  1) By \scite{1.6}(5) (really this is a particularly simple case).
I.e. let $A = A_{\lambda + \xi} = E^\lambda_\xi \cup 
[\lambda,\lambda + \xi)$,
so by \scite{1.6}(5) it is enough to check $\otimes_{Q \restriction (\lambda
+ \xi) \backslash A}$ there.  So let $\alpha \in A$ be such that $|A_\alpha| \ge \theta$
and $B \subseteq \alpha$ is countable, and we should find $f$ as there.
As $|A_\alpha| \ge \theta$ necessarily $\alpha \ge \lambda$ so for some
$\zeta < \xi,\alpha = \lambda + \zeta$.  By \scite{1.7A}(3) the existence of
$f$ is immediate. \nl
2) Straight.  \nl
3) By \cite[2.3(7)]{Sh:592} and clause (a) of part (1).
 \hfill$\square_{\scite{1.8}}$\margincite{1.8}
\enddemo
\bigskip

\demo{\stag{1.9} Fact}  $X^* = \{{\underset\tilde {}\to \tau_{\lambda + \xi}}:
\xi < \kappa\}$ is not null (in $\bold V^{P_{\lambda + \kappa}}$).
\enddemo
\bigskip

\demo{Proof}  It if is null, it is included in a Borel null set $X$ which is
coded by a real $s$ which is determined by $\langle$truth value$(p_\ell
\in G_{P_{\lambda + \kappa}}):\ell < \omega \rangle$ for some $\langle p_\ell:
\ell < \omega \rangle$ where $p_\ell \in P'_{\lambda + \kappa}$.  Let
$w = \dbcu_{\ell < \omega} \text{ dom}(p_\ell) \cup \bigcup\{\text{supp}
(p_\ell(\alpha)):\text{for some } \ell < \omega,\alpha \in \text{dom}
(p_\ell)\} \in [\lambda + \kappa]^{\le \aleph_0}$, so $w \cap \lambda \in
[\lambda]^{\aleph_0}$. \nl
Now
\mr
\item "{$\otimes$}"  there is $\xi < \kappa$ such that ($\xi > 0$ and)
{\roster
\itemitem{ $(i)$ }  $(\forall \gamma)(\lambda + \gamma \in w \Rightarrow
\gamma < \xi)$
\sn
\itemitem{ $(ii)$ }  $\alpha \in w \cap \lambda \Rightarrow \alpha \in 
E^\lambda_\xi (\subseteq A_{\lambda + \xi})$.
\endroster}
\ermn
This is possible as we have $\aleph_0$ demands; for each one the set of
$\xi < \kappa$ satisfying it is in $D$.  Now ${\underset\tilde {}\to
\tau_{\lambda + \xi}}$ is random over $\bold V^{P'_{A_{\lambda + \xi}}}$, and
$X$ is (definable) in $\bold V^{P'_{A_{\lambda + \xi}}}$ (see
\cite[2.3]{Sh:592}).  
Hence ${\underset\tilde {}\to \tau_{\lambda + \xi}}$ is not in
the Borel set $X$, a contradiction.  (Alternatively follows the proof
of \cite[2.3(2)]{Sh:592}).  \hfill$\square_{\scite{1.9}}$\margincite{1.9}
\enddemo
\bigskip

\demo{\stag{1.10} Fact}  If $D$ is a $\kappa$-complete ultrafilter \ub{then} in 
$\bold V^P,D$ is a $\kappa$-complete \nl
$\aleph_1$-saturated filter.
\enddemo
\bigskip

\demo{Proof}  Well known (see \cite{J}), as $P$ is a c.c.c. forcing notion.
\enddemo
\bn
So the ``only" point left:
\demo{\stag{1.11} Fact} If $Y \in D$ \ub{then} $\Vdash_P 
``\{{\underset\tilde {}\to \tau_{\lambda + \xi}}:\xi \in \kappa \backslash
Y\}$ is null".
\enddemo
\bn
This follows from
\demo{\stag{1.12} Main Fact}  If $\xi < \kappa,\alpha < \lambda,\alpha \notin
E_\xi$ then $\Vdash_{P'_{\lambda,\lambda + \kappa}} 
``{\underset\tilde {}\to \tau_{\lambda + \xi}} \in
N({\underset\tilde {}\to \tau_\alpha})"$, where 
$N({\underset\tilde {}\to \tau_\alpha})$ is as in
\cite[Definition 2.4]{Sh:592}.
\enddemo
\bigskip

\demo{Proof}  Note that
\mr
\item "{$(*)$}"  for $\xi < \kappa$
$$
\align
\text{EAUT}(\bar Q \restriction (\lambda + \xi)) \supseteq \{f:&f
\text{ a permutation of } \lambda + \xi, \\
  &f \restriction [\lambda,\lambda + \xi) \text{ is identity and } \\
  &f \restriction \lambda \text{ maps } E^\lambda_\gamma \text{ onto }
E^\lambda_\gamma \text{ for } \gamma < \xi\}.
\endalign
$$
\ermn
Assume toward a contradiction that the desired conclusion fails, hence for some
$p \in P_{\lambda + \kappa}$ we have $p \Vdash 
``{\underset\tilde {}\to \tau_{\lambda + \xi}} \notin N(
{\underset\tilde {}\to \tau_\alpha})"$.  In fact, possibly increasing $p$,
without loss of generality for some $\ell^*$

$$
p \Vdash ``\dsize \bigwedge_{\ell \ge \ell^*}
{\underset\tilde {}\to \tau_{\lambda + \xi}} \restriction 
{\underset\tilde {}\to n^\alpha_\ell} \in
{\underset\tilde {}\to a^\alpha_\ell}"
$$
\mn
where $\langle {\underset\tilde {}\to n^\alpha_\ell},
{\underset\tilde {}\to a^\alpha_\ell}:\ell < \omega \rangle$ is as in
\cite[Definition 2.4]{Sh:592}.
\nl
We may assume that $p \in P_{\lambda + \xi + 1}$.  Let $G \subseteq
P_{\lambda + \xi}$ be generic over $\bold V$ such that $p \restriction
(\lambda + \xi) \in G$.  So in $\bold V[G]$:

$$
\align
p(\lambda + \xi) \Vdash_{Q_{\lambda + \xi}}
``&{\underset\tilde {}\to \tau_{\lambda + \xi}} \in
{\underset\tilde {}\to T^\alpha_{\ell^*}}" \text{ where }
{\underset\tilde {}\to T^\alpha_{\ell^*}} = T_{\ell^*}
({\underset\tilde {}\to {\bar a}^\alpha}) \\
  &=: \{\eta \in {}^{\omega >} 2:\text{ if } \ell \in [\ell^*,\omega)
\text{ and } {\underset\tilde {}\to n^\alpha_\ell} \le \ell g(\eta) \\
  &\text{then } \eta \restriction {\underset\tilde {}\to n^\alpha_\ell} \in
{\underset\tilde {}\to a^\alpha_\ell}\}".
\endalign
$$
\mn
Now for some $q$ we have:

$$
p \restriction (\lambda + \xi) \le q \in G \text{ and}
$$

$$
q \Vdash ``p(\lambda + \xi) = {\Cal B}(\ldots,\text{truth value}
(\varepsilon_\ell \in {\underset\tilde {}\to \tau_{\beta_\ell}}),\ldots)
_{\ell < \omega}",
$$
\mn
where ${\Cal B} \in \bold V$ is a Borel function, $\beta_\ell \in
A_{\lambda + \xi}$ and $\varepsilon_\ell < \mu_{\beta_\ell}$ for $\ell <
\omega$.  So (see \cite[Definition 2.2(1)(F)$(\alpha)$]{Sh:592} 

$$
{\underset\tilde {}\to {\bold t}} =: {\Cal B}(\ldots,
\text{truth value}(\varepsilon_\ell \in {\underset\tilde {}\to 
\tau_{\beta_\ell}}),\ldots)
$$
\mn
is a perfect subtree of $({}^{\omega >} 2,
\triangleleft)$ of positive measure above every node, i.e. $\eta \in
{\underset\tilde {}\to {\bold t}} \Rightarrow 0 < \text{ Leb}(\{\nu \in
{}^\omega 2:\eta \triangleleft \nu \text{ and } n < \omega \Rightarrow \nu
\restriction n \in {\underset\tilde {}\to {\bold t}}\})$.  

So by the choice of
$p$ and $q$ we have $q \Vdash ``{\underset\tilde {}\to {\bold t}} \subseteq
{\underset\tilde {}\to T^\alpha_{\ell^*}}"$.  Without loss of generality
$q \in P'_{\lambda + \xi}$.  Let $w^* = \text{ dom}(q) \cup \bigcup\{
\text{supp}(q(\zeta)):\zeta \in \text{ dom}(q)\} \cup \{\beta_\ell:\ell < \omega\}
\cup \{\alpha\}$.  We can choose $p_\zeta \in P_{\lambda + \xi},f_\zeta,
\alpha_\zeta$ (by induction on $\zeta < (2^{\aleph_0})^+)$ such that
\mr
\item "{$(a)$}"  $f_\zeta \restriction (w^* \backslash \{\alpha\}) =$ the
identity
\sn
\item "{$(b)$}"  $f_\zeta \restriction [\lambda,\lambda + \xi +1)$ is the
identity
\sn
\item "{$(c)$}"  $f_\zeta \in \text{ EAUT}(\bar Q \restriction (\lambda +
\xi +1))$,
\sn
\item "{$(d)$}"  $\alpha_\zeta = f_\zeta(\alpha)$,
\sn
\item "{$(e)$}"  $\alpha_\zeta \notin \{\alpha_\varepsilon:\varepsilon 
< \zeta\}$,
\sn
\item "{$(f)$}"  $q_\zeta = \hat f_\zeta(q)$ and $\beta^\zeta_\ell = f_\zeta
(\beta_\ell)$ for $\ell < \omega$.
\ermn
Hence
\mr
\item "{$(A)$}"  $q_\zeta \Vdash ``{\Cal B}(\ldots,\text{truth value}
(\varepsilon_\ell \in {\underset\tilde {}\to \tau_{\beta^\zeta_\ell}}),
\ldots)_{\ell < \omega} \in {\underset\tilde {}\to Q_{\lambda + \xi}}$ is
perfect, of positive measure (above every node), and $\subseteq 
T_{\ell^*} ({\underset\tilde {}\to {\bar a}^{\alpha_\zeta}})"$ \nl
(also for the last statement we need ``$p_\zeta \Vdash$")
\sn
\item "{$(B)$}"  $\zeta_1 < \zeta_2 \Rightarrow \alpha_{\zeta_1} \ne
\alpha_{\zeta_2}$,
\sn
\item "{$(C)$}"  $\alpha_\zeta \notin A_{\lambda + \xi}$.
\ermn
So if $G_{P_{\lambda + \xi}} \subseteq P_{\lambda + \xi}$ is generic over
$\bold V$ and $E = \{\zeta < (2^{\aleph_0})^+:q_\zeta \in
G_{P_{\lambda + \xi}}\}$ is unbounded we get that in the universe $\bold V[
G_{P_{\lambda + \xi}}]$ for some unbounded subset $E$ of $(2^{\aleph_0})^+$
the set $Y = \dbca_{\zeta \in E} T_{\ell^*} \, 
({\underset\tilde {}\to {\bar a}^{\alpha_\zeta}})[G_{P_{\lambda + \xi}}]$
has an $\omega$-branch ${\underset\tilde {}\to \tau_{\lambda + \xi}}$ in
($\bold V[G_{P_{\lambda + \xi}}])^{Q_{\lambda + \xi}}$ which is not in
$\bold V[G_{P_{\lambda + \xi}}]$ (think), but this set $Y$ is a 
subtree of $({}^{\omega >}2,\triangleleft)$.  Hence \footnote{In fact,
we could have demanded $\beta^\ell_\zeta = \beta_\ell$ for $\ell <
\omega,\zeta < (2^{\aleph_0})^+$ so $\bold t$ hence $Y$ includes
$\underset\tilde {}\to t$ itself}
$Y$ contains a perfect
subtree, which was exactly our problem in proving theorem \scite{0.1} in
\cite{Sh:592}.  So we would like to
continue as in the proof of \cite[0.1]{Sh:592},
but for this we need a suitable ${\underset\tilde {}\to \Xi}$.  
\bn
\ub{First Proof}:

We start repeating the proof of \cite[3.3]{Sh:592}.  As there we choose
$\bar \varepsilon = \langle \varepsilon_\ell:\ell < \omega \rangle$ a
sequence of positive reals satisfying $\Sigma_{\ell < \omega} {\sqrt \varepsilon_\ell}
< 1/10$.  Letting $\kappa = 2^{\aleph_0}$, for each $\zeta < \kappa^+$
we choose $p_\zeta \in {\Cal I}_{\bar \varepsilon} \, (\subseteq
P'_{\lambda + \xi},{\Cal I}_\varepsilon$ defined in \cite[3.1]{Sh:592})
such that $q_\zeta \le p_\zeta$ and let $\bar \nu^\zeta = \langle
\nu^\zeta_\beta:\beta \in \text{ dom}(p_\zeta),|A_\beta| \ge \theta
\rangle$ witness $p_\zeta \in {\Cal I}_{\bar \varepsilon}$.  Replacing
$\langle p_\zeta:\zeta < \kappa^+ \rangle$ by a subsequence of the
same length, \wilog \, the bullets $(\bullet)$ (i.e. $\boxtimes$) of the
proof of \cite[3.3]{Sh:592} holds, so we have $i^*,\gamma^\zeta_i \, (i
< i^*,\zeta < \kappa^+),v_0,v_1,z,\langle \gamma_i:i \in v_0
\rangle,\bar \nu = \langle \nu_i:i < v_1 \rangle,s^*,m^*$ as there.

We proceed to define $p'_\zeta \ge p_\zeta$ as in the proof of
\cite[3.3]{Sh:592} and exactly as there it suffices to prove the
parallel of \cite[3.4]{Sh:592}.

The difficulty in adapting the proof of \cite[3.4]{Sh:592} is that it
uses a blueprint defined from $\langle p_\zeta:\zeta < \omega \rangle$, so
using the $\langle \eta_\beta:\beta < \ell g(\bar Q) \rangle$, which
does not exist here.  The following two possibilities seem natural.
\bn
\ub{Possibility A}:

Let $\chi > \lambda^{+ \omega}$ and let $R \in \bold V$ be $\{f:f$ a partial
function from $\chi$ to $\{0,1\}$ with countable domain$\}$.  The iterations
$\bar Q^{\lambda^{+k}}$ and the properties we are interested in, are
the same in $\bold V$ and in $\bold V^R$, so we can work in $\bold V^R$. 

So we can choose $\langle {\underset\tilde {}\to \eta_\alpha}:\alpha <
\chi \rangle$, forced to be pairwise distinct members of ${}^\kappa
2$.
\bn
\ub{Possibility B}:  (Definition)  The set ${\Cal T}^-$ of weak
blueprints as in \cite[Definition 2.9]{Sh:592} replacing $w^t$ and the
$\bar \eta^t$-s by $\bar \gamma^t$ so replacing clauses (a), (c), (i),
(j) by
\mr
\item "{$(a)'$}"  $\bar \gamma^t = \langle \gamma^t_{\bold n,k}:\bold
n < \bold n^t,k < \omega \rangle$,
\sn
\item "{$(c)'$}"  $\gamma^t_{\bold n,k}$ an ordinal increasing with
$\bold n$
\sn
\item "{$(i)'$}"  $\gamma^t_{\bold n_1,k_1} = \gamma^t_{\bold n_2,k_2}
\Rightarrow \bold n_1 = \bold n_2$ and for $k_1$
\sn
\item "{$(j)'$}"  for each $\bold n \le \bold n^t$ the sequence
$\langle \gamma^t_{\bold n,k}:k < \omega \rangle$ is constant or is
strictly increasing; if it is constant and $\bold n \in \text{
Dom}(h^t_0)$ then $h^t_0(n)$ is constant
\sn
\item "{$(k)'$}"  for $k_2 < \bold n_1 < \bold n_2 < \bold n^t$ one of
the following occurs as follows:
{\roster
\itemitem{ $(\alpha)$ }  $\bold n^{t^*},\bold
m^{t^*},h^{t^*}_0,h^{t^*}_1,h^{t^*}_2,\bar n^{t^*}$ is as in the proof
of \cite[3.4]{Sh:592}
\sn
\itemitem{ $(\beta)$ }  $\gamma^{t^*}_{i,\zeta}$ is $\gamma_{i,\zeta}
\in \text{ Dom}(p'_\zeta)$ above (for $i < i^* = \bold n^t,\zeta < \omega$).
\endroster}
\ermn
Nothing relevant to us changes.

The actual difference between the two possibilities for the rest of
the proof is small and we shall use the second.  We define the weak
blueprint $t^* = t(*) = (\bar \gamma^{t^*},\bold n^{t^*},\bold
m^{t^*},h^{t^*}_0,h^{t^*}_1,h^{t^*}_2,\bar n^{t^*})$. 

Clearly $t^* \in {\Cal T}^-$; we would like to proceed and choose
$\underset\tilde {}\to \Xi$ (but first choose $B_{\lambda,n} \subseteq
\lambda$ for $n < \omega$ such that: if $Y_\xi \in \{A_{\lambda,\lambda +
\xi} \cap \lambda,\lambda \backslash A_{\lambda,\lambda +
\varepsilon}\}$ for $\xi < \kappa,X_n \in \{B_{\lambda,n},\lambda
\backslash B_{\lambda,n}\}$ and $\dbca_{\xi < \kappa} Y_\xi \ne
\emptyset$ then $\dbca_{\xi < \kappa} Y_\xi \cap \dbca_{n < \omega} X_n
\subseteq \lambda$ has cardinality $\lambda$ and the $B_n$ separates
the $\gamma^{t^*}_{\bold n,k},\bold n < \bold m^{t^*},k < \omega$.
Now for $\xi < \lambda,\bold n < \bold m^{t(*)}$ let
$\Delta_{\lambda,\xi,\bold n,k} = \{\gamma < \lambda:$ for every
$\zeta < \xi$ we have $\gamma^{t(*)}_{\bold n,k} \in
A_{\lambda,\lambda + \zeta} \equiv \gamma \in A_{\lambda,\lambda +
\zeta}$ and $\gamma^{t^*}_{\bold n,k} \in B_{\lambda,n} \equiv \gamma
\in B_{\lambda,n}\}$. 

For $\xi < \lambda$ let ${\Cal T}_{\lambda,\xi} = \{t^* \in {\Cal T}^-:(\bold
d^{t(*)},\bold m^{t(*)},h^{t(*)}_0,h^{t(*)}_1,h^{t(*)}_2,\bar
n^{t(*)}) =$ \nl
$(\bold n^t,\bold m^t,h^t_0,h^t_1,h^t_2,\bar n^t),\bold n
\in [\bold m^t,\bold n^t) \and k < \omega \Rightarrow \gamma^{t(*)}_{\bold
n,k} = \gamma^t_{\bold n,k}$ and $\bold n < \bold m^t \Rightarrow
\gamma^t_{n,k} \in \Delta_{\lambda,\xi,\bold n,k}\}$.

Let $\Gamma = \{\gamma^{t(*)}_{\bold n,0},\gamma^t_{\bold n,0}
+1:\gamma^{t^*}_{\bold n,0} = \gamma^t_{\bold n,1}\} \cup \{\lambda\}
\cup \{\dbcu_{k < w} \gamma^{t^*}_{\bold n,k}:\gamma^t_{n,0} <
\gamma^{t^*}_{n,k}\}$ for any (equivalently some) $t \in
{\Cal T}_{\lambda,0}$; and let $\{\gamma^j_n:n \le n(*)\}$ list $\Gamma
\backslash \lambda$ in increasing order.

Now we choose the sequence $\langle {\underset\tilde {}\to \Xi^j}:j <
\omega \rangle$ by induction on $n \le n(*)$ similarly to
\cite[Definition 2.11]{Sh:592} but just for the $t \in {\Cal
T}_{\lambda^{+j},\gamma^j_n}$
\mr
\item "{$(a)$}"   ${\underset\tilde {}\to \Xi^j_n}$ is a
$P_{\lambda^{+j},\gamma^j_n}$-name of a full finitely additive measure
(on ${\Cal P}(w)$)
\sn
\item "{$(b)$}"  if $\bold n < \bold n^{t^*},\gamma^j_{\bold n,k} <
\gamma^j_n$ for every $k < \omega$ and $\gamma^{t^*}_{\bold n,0} <
\gamma^{t^*}_{\bold n,1}$ and $t \in {\Cal T}_{\lambda^{+j},\gamma^j_n}$
\ub{then} $\Vdash_{P_{\lambda^{+j},\gamma^j_n}}$ ``the following set has
${\underset\tilde {}\to \Xi^j_n}$-measure 1" \nl
$\{k < \omega$: if $\ell \in [n^t_k,n^t_{k+1})$ then $(h^t_0(\bold
n))(\ell) \in {\underset\tilde {}\to G_{\gamma^j_{\bold n,\ell}}}\}$
\sn
\item "{$(c)$}"  \ub{if}
{\roster
\itemitem{ $(\alpha)$ }  $t \in T_{\lambda^{+j},\gamma^j_n},\bold n <
\bold n^t,\gamma^j_{\bold n,k} < \gamma^j_n$ for $k < \omega$ and
$\gamma^j_{\bold n,0} = \gamma^j_{\bold n,1},n \in \text{ Dom}(h^t_1)$
\sn
\itemitem{ $(\beta)$ }  $\underset\tilde {}\to r,{\underset\tilde
{}\to r_\ell}$ (for $\ell < \omega$) are $P'_{\lambda^{+j},\gamma^j_{\bold
n,0}}$-names of members of ${\underset\tilde {}\to
Q_{\lambda^{+j},\gamma^j_{\bold n,0}}}$ satisfying
$(**)_{\underset\tilde {}\to r,
\langle {\underset\tilde {}\to r_\ell}:\ell < \omega \rangle}$ below,
\ub{then} \nl
$\Vdash_{P_{\lambda^{+j},\gamma^j_n}}$ ``if $\underset\tilde {}\to r
\in {\underset\tilde {}\to G_{\gamma^j_n}}$ then $1 - h^t_1(\bold n)
\le \text{ Av}_{\underset\tilde {}\to \Xi^j_n}(\langle |\{\ell \in 
[n^t_k,n^t_{k+1}):{\underset\tilde {}\to r_\ell} \in
{\underset\tilde {}\to G_{Q_{\gamma^j_n}}}\}|/(n^t_{k+1} - n^t_k):k <
\omega \rangle)$ where
\sn
\itemitem{ ${{}}$ }  $(**)^{\bar Q}_{\underset\tilde {}\to \tau,\langle
{\underset\tilde {}\to \tau_\ell}:\ell < \omega \rangle} \qquad
\underset\tilde {}\to \tau,{\underset\tilde {}\to \tau_\ell}$ are
$P'_{A_\alpha}$-names of members of \nl

\hskip60pt  $Q_\alpha,
(\langle{\underset\tilde {}\to r_\ell}:\ell < \omega \rangle \in \bold V)$
and, in $\bold V^{P_\alpha}$, for every $\tau^\ell \in Q_\alpha$ \nl

\hskip60pt satisfying $\tau \le \tau'$ we have
$$
\text{Av}_{\Xi^\tau_\alpha} (\langle a_k(\tau'):k < \omega \rangle)
\ge 1 - h^t_1(n)
$$
\nl
\hskip40pt where
\sn
\itemitem{ ${{}}$ } $\boxtimes \qquad a_k(\tau') = a_k(\tau',\bar \tau) =
a_k(\tau',\bar \tau,\bar n^t)$ \nl

$\qquad \quad \bigl( \dsize \sum_{\ell \in [n^\ell_k,n^t_{k+1})}
\frac{\text{Leb}(\text{lim}(r') \cap \text{ lim}(r_\ell))}
{\text{Leb}(\text{lim}(r'))} \bigr) \,\, \frac{1}{n^t_{k+1} - n^t_k} \bigr)$
\endroster}
\ermn
(so $a_k(r',\bar r,\bar n) \in [0,1]$ is well defined for $k <
\omega,\bar t = \langle r_\ell:\ell < \omega \rangle,\{r,r_\ell\}
\subseteq \text{ Random}, \bar n = \langle n_\ell:\ell < \omega
\rangle,n_\ell < n_{\ell +1} < \omega)$.
\sn
Let us carry the induction. 
\bn
\ub{Case 1}:  $n=0$.

Necessarily $\gamma^j_n = \lambda^{+j}$.  Note that $\alpha <
\gamma^j_n \Rightarrow A_{\lambda,\alpha} = \emptyset$ and the sets
$\{\gamma^t_{\bold n,k}:t \in {\Cal T}_{\lambda,0}\} \cap \lambda$ for $\bold
n < \bold n^{t^*},k < \omega$ are pairwise disjoint.  Hence the proof
if easy (similar to \cite[Lemma 2.14]{Sh:592}).
\bn
\ub{Case 2}:  $n+1,\gamma^j_{n+1} > \gamma^j_n +1$.

Similar to Case 1.
\bn
\ub{Case 3}: $n+1$ and $\gamma^j_{n+1} = \gamma^j_n+1$.

Let $\bold n$ be such that $\gamma^j_{\bold n,k} = \gamma^j_n$ for $k
< \omega$.  If $\bold n \in \text{ Dom}(h_0)$ nothing to do, so assume
$\bold n \in \text{ Dom}(h_1)$.  We would like to repeat the proof of
\cite[Lemma 2.16(1)]{Sh:592}, \ub{but} this proof needs, in our
notation, not just $P'_{\lambda^{+j},A_{\gamma^j_n}} \lessdot
P'_{\lambda^{+j},\gamma^j_n}$ (which we have proved) but also
``${\underset\tilde {}\to \Xi^j_n} \cap {\Cal P}
(\omega)^{P'_{\lambda^{+j},A_{\gamma^j_n}}}$ is a
$P'_{\lambda^{+j},A_{\gamma^j_n}}$-name", which is not part of our
induction hypothesis here.  But we have in our induction hypothesis 
${\underset\tilde {}\to \Xi^{j+1}_n}$.  We let $\chi$ be large enough
and choose $M_0,M_1$ elementary submodels of $({\Cal
H}(\chi),\in,<^*_\chi)$ of cardinality $\lambda^{+j}$, such that
$\lambda^{+j} \subseteq M_\ell,\{\bar Q^{\lambda^{j+1}},
{\underset\tilde {}\to \Xi^{j+1}_n},{\Cal T}^1_{\lambda^{+j+1},\gamma^j_n}\}
\subseteq M_\ell$ and $M_0 \in M_1$.

Now we can find a one-to-one function $h$ from $\lambda^{+j} + \kappa$
onto $(\lambda^{+j+1} + \kappa) \cap M_1$ such that:
\mr
\item  for $\xi < \kappa,h(\lambda^{+j} + \xi) = \lambda^{+j+1} + \xi$
\sn
\item  for $\lambda \in \lambda^{+j}$, (and $\beta \in
[\lambda^{+j+1},\gamma^j_n)$ we have $\gamma \in
A_{\lambda^{+j},\beta} \Leftrightarrow h(\gamma) \in
A_{\lambda^{+j+1},h(\beta)}$
\sn
\item   for $\lambda < \lambda^{+j},n < \omega$ we have $\gamma \in
B_{\lambda^{+j},n} \Leftrightarrow h(\gamma) \in B_{\lambda^{+j+1},n}$
\sn
\item  for $\gamma < \lambda^{+j}$ we have $\gamma \in
A_{\lambda^{+j},\gamma^j_n} 
\Leftrightarrow h(\gamma) \in M_0$.
\ermn
This should be clear.  So we have finished the induction step in
constructing the $\Xi^j_{n^*}$, so we have completed the missing point
in the proof of the parallel of \cite[3.4]{Sh:592}, hence of the
parallel to \cite[3.3]{Sh:592}, so we are done proving our main
theorem.
\bn
\ub{Alternate Proof}:

Let $\bar Q^{n^*} = \langle P_\alpha,{\underset\tilde {}\to Q_\beta},A_\beta,
{\underset\tilde {}\to \tau_\beta}:\alpha < \lambda^{+n^*} + \kappa \rangle$
be as above except that with $A^*_\beta$ being  $\emptyset$ if $\beta <
\lambda^{+n^*}$ and $\{\gamma < \lambda^{+n^*}:\beta - \lambda \notin
E^{\lambda^{+n^*}}_\gamma\} \cup [\lambda,\beta)$ if $\beta \in
[\lambda^{+n^*},\lambda^{+n^*} + \kappa) \backslash \{\gamma_n:n < n^*\}$
and $\{\gamma < \lambda^{+n^*}:\beta - \lambda \notin E^{\lambda^{+n^*}}
_\gamma\} \cup [\lambda^{+n^*},\beta)$ if $\beta = \lambda + n^* + \gamma_n$. \nl
It suffices to prove that for some we do not have $\langle
\lambda^{+n} + \gamma_\ell:\ell < n^* \rangle,
\bar \varepsilon,\langle p'_\ell:\ell < \omega \rangle$ as 
above for the $\bar Q^*$ defined above.  We take $k=n^*$. \nl
[Why?  Choose $\chi^* >> \chi$, and an elementary submodel $M$ of $({\Cal H}
(\chi^*),\in,<^*_\chi)$ of cardinality $\lambda$ to which $\bar Q^*$ belongs,
$M^\omega \subseteq M$.  Now $\bar Q^*$ as interpreted in $M$ is just
$\bar Q^\lambda$ with changes to names, and we have enough absoluteness.]
\nl
Now, fixing $\langle \gamma_\ell:\ell < n^* \rangle$ (not the others!) we
prove by induction on $n \le n^*$ that we can find a suitable $\langle
\eta_\beta:\beta < \lambda^{+n^*} + \gamma_n+1+\kappa \rangle,\left< \langle
{\underset\tilde {}\to \Xi^t_\alpha}:\alpha < \lambda^{+n^*} + \gamma_n+1
\rangle:t \in T\right>$, \nl
$P_{\lambda^{+n^*} + 
\gamma_n}$-name ${\underset\tilde {}\to \Xi}$ (stipulating $\gamma_0+1 =
\gamma_0,\gamma_{n^*} = \kappa$).  Which means: letting $P_\alpha,
{\underset\tilde {}\to Q_\beta},A_\beta,\mu_\beta,
{\underset\tilde {}\to \tau_\beta}$ be as in $\bar Q^*$ there are $\eta_\beta,
{\underset\tilde {}\to \Xi^t_\beta}$ such that $\langle P_\alpha,
{\underset\tilde {}\to Q_\beta},A_\beta,\mu_\beta,
{\underset\tilde {}\to \tau_\beta},\eta_\beta,
({\underset\tilde {}\to \Xi^t_\alpha})_{t \in {\Cal T}}:\alpha \le
\lambda^{+n^*} + \gamma_n,\beta < \lambda^{+n^*} + \gamma_n \rangle \in
{\Cal K}^3_\theta$, i.e. satisfies \cite[Definition 2.11]{Sh:592} except
that in clause (d) there we demand that if $\bold n \in \text{ dom}(h^t_1)$
and $\langle \eta^t_{\bold n,\ell}:\ell < \omega \rangle$ is constant, then
$\alpha_\ell \in \{\gamma_m:n \le m < n^*\}$.  For $n$ we use $n+1$,
L\"owenheim-Skolem argument and uniqueness in Definition 2.11,
\cite{Sh:592}.

In detail:  for $n = 0$, and just let $\underset\tilde {}\to \Xi$ be a
$P_{\lambda^{+n^*} + \gamma_0}$-name for a finitely additive measure on
$\omega$ (can be Ramsey if $\bold V \models CH$) such that conditions
2.11(e)+(f) in \cite{Sh:592} are satisfied.  For $n$, by the induction
hypothesis we have $\langle \eta^n_\beta:\beta < \lambda^{+n^*} + \gamma_n
\rangle, \left< \langle {\underset\tilde {}\to \Xi^t_\alpha}:\alpha \le
\lambda^{+n^*} + \gamma_n \rangle:t \in {\Cal T} \right>$, we can
find $A \subseteq \lambda^{+n^*} + \gamma_n$ of cardinality $\lambda^{+n}$
such that $\dbcu_{m \ge n} A^*_{\gamma_m} \subseteq A,[\lambda^{+n^*},
\lambda^{+n^*} + \gamma_n) \subseteq A,P'_A \lessdot P_{\lambda^{+n^*} +
\gamma_n}$ and $t \in {\Cal T} \Rightarrow 
{\underset\tilde {}\to \Xi^t_{\lambda^{+n^*} + \gamma_n}}
\restriction {\Cal P}(\omega)^{{\bold V}^{P'_A}}$ is a
$P'_A$-name.  If $A = \lambda^{+n^*}$, then as in \cite{Sh:592}, we can
define ${\underset\tilde {}\to \Xi^t_{\lambda + \gamma_n+1}}$ as required. 
This is not necessarily true, but some $f \in \text{EAUT}(\bar a \restriction
\lambda^{+n^*} + \gamma_n)$ maps $A$ onto $A_{\lambda^{+n^*},\lambda + \gamma_n}$, and maps
$A_{\lambda^{+n^*},\beta}$ (for $\beta \ne \lambda^{+n^*} 
+ \gamma_n$) onto itself, so possibly
changing the $\Xi^t_{\lambda^{+n^*} + \gamma_n},\eta_\beta(\beta <
\lambda^{+n^*} + \gamma_n)$ we can define ${\underset\tilde {}\to \Xi^t
_\alpha} \, (\alpha \le \lambda^{+n^*} + \gamma_n +1),\eta_\beta(\beta \le
\lambda^{+n^*} + \lambda^{+n^*})$.  The advance from $\lambda^{+n^*} + \gamma_n +1$ to
$\lambda^{+n^*} + \gamma_{n+1}$ is as for $n=0$: as we are proving $\in {\Cal K}'_\theta$
(not $\in {\Cal K}_3$) there are no problems repeating the proof from
\cite{Sh:592}.

So we have finished the proof of \scite{1.1}.
\hfill$\square_{\scite{1.12}}$\margincite{1.12}
\enddemo
\newpage

\head {\S2 Non-null set with no non-null function} \endhead  \resetall \sectno=2
\bn
Peter Komj\'ath referred me to the following problem, answered below.
\proclaim{\stag{2.1} Theorem}  It is consistent that:
\mr
\item "{$\oplus$}"  there is a non-null $A \subseteq \Bbb R$ such that:
for every $f:A \rightarrow \Bbb R$, the function $f$ as a subset of the plane
$\Bbb R \times \Bbb R$ is null
\ermn
provided that ``ZFC + there is a measurable cardinal" is consistent.
\endproclaim
\bigskip

\demo{Proof}  We can use ${}^\omega 2$ instead of $\Bbb R$.  Let $\kappa$ be
measurable, $D$ a normal ultrafilter on $\kappa$.  Let $\lambda = 
\lambda^{\aleph_0},\lambda \ge 2^\kappa$ for simplicity, and we use the same
forcing as in \S1.  Now we interpret the Cohen forcing
${\underset\tilde {}\to Q_\alpha}$ (for $\alpha < \lambda)$ as

$$
\align
\bigl\{ \langle (n_\ell,a_\ell):\ell < \ell^* \rangle:&\ell^* < \omega,n_\ell
< \omega,a_\ell \text{ is a subset of} \\
  &({}^{n_\ell}2) \times ({}^{n_\ell}2), \\
  &\eta \in {}^{n_\ell}2 \Rightarrow |\{\nu \in {}^{n_\ell}2:(\eta,\nu) 
\in a_\ell\}|/2^{n_\ell} \ge 1 -4^{-\ell} \bigr\}.
\endalign
$$
\mn
We interpret the generic sequence $\bar a^\alpha = \langle (n^\alpha_\ell,
a^\alpha_\ell):\ell < \omega \rangle$ as the following null subset of the
plane:

$$
\align
\text{null}(\bar a^\alpha) =: \{(\eta,\nu):&\eta \in {}^\omega 2
\text{ and for infinitely many} \\
  &\ell < \omega \text{ we have } (\eta \restriction n^\alpha_\ell,\nu
\restriction n^\alpha_\ell) \notin a^\alpha_\ell\}.
\endalign
$$
\mn
Let ${\underset\tilde {}\to X^*} = \{{\underset\tilde {}\to \tau_{\lambda +i}}:i < \kappa\}$, as in \S1, it is not null (in fact everywhere).  Now suppose
$p \Vdash ``\underset\tilde {}\to g:{\underset\tilde {}\to X^*} \rightarrow
{}^\omega 2"$; of course, this $\underset\tilde {}\to g$ is unrelated
to the $g_\lambda$'s from \S1.  So for each $i < \kappa,\underset\tilde {}\to g
({\underset\tilde {}\to \tau_{\lambda+i}})$ is a $P_{\lambda + \kappa}$-name
involving the conditions $p_{i,\ell} \in P'_{\lambda + \kappa}$ (for $\ell <
\omega$).  Let for $\xi < \kappa,M_\xi \prec ({\Cal H}(\beth_\omega(
\lambda^{+ \omega})^+),\in,<^*)$ be such that $\bar Q,\langle
{\underset\tilde {}\to \tau_{\lambda +i}}:i < \kappa \rangle,p,
\underset\tilde {}\to g,\left< \langle p_{i,\ell}:\ell < \omega \rangle:i <
\kappa \right>,\lambda,\xi$ belongs to $M_\xi$ and $\|M_\xi\| = 2^{\aleph_0}$
and ${}^\omega(M_\xi) \subseteq M_\xi$ and let $w_\xi = M_\xi \cap (\lambda
+ \kappa)$.  Then for some $A \in D$ we have
\mr
\item "{$(a)$}"  $\gamma \ne \xi \in A \Rightarrow \lambda + \xi \notin
w_\gamma$,
\sn
\item "{$(b)$}"  $\langle w_\gamma:\gamma \in A \rangle$ and $\langle
M_\gamma:\gamma \in A \rangle$ are $\Delta$-systems with hearts $w^*,M^*$
respectively,
\sn
\item "{$(c)$}"  $w_\gamma \cap [\lambda,\lambda + \gamma)$ is constant,
$w_\gamma \cap \{\lambda + \xi:\xi \in A\} = \{\lambda + \gamma\}$, moreover
sup$(w_\gamma) < \lambda + \text{ min}(A \backslash (\gamma +1))$ and
$w_\gamma \cap [\lambda,\lambda + \gamma) \subseteq w^*$,
\sn
\item "{$(d)$}"  for $\gamma,\xi \in A$, the models $M_\gamma,M_\xi$ are
isomorphic, with the isomorphism mapping $\gamma$ to $\xi$, and $\bar Q,
\langle {\underset\tilde {}\to \tau_{\lambda + \gamma}}:\gamma < \kappa
\rangle,\underset\tilde {}\to g$ to themselves, in fact is the identity on
$M^*$.
\ermn
So for all $\gamma \in A$ we have an isomorphic situation.

For $\gamma \in A$ let $\bar r^\gamma = \langle r_{\gamma,m}:m < \omega
\rangle$ be the $<^*$-first maximal antichain of $P_{\lambda + \kappa}$ above
$p$, of forcing conditions deciding whether $\underset\tilde {}\to g
({\underset\tilde {}\to \tau_{\lambda + \gamma}})$ is in $\bold V[\langle
{\underset\tilde {}\to \tau_\beta}:\beta \le \lambda + \gamma +1 \rangle) =
\bold V^{P_{\lambda + \gamma +1}}$, and deciding whether
$\underset\tilde {}\to g({\underset\tilde {}\to \tau_{\lambda + \gamma}})$ is
in $\bold V[\langle {\underset\tilde {}\to \tau_\beta}:\beta \in M^* \rangle
\char 94 \langle {\underset\tilde {}\to \tau_{\lambda + \gamma}} \rangle]$ and
in both cases if yes, \wilog, \, it forces for some Borel function,
${\Cal B} = {\Cal B}_{\gamma,m}$ that $r_{\gamma,m} \Vdash
``\underset\tilde {}\to g({\underset\tilde {}\to \tau_{\lambda + \gamma}}) =
{\Cal B}(\ldots,{\underset\tilde {}\to \tau_\beta},\dotsc,
{\underset\tilde {}\to \tau_{\lambda + \gamma}})_{\beta \in w}"$ where in the
first case $w = w_\gamma \cap \gamma$ and in the second case $w = w^*$.  As
we can shrink $A$ (as long as it is in $D$) \wilog \, ${\Cal B}_{\gamma,m}$ does not
depend on $\gamma$, i.e. ${\Cal B}_{\gamma,m} = {\Cal B}_m$.  Also, \wilog \,
the answer (= decision) for each $m$ of $r_{\gamma,m}$ does not depend on
$\gamma$.

Choose $\alpha < \lambda$ such that $\lambda(\alpha) = A$ where
$g_\lambda$ is from \scite{1.3} (and is not related to
$\underset\tilde {}\to g$)  and $\alpha \notin
\dbcu_{\gamma \in A} M_\gamma$.  It is enough to prove the following two
statements:
\mr
\item "{$(*)_1$}"  $\{({\underset\tilde {}\to \tau_{\lambda + \gamma}},
\underset\tilde {}\to g({\underset\tilde {}\to \tau_{\lambda + \gamma}})):
\gamma \in \kappa \backslash A\}$ is null (subset of the plane)
\sn
\item "{$(*)_2$}"  $\{({\underset\tilde {}\to \tau_{\lambda + \gamma}},
\underset\tilde {}\to g({\underset\tilde {}\to \tau_{\lambda + \gamma}})):
\gamma \in A\}$ is null.
\endroster
\enddemo
\bigskip

\demo{Proof of $(*)_1$}  Trivial as 
$\{{\underset\tilde {}\to \tau_{\lambda + \gamma}}:\gamma \in \kappa
\backslash A\}$ is null (as proved in \scite{1.11}).
\enddemo
\bigskip

\demo{Proof of $(*)_2$}  we shall use the proof of Fact \scite{1.12}
and the choice of $\alpha$ to show that

$$
\{({\underset\tilde {}\to \tau_{\lambda + \gamma}},
\underset\tilde {}\to g({\underset\tilde {}\to \tau_{\lambda + \gamma}})):
\gamma \in A\} \subseteq \dbcu_\ell ({}^\omega 2 \times {}^\omega 2 \backslash \text{ lim
tree}_\ell({\underset\tilde {}\to {\bar a}^\alpha}))
$$
\mn
where here tree$_{\ell^*}({\underset\tilde {}\to {\bar a}^\alpha}) =
\{(\eta,\nu) \in {}^\omega 2 \times {}^\omega 2$
if $\ell \in [\ell^*,\omega)$ then 
$(\eta \restriction n_\ell,\nu \restriction n_\ell) \notin
{\underset\tilde {}\to a^\alpha_\ell}\}$. \nl
Let ${\underset\tilde {}\to A_m} = \{\xi:r_{\xi,m} \in
{\underset\tilde {}\to G_{P_{\lambda + \kappa}}}\}$ for $\xi < \kappa$.

It suffices to prove, for each $m$, that
\mr
\item "{$(*)_{2,m}$}"  
$p \Vdash ``\{({\underset\tilde {}\to \tau_{\lambda + \gamma}},
\underset\tilde {}\to g({\underset\tilde {}\to \tau_{\lambda + \gamma}})):
\gamma \in {\underset\tilde {}\to A_m}\}$ is null".
\ermn
We do it by cases (note that in each case ``for some $\gamma \in
{\underset\tilde {}\to A_m}$" is equivalent 
to ``for every $\gamma \in {\underset\tilde {}\to A_m}$").
\enddemo
\bn
\ub{Case 1}:  $r_{\gamma,m} \Vdash ``\underset\tilde {}\to g
({\underset\tilde {}\to \tau_{\lambda + \gamma}})$ is in $\bold V
[\langle {\underset\tilde {}\to \tau_\beta}:\beta \in M^* \rangle \char 94
\langle \tau_{\lambda + \gamma} \rangle]"$.

Clearly $(*)_{2,m}$ holds as the graph of a Borel function is null and
apply this to the function $\rho \mapsto {\Cal
B}(\dotsc,{\underset\tilde {}\to
\tau_\beta},\ldots,\ldots;\rho)_{\beta \in M^*}$.
\bn
\ub{Case 2}:  $r_{\gamma,m} \Vdash 
``({\underset\tilde {}\to \tau_{\lambda + \gamma}},\underset\tilde {}\to g
({\underset\tilde {}\to \tau_{\lambda + \gamma}})) \in ({}^\omega 2
\times {}^\omega 2 \backslash \dbcu_{\ell < \omega} \text{ lim tree}_\ell
({\underset\tilde {}\to {\bar a}^\alpha}))"$.
\nl
As the set on the right side is null this is trivial.
\bn
\ub{Case 3}:  Not Case 2 and
$r_{\gamma,m} \Vdash ``\underset\tilde {}\to g
({\underset\tilde {}\to \tau_{\lambda + \gamma}})$ is not in 
$\bold V^{P_{\lambda + \gamma +1}}"$; remember $\alpha < \lambda$ was chosen
such that $g_\lambda(\alpha) = A$.

Fix $\xi = \text{ min}(A)$, and let $\xi_1 = \text{ min}(A \backslash
(\xi +1))$.  Let $p^*,\ell^*$ be such that:
$r_{\xi,m} \le p^* \in P_{\lambda + \xi_1}$ and $\ell^* < \omega$ and $p^*
\Vdash ({\underset\tilde {}\to \tau_{\lambda + \xi}},
\underset\tilde {}\to g({\underset\tilde {}\to \tau_{\lambda + \gamma}})) \in
\text{ lim tree}_{\ell^*}({\underset\tilde {}\to {\bar a}^\alpha})"$ (note:
both are $P_{\lambda + \xi_1}$-names and $r_{\xi,m} \in 
P_{\lambda + \xi_1})$.

As in the proof of \scite{1.12} we can for $\zeta < (2^{\aleph_0})^+$ choose
by induction $f_\zeta,a_\zeta,p_\zeta$ such that

$$
f_\zeta \in \text{ EAUT}(\bar Q \restriction (\lambda + \xi_1)),\alpha_\zeta =
f_\zeta(\alpha),\alpha_\zeta \notin \{\alpha_{\zeta_1}:\zeta_1 < \zeta\}
$$

$$
\align
\text{ and } p_\zeta = &\hat f_\zeta(p^*) \text{ and } f_\zeta \restriction
(\dbcu_{\gamma \in A} M_\gamma \cap (\lambda + \xi_1)), \\
  &f_\zeta \restriction [\lambda,\lambda + \xi_1) \text{ are the identity}.
\endalign
$$
\mn
So $\hat f_\zeta$ maps ${\underset\tilde {}\to \tau_{\lambda + \xi}},\
\underset\tilde {}\to g({\underset\tilde {}\to \tau_{\lambda + \xi}})$ to
themselves.  Let $G_{P_{\lambda + \xi +1}} \subseteq P_{\lambda + \xi +1}$ be
generic over $\bold V$ such that $E = \{\zeta:p_\zeta \restriction (\lambda +
\xi + 1) \in G_{P_{\lambda + \xi +1}}\}$ is unbounded in $(2^{\aleph_0})^+$.
So in $\bold V[G_{P_{\lambda + \xi +1}}]$ we have for $\zeta \in E$, a
$(P_{\lambda + \xi_1}/G_{P_{\lambda + \xi+1}})$-name of a real
$\underset\tilde {}\to g({\underset\tilde {}\to \tau_{\lambda + \gamma}})$,
which is not in $\bold V^{P_{\lambda + \xi +1}}$ and $p_\zeta \Vdash
``({\underset\tilde {}\to \tau_{\lambda + \xi}},
\underset\tilde {}\to g({\underset\tilde {}\to \tau_{\lambda + \xi}})) \in
\text{ lim tree}_{\ell^*}(\bar a^{\alpha_\zeta})"$.

Now, still as in the proof of \scite{1.12}, we can prove by induction on
$\beta \in [\xi +1,\xi_1]$ that in $\bold V^{P_\beta}$, no $(\eta_0,\eta_1)
\in ({}^\omega 2 \times {}^\omega 2)^{V^{P_{\lambda + \beta}}} \backslash ({}^\omega 2
\times {}^\omega 2)^{{\bold V}^{P_{\lambda + \xi +1}}}$ 
belongs to lim tree$_{\ell^*}(\bar a^{\alpha_\zeta})$ 
for unboundedly many $\zeta \in E$.  The induction
is straight;  for successor case
$\beta +1$ using $\beta \notin A = g_\lambda(\alpha_\zeta)$ (and $\hat f$ for
some $f \in \text{ EAUT}(\bar Q \restriction (\lambda +i)))$.  Contradiction
as $p^*$ forces that $({\underset\tilde {}\to \tau_{\lambda + \xi}},
\underset\tilde {}\to g({\underset\tilde {}\to \tau_{\lambda + \xi}}))$ is not
in $({}^\omega 2 \times {}^\omega 2)^{{\bold V}^{P_{\lambda + \xi +1}}}$ as the second 
coordinate is not in $\bold V^{P_{\lambda + \xi +1}}$.
\bn
\ub{Case 4}:  Neither Case 1 nor Case 2, nor Case 3.

As in the Case 3, but we can have $f_\zeta \restriction ((\lambda + \xi_1)
\cap M^*)$ is the identity.  We let hence $\hat f_\zeta \restriction
(P_{\lambda + \xi_1} \cap M^*)$ is the identity.  We let $B_1 = \{\alpha:\alpha <
\lambda + \xi_1 \text{ and } \alpha \notin M_\xi \cap (\lambda + \xi)
\backslash M^*\}$, using $\hat f$ for $f \in \text{ EAUT}(\bar Q)$ we can
easily show that $P'_{B_1} \lessdot P'_{\lambda + \kappa}$, and
$G_{B_1} = G_{P_{\lambda + \xi+1}} \cap P'_{B_1}$ (or see \S3).  Now we have
in $\bold V[G_{P_{\lambda + \xi+1}}]$:

$$
\underset\tilde {}\to E[G_{B_1}] = \{\zeta:p_\zeta \restriction (\lambda +
\xi +1) \in G_{B_1}\} \text{ is unbounded in } (2^{\aleph_0})^+.
$$
\mn
We can assume that $p_\zeta \Vdash ``\underset\tilde {}\to g
({\underset\tilde {}\to \tau_{\lambda + \xi}}) \notin \bold V
[{\underset\tilde {}\to G_{B_1}}]",p_\zeta \in {\Cal I}_{\bar \varepsilon}$
for some $\bar \varepsilon$ (as in \S1), as if not \wilog \, $p_\zeta$ forces
$\underset\tilde {}\to g({\underset\tilde {}\to \tau_{\lambda + \xi}}) =
{\underset\tilde {}\to \tau'}$ such that
${\underset\tilde {}\to \tau'} = {\underset\tilde {}\to {\Cal
B}}(\ldots$,truth value$(\xi^*_\ell \in 
{\underset\tilde {}\to \tau_{\beta^*_\ell}}),\ldots)
_{\ell < \omega}$ and $\beta^*_\ell \in B_1$ and then we can find $f \in
\text{ EAUT}(\bar Q)$ such that 
$\hat f(\underset\tilde {}\to g({\underset\tilde {}\to \tau_{\lambda + \xi}}
)) = \underset\tilde {}\to g({\underset\tilde {}\to \tau_{\lambda + \xi}})
\in M^*,\hat f(p_\zeta) \in M_\zeta$ is compatible with $p_\zeta$, and we get
an easy contradiction.

Note: in $\bold V[\langle {\underset\tilde {}\to \tau_\beta}:\beta \in B_1
\rangle]$ we can compute $\underset\tilde {}\to E[G_{P_{\lambda + \xi+1}}]$,
but we do not have $\underset\tilde {}\to g
({\underset\tilde {}\to \tau_{\lambda + \xi}})$ by the previous sentence so
easy contradiction as in the proof of \scite{1.12}.  \nl
${{}}$  \hfill$\square_{\scite{2.1}}$\margincite{2.1}
\newpage

\head {\S3 The $L_{\aleph_1,\aleph_1}$-elementary submodels and the forcing}
\endhead  \resetall 
\bn
We may wonder what is really needed in \S1, \S2 (and \cite{Sh:592}, \S2, \S3).
Here we generalize one feature: iterating with partial memory but without
transitivity in the memory: not restricting ourselves to Cohen and random.
\definition{\stag{3.1} Definition}  ${\Cal K}^s_\kappa$ is the family of
sequences of $\bar Q$ of this form

$$
\bar Q = \langle P_\alpha,{\underset\tilde {}\to Q_\alpha},A_\alpha,
B_\alpha,\bar \mu_\alpha,\bar \varphi_\alpha,{\underset\tilde {}\to Y_\alpha},
{\underset\tilde {}\to \tau_\alpha}:\alpha < \beta \rangle
$$
\mn
(we write $\beta = \ell g(\bar Q))$ such that:
\mr
\item "{$(a)$}"  $\langle P_\alpha,{\underset\tilde {}\to Q_\alpha}:\alpha
< \beta \rangle$ is an FS iteration of c.c.c. forcing notions (so $P_\beta$
denotes the limit)
\sn
\item "{$(b)$}"  $\mu_\alpha < \kappa$ (see clause (d) below),
${\underset\tilde {}\to \tau_\alpha}$ is a 
${\underset\tilde {}\to Q_\alpha}$-name (i.e. a $P_\alpha$-name of one),
$\Vdash_{P_{\alpha +1}} ``{\underset\tilde {}\to \tau_\alpha} \subseteq
\mu_\alpha$ and ${\underset\tilde {}\to \tau_\alpha}$ is the generic of
${\underset\tilde {}\to Q_\alpha}"$
\sn
\item "{$(c)$}"  $B_\alpha \subseteq A_\alpha \subseteq \alpha$ and $\beta
\in B_\alpha \Rightarrow B_\beta \subseteq B_\alpha$ and $|B_\alpha| <
\kappa$.
\endroster
\bn
\ub{First simple version}:
\mr
\item "{$(d)$}"  $\bar \mu = \langle
\mu^0_\alpha,\mu^1_\alpha,\mu^2),\bar \varphi = \langle
\varphi^0_\alpha,\varphi^1_\alpha \rangle,\varphi^0_\alpha =
\varphi^0_\alpha(x,{\underset\tilde {}\to Y_\alpha}),\varphi^1_\alpha =
\varphi^1_\alpha(x,y,{\underset\tilde {}\to Y_\alpha})$, let
$\mu_\alpha = \mu^2_\alpha$
\sn
\item "{$(e)$}"  the set of elements of ${\underset\tilde {}\to
Q_\alpha}$ is a non empty subset of $Q^{\text{pos}}_\alpha =
\{\eta:\eta \in {}^\omega(\mu^0_\alpha)$ in the universe $\bold
V[\langle
{\underset\tilde {}\to \tau_\gamma}:\gamma \in A_\alpha \rangle]$
\sn
\item "{$(f)$}"  ${\underset\tilde {}\to Y_\alpha}$ is a
$P_\alpha$-name of a subset of $\mu^1_\alpha$, moreover we have
sequences $\bar{\Cal B}_\alpha = \langle {\Cal B}_{\alpha,\zeta}:\zeta
< \mu^1_\alpha \rangle,\bar \xi_\alpha = \langle
\xi_{\alpha,\zeta,n}:\zeta < \mu^1_\alpha,n < \omega \rangle$ and
$\bar \beta_{\alpha,\zeta} = \langle \beta_{\alpha,\zeta,n}:\zeta <
\mu^1_\alpha,n < \omega \rangle$ all three from $\bold V$ such that:
{\roster
\itemitem{ $(\alpha)$ }  ${\Cal B}_{\alpha,\zeta} = {\Cal
B}_{\alpha,\zeta}(\ldots,x_n,\ldots)_{n < \omega}$ is a Borel function
from ${}^\omega\{\text{true,false}\}$ to $\{\text{true,false}\}$
\sn
\itemitem{ $(\beta)$ }  $\xi_{\alpha,\zeta,n} <
\mu_{\beta_{\alpha,\zeta,n}},\beta_{\alpha,\zeta,n} \in B_\alpha$
\sn
\itemitem{ $(\gamma)$ }  $\underset\tilde {}\to Y = \{\zeta <
\mu^1_\alpha:{\Cal B}_{\alpha,\zeta}(\ldots$,truth
value$(\xi_{\alpha,\zeta,n} \in 
{\underset\tilde {}\to \tau_{\beta_{\alpha,\zeta,n}}}),
\ldots)_{n < \omega} = \text{ truth}\}$
so ${\underset\tilde {}\to Y_\alpha} \in \bold V[\langle
{\underset\tilde {}\to \tau_\gamma}:\gamma \in B_\alpha \rangle)$
\endroster}
\sn
\item "{$(g)$}"  $p \in {\underset\tilde {}\to Q_\alpha}$ iff $\{\eta
\in {\underset\tilde {}\to Q^{\text{pos}}_\alpha}:\bold V[\langle
{\underset\tilde {}\to \tau_\gamma}:\gamma \in B_\alpha \rangle]
\models \varphi^0_\alpha(p,\underset\tilde {}\to Y)\}$ and
${\underset\tilde {}\to Q_\alpha} \models p \le q$ iff $\bold V[\langle
{\underset\tilde {}\to \tau_\gamma}:\gamma \in B_\alpha \rangle]
\models \varphi^1_\alpha[p,q,\underset\tilde {}\to Y]$.
\endroster
\bn
\ub{Second, non simple version}:
\mr
\item "{$(d)$}"  $\bar \mu_\alpha = \langle \mu^\ell_\alpha:\ell < 7 \rangle$,
we let $\mu_\alpha = \mu^4_\alpha,\bar \varphi_\alpha = \langle
\varphi^\ell_\alpha:\ell < 6 \rangle,\varphi_\alpha =
\varphi^0_\alpha,
\varphi^\ell_\alpha$ is a 3-place
relation on ${}^{\omega >}(\mu^\ell_\alpha)$ if $\ell = 0,3$, a 4-place relation
on ${}^{\omega >}(\mu^\ell_\alpha)$ if $\ell = 1,2,\varphi^4_\alpha$ is
binary, $\varphi^5_\alpha$ is a 5-place relation, and for $i<6$ we
have $\varphi^i_\alpha$ is a subset \footnote{note: lim$(\varphi^0_\alpha)$
describes $Q_\alpha$, lim$(\varphi^1_\alpha)$, lim$(\varphi^2_\alpha)$
describes $\le_{Q^\alpha},+_{Q_\alpha}$, lim$(\varphi^3_\alpha)$ describes
the maximal antichain of $Q_\alpha$, lim$(\varphi^4_\alpha)$ describes
``$\underset\tilde {}\to \xi \in {\underset\tilde {}\to \tau_\alpha}$" and
lim$(\varphi^5_\alpha)$ describes $G_\alpha(2)$} of $\{\bar \eta:\bar \eta =
(\eta_0,\dotsc,\eta_{j-1})$, for some $n$ each $\eta_\ell$ is a sequence of
ordinals of length $n$, and $(i,j) \in \{(0,3),(1,4),(2,4),(3,3),(4,2),
(5,5)\}\}$ (can read more natural restrictions), $\varphi^i_\alpha$ is closed
under initial segments (i.e. letting $i=0, \bar \eta = (\eta_0,\eta_1,\eta_2)
\in \varphi^0_\alpha \Rightarrow \bar \eta \restriction n = (\eta_0
\restriction n,\eta_1 \restriction n,\eta_2 \restriction n) \in \varphi^0
_\alpha$, of course $\bar \eta \restriction n$ is an abuse of notation),
$$
\align
\text{lim}(\varphi_\alpha) = \{(\eta,\nu,\rho):&\eta,\nu,\rho
\text{ are sequences of length } \omega \\
  &\text{such that } n < \omega \Rightarrow (\eta \restriction n,\nu
\restriction n,\rho \restriction n) \in \varphi_\alpha\}
\endalign
$$
(depend on the universe) and similarly for other $\varphi^i_\alpha$
\sn
\item "{$(e)$}"  the set of elements of $Q_\alpha$ is a (nonempty) subset
of $Q^{\text{pos}}_\alpha = \{\eta:\eta \in {}^\omega(\mu^0_\alpha) \cap
\bold V[{\underset\tilde {}\to \tau_\gamma}:\gamma \in A_\alpha]\}$
\sn
\item "{$(f)$}"  ${\underset\tilde {}\to {\bar \rho}_\alpha} =
\langle {\underset\tilde {}\to \rho^\ell_{\alpha,\zeta}}:\zeta < \zeta^\ell
_\alpha,\ell < 6 \rangle,{\underset\tilde {}\to \rho^\ell_{\alpha,\zeta}} =
\langle {\underset\tilde {}\to \rho^\ell_{\alpha,\zeta,m}}:m < \omega 
\rangle$, and ${\underset\tilde {}\to \rho^\ell_{\alpha,\zeta,m}}$ is a
$P_\alpha$-name of a natural number (or of a member of $\{0,1\}$), it has the form
${\Cal B}^\ell_{\alpha,\zeta,m}(\ldots$,truth value$(\xi^\ell_{\alpha,\zeta,m,k} \in
{\underset\tilde {}\to \tau_{\beta^\ell_{\alpha,\zeta,m,k}}}),\ldots))_{k < \omega},
{\Cal B}^\ell_{\alpha,\zeta,m}$ is a Borel function, $\beta^\ell_{\alpha,\zeta,m,k} \in
B_\alpha$ and $\xi^\ell_{\alpha,\zeta,m,k} < \mu_\alpha$
\sn
\item "{$(g)$}"  in $\bold V[\langle {\underset\tilde {}\to \tau_\beta}:\beta \in
A_\alpha \rangle],{\underset\tilde {}\to Q_\alpha}$ is the set of
elements of $\{\eta \in
{\underset\tilde {}\to Q^{\text{pos}}_\alpha}$: there are $\nu \in
{}^\omega(\text{Ord})$ and $\zeta < \zeta^0_\alpha$ such that $(\eta,\nu,
\rho^0_{\alpha,\zeta}) \in \text{ lim}(\varphi^0_\alpha)$ and $\nu(0) =
\zeta$; so $\eta,\nu,{\underset\tilde {}\to \rho^0_{\alpha,\zeta}}
[{\underset\tilde {}\to G_{P_\alpha}}]$ belongs to $\bold V[\langle
{\underset\tilde {}\to \tau_\beta}:\beta \in A_\alpha \rangle]\}$, we may
write $\zeta = \zeta^0_\alpha(p)$.  Similarly $\le_{\underset\tilde
{}\to Q_\alpha}$ is defined by $\varphi^1_\alpha$
\sn
\item "{$(h)$}"  in $V[\langle {\underset\tilde {}\to \tau_\beta}:\beta \in
A_\alpha \rangle]$ we have $\{(p,q):p \le_{{\underset\tilde {}\to Q_\alpha}}
q\} = \{(p,q):p,q \in {\underset\tilde {}\to Q_\alpha}$ and there are $\nu \in
{}^\omega(\text{Ord})$ and $\zeta < \zeta^1_\alpha$ such that $(p,q,\nu,
\rho^1_{\alpha,\zeta}) \in \text{ lim}(\varphi^1_{\alpha,\zeta}) \in
\text{lim}(\varphi^1_\alpha),\nu(0) = \zeta\}$ \nl
We write $\zeta^1_\alpha(p,q)$ for the minimal $\zeta$
\sn
\item "{$(i)$}"  in $V[\langle {\underset\tilde {}\to \tau_\beta}:\beta \in
A_\alpha \rangle]\{(p,q):p,q \in {\underset\tilde {}\to Q_\alpha}$
are compatible in ${\underset\tilde {}\to Q_\alpha}\} = \{(p,q):p,q \in
{\underset\tilde {}\to Q_\alpha}$ and there are
$\nu \in {}^\omega(\text{Ord})$ and $\zeta < \zeta^2_\alpha$ such that 
$(p,q,\nu,\rho^2_{\alpha,\zeta}) \in \text{lim}(\varphi^2_\alpha)$ and 
$\nu(0) = \zeta\}$; 
\nl
(we write $\zeta^2_\alpha(p,q)$ for the minimal such $\zeta$)
\sn
\item "{$(k)$}"  in $V[\langle {\underset\tilde {}\to \tau_\beta}:\beta \in
A_\alpha \rangle]$ we have $\{ \langle p_n:n < \omega \rangle:p_n \in
{\underset\tilde {}\to Q_\alpha}$ for $n < \omega,\{p_n:n < \omega\}$ is
predense in ${\underset\tilde {}\to Q_\alpha}\} = \{ \langle p_n:n < \omega
\rangle$: for some $\nu \in {}^\omega(\text{Ord})$ and 
$\zeta < \zeta^3_\alpha$ we have $(p,\nu,\rho^3_{\alpha,\zeta}) \in 
\text{ lim}(\varphi^3_\alpha)$ [where $p = \langle p_{n - [\sqrt n]^2}
([\sqrt n]):n < \omega]\}$, and we write $\zeta^2_\alpha(\bar p)$ for the
minimal such $\zeta$
\sn
\item "{$(l)$}"  in $V[\langle {\underset\tilde {}\to \tau_\beta}:\beta \in
A_\alpha \rangle]$ we have $\{(p,\xi):p \in {\underset\tilde {}\to Q_\alpha},
\xi < \mu^4_\alpha$ and $p \Vdash_{\underset\tilde {}\to Q_\alpha} ``\xi \in
{\underset\tilde {}\to \tau_\alpha}"\} = \{(p,\xi)$: for some 
$\nu \in {}^\omega(\mu^4_\alpha)$ and $\zeta < \zeta^4_\alpha$ we have 
$\left< \langle \xi \rangle \char 94 p,\langle \zeta \right> \in 
\varphi^4_\alpha\}$ 
\sn
\item "{$(m)$}"  if $G_{< \alpha} \subseteq P_\alpha$ is generic over $\bold V,
G_\alpha \subseteq Q_\alpha[G_{< \alpha}]$ is generic over $\bold V
[G_{< \alpha}]$ and $\tau_\alpha = {\underset\tilde {}\to \tau_\alpha}
[G_\alpha]$, \ub{then} in $V[\langle {\underset\tilde {}\to \tau_\beta}
[G_\alpha]:\beta \in A_\alpha \rangle]$ we have
\endroster
$$
\align
G_\alpha = \{\eta \in {}^\omega(\text{Ord}):&\text{there are } \nu \in
{}^\omega(\text{Ord) and } \zeta < \zeta^5_\alpha \text{ and} \\
  &\nu^1 \in {}^\omega(\mu_\alpha), \text{ and } \rho \in {}^\omega 2
\text{ such that}: \\
  &\rho(n) = 0 \Leftrightarrow \nu^1(n) \in \tau_\alpha \text{ and} \\
  &(\eta,\nu,\nu^1,\rho,\rho_{\alpha,\zeta}) \in \text{ lim}(\varphi^5
_\alpha)\}.
\endalign
$$
\enddefinition
\bn
We can simplify Definition \scite{3.1} by
\definition{\stag{3.2} Definition}  1) We say ${\Cal B}$ is a Borel function
to $\bold V$, if it is a function from ${}^\omega\{0,1\}$ to $\bold V$
(identifying sometimes 1 with truth, 0 with false) with countable range,
such that each $\{\eta \in {}^\omega 2:{\Cal B}(\eta) = x\}$ is a Borel set.
\nl
2) We say ${\Cal B}$ is a Borel function to ${}^\omega \bold V$ if
${\Cal B}(\eta) = \langle {\Cal B}_n(\eta):n < \omega \rangle$ with each
${\Cal B}_n$ being a Borel function to $\bold V$. \nl
3) $\bar Q$ is simple if in Definition \scite{3.1} clause (e)

$$
\text{lim}(\varphi^0_\alpha) \subseteq \{(\eta,\nu,\rho):\rho \in
{}^\omega \omega\}.
$$
\enddefinition
\bigskip

\definition{\stag{3.3} Definition} Let $\bar Q \in {\Cal K}^s_\kappa,\alpha =
\ell g(\bar Q)$. \nl
1) Let

$$
\align
P'_\beta = \biggl\{p \in P_\beta:&\text{for every } \gamma \in \text{ dom}(p),
p(\gamma) \text{ is computed by some} \\
  &\text{Borel function } {\Cal B} \text{ (so } {\Cal B} \in \bold V, \text{ and is
an object not a name),} \\
  &{\Cal B} \text{ is a function to } {}^\omega \bold V, \text{ from some} \\
  &\langle \text{truth value} (\xi_\ell \in
{\underset\tilde {}\to \tau_{\alpha_\ell}}):\ell < \omega \rangle
\text{ where} \\
  &\xi_\ell < \mu_{\alpha_\ell} \text{ and } \alpha_\ell < \gamma \text{ and we
say in this case that} \\
  &\text{the truth value of } (\xi_\ell \in
{\underset\tilde {}\to \tau_{\alpha_\ell}}) \text{ appear in } p(\gamma)
\text{ for } \ell < \omega \\
  &\text{and } p \restriction \gamma \text{ forcing a value to }
\zeta^0_\gamma(p(\gamma)) \biggr\}.
\endalign
$$
\noindent
2) For $p \in P'_\beta,\alpha \in \text{ dom}(p)$ let

$$
\text{supp}(p(\alpha)) = \{\gamma:\text{for some } \xi \text{ the truth value
of } (\xi \in {\underset\tilde {}\to \tau_\gamma}) \text{ appear in }
p(\alpha)\}
$$
\noindent
and let supp$(p) = \text{ dom}(p) \cup \dbcu_{\alpha \in \text{ dom}(p)}
\text{ supp}(p(\alpha))$. \nl
3) For $A \subseteq \alpha$ let $P'_A = \{p \in P'_\alpha:\text{dom}(p)
\subseteq A$ and $\gamma \in \text{ dom}(p) \Rightarrow \text{ supp}(p(\gamma))
\subseteq A\}$, i.e. $P'_A = \{p \in P'_\alpha:\text{supp}(p) \subseteq A\}$
with the order inherited from $P'_\alpha$ which is inherited from
$P_A$ (recall: only for some $A$'s,
$P'_A \lessdot P'_\alpha$). \nl
4)  $A$ is called $\bar Q$-closed if $A \subseteq \ell g(\bar Q)$, and 
$\alpha \in A \Rightarrow B_\alpha \subseteq A$.  We call $A$ strongly
$\bar Q$-closed if $A \subseteq \ell g(\bar Q),\alpha \in A \Rightarrow
A_\alpha \subseteq A$.  We let $c \ell_{\bar Q}(A)$ be the $\bar Q$-closure
of $A$ and $sc \ell_{\bar Q}(A)$ be the strong $\bar Q$-closure of $A$. \nl
5) Let
\sn
\ub{Simple version}:

$$
\align
\text{PAUT}(\bar Q) = \biggl\{ f:&(i) \quad f 
\text{ is a one-to-one function}, \\
  &(ii) \quad \text{domain and range of } f \text{ are }
\bar Q \text{-closed} \\
  &\quad \quad \text{(in particular are } \subseteq \ell g(\bar Q)) \\
  &(iii) \quad \text{for } \alpha_1,\alpha_2 \in \text{ dom}(f)
\text{ we have:} \\
  &\quad \quad \,\alpha_1 \in A_{\alpha_2} \Leftrightarrow f(\alpha_1) \in
A_{f(\alpha_2)} \\
  &(iv) \quad \text{for } \alpha_1,\alpha_2 \in \text{ dom}(f)
\text{ we have} \\
  &\quad \quad \alpha_1 \in B \alpha_2 \Leftrightarrow (f_1) \in
B_{f(\alpha_2)} \\
  &(v) \quad \text{if } f(\alpha_1) = \alpha_2 \text{ then} \\
   &\quad \quad f \text{ maps } {\underset\tilde {}\to Y_{\alpha_1}}
\text{ to } {\underset\tilde {}\to Y_{\alpha_2}}, \text{ i.e. } \bar
\mu_{\alpha_1} = \bar \mu_{\alpha_2}, \\
  &\quad \quad \bar\varphi_{\alpha_1} = \bar\varphi_{\alpha_2},
{\Cal B}_{\alpha_1,\zeta} = {\Cal B}_{\alpha_2,\zeta}, \\
  &\quad \quad \xi_{\alpha_1,\zeta,n} =
\xi_{\alpha_2,\zeta,n},\beta_{\alpha_1,\zeta,n} =
\beta_{\alpha_2,\zeta,n} \biggr\}.
\endalign
$$
\sn
\ub{Non simple version}:

$$
\align
\text{PAUT}(\bar Q) = \biggl\{ f:&(i) \quad f 
\text{ is a one-to-one function}, \\
  &(ii) \quad \text{domain and range of } f \text{ are }
\bar Q \text{-closed} \\
  &\quad \quad \text{(in particular are } \subseteq \ell g(\bar Q)) \\
  &(iii) \quad \text{for } \alpha_1,\alpha_2 \in \text{ dom}(f)
\text{ we have}: \\
  &\quad \quad \,\alpha_1 \in A_{\alpha_2} \Leftrightarrow f(\alpha_1) \in
A_{f(\alpha_2)} \\
  &(iv) \quad \text{for } \alpha_1,\alpha_2 \in \text{ dom}(f)
\text{ we have} \\
  &\quad \quad \alpha_1 \in B_{\alpha_2} \Leftrightarrow (f_1) \in
B_{f(\alpha_2)} \\
  &(v) \quad \text{if } f(\alpha_1) = \alpha_2 \text{ then} \\
  &\quad \quad \mu^\ell_{\alpha_1} = \mu^\ell_{\alpha_2},\zeta^\ell
_{\alpha_1} = \zeta^\ell_{\alpha_2},\varphi^\ell_{\alpha_1} = \varphi^\ell
_{\alpha_2} \\
  &\quad \quad \text{and } f \text{ maps }
{\underset\tilde {}\to \rho^\ell_{\alpha_1,\zeta}} \text{ to }
{\underset\tilde {}\to \rho^\ell_{\alpha_2,\zeta}} \\
  &\quad \quad \text{(i.e. } {\Cal B}^\ell_{\alpha_1,\zeta,m} = {\Cal
B}^\ell_{\alpha_2,\zeta,m} \\
  &\quad \quad \text{ for } \ell < 6,\zeta < \zeta^\ell_{\alpha_1} =
\zeta^\ell_{\alpha_2},m < \omega \\
  &\quad \quad \text{ and } \xi^\ell_{\alpha_1,\zeta,m,k} =
\xi^\ell_{\alpha_2,\zeta,m,k} \\
 &\quad \quad \text{ for } \ell < 6,\zeta < \zeta^\ell_{\alpha_1} =
\zeta^\ell_{\alpha_2},m < \omega,k < \omega \text{ but, of course } \\
  &\quad \quad ``\xi \in {\underset\tilde {}\to \tau_\beta}"
\text{ is replaced by } ``\xi \in {\underset\tilde {}\to
\tau_{f(\beta)}}" \\
  &\quad \quad \text{ that is } \beta^\ell_{\alpha_2,\zeta,m,k} =
f(\beta^\ell_{\alpha_1,\zeta,m,k}) \biggr\}.
\endalign
$$
6) For $f \in \text{ PAUT}(\bar Q),\hat f$ is the natural map which $f$
induces from $P'_{\text{dom}(f)}$ onto $P'_{\text{rang}(f)}$;
similarly for part (7). \nl
7) For $\bar Q^1,\bar Q^2 \in {\Cal K}^s_\kappa$ let PAUT$(\bar Q^1,\bar Q^2)$
be defined similarly:

$$
\align
\text{PAUT}(\bar Q^1,\bar Q^2) = \biggl\{ f:&(i) \quad f 
\text{ is a one-to-one function}, \\
  &(ii) \quad \text{domain of } f \text{ is } \bar Q^1\text{-closed} \\
  &\quad \quad \text{(in particular is } \subseteq \ell g(\bar Q^1)) \\
  &\quad \quad \text{and range of } f \text{ is } \bar Q^2\text{-closed} \\
  &\quad \quad \text{(in particular is } \subseteq \ell g(\bar Q^2)) \\
  &(iii) \quad \text{for } \alpha_1,\alpha_2 \in \text{ dom}(f)
\text{ we have}: \\
  &\quad \quad \,\,\alpha_1 \in A^1_{\alpha_2} \Leftrightarrow f(\alpha_1) \in
A^2_{f(\alpha_2)} \\
  &(iv) \quad \text{for } \alpha_1,\alpha_2 \in \text{ dom}(f)
\text{ we have} \\
  &\quad \quad \,\,\alpha_1 \in B^1_{\alpha_2} \Leftrightarrow f(\alpha_1) \in
B^2_{f(\alpha_2)} \\
  &(v) \quad \text{if } f(\alpha_1) = \alpha_2 \text{ then} \\
  &\quad \quad \text{ \ub{simple version}}: \\
  &\quad \quad {}^1 \bar \mu_{\alpha_1} = 
{}^2 \bar \mu_{\alpha_1},{}^1 \bar \varphi_{\alpha_1} = {}^2 \bar
\varphi_{\alpha_2},\\
  &\quad \quad  \text{ and } f \text{ maps } {\underset\tilde {}\to Y^1_{\alpha_1}}
\text{ to } {\underset\tilde {}\to Y^2_{\alpha_2}}, \\
  &\quad \quad \text{ i.e. }
{\Cal B}^1_{\alpha_1,\zeta} = {\Cal B}^2_{\alpha_2,\zeta},\xi^1_{\alpha_1,\zeta,m} = 
\xi^2_{\alpha_2,\zeta,m},\beta^1_{\alpha_1,\zeta,m} =
\beta^2_{\alpha_2,\zeta,m} \\
    &\quad \quad \text{ \ub{non simple version}}: \text{ (for } \ell < 6) \\
    &\quad \quad {}^1\mu^\ell_{\alpha_1} = {}^2\mu^\ell_{\alpha_2},
{}^1\zeta^\ell_{\alpha_1} = {}^2\zeta^\ell_{\alpha_2},{}^1\varphi^\ell
_{\alpha_1} = {}^2\varphi^\ell_{\alpha_2} \\
  &\quad \quad \text{and } f \text{ maps }
{}^1{\underset\tilde {}\to \rho^\ell_{\alpha_1,\zeta}} \text{ to }
{}^2{\underset\tilde {}\to \rho^\ell_{\alpha_2,\zeta}} \\
  &\quad \quad \text{(i.e. } 
   {\Cal B}^\ell_{\alpha_2,\zeta,m} = {\Cal B}^\ell_{\alpha_1,\zeta,m} \\
  &\quad \quad \xi^\ell_{\alpha_2,\zeta,m,k} =
\xi^\ell_{\alpha_1,\zeta,m,k} \\
 &\quad \quad \beta^\ell_{\alpha_2,\zeta,m,k} =
f(\beta^\ell_{\alpha_1,\zeta,m,k}) \text{ for } \zeta < \zeta^\ell_{\alpha_1} =
\zeta^\ell_{\alpha_2},m < \omega,k < \omega) \biggr\}.
\endalign
$$
\enddefinition
\bigskip

\proclaim{\stag{3.4} Claim}  Let $\bar Q \in {\Cal K}^s_\kappa$ be of
length $\alpha^*$. \nl
1) $P'_{\alpha^*}$ is a dense subset of $P_{\alpha^*}$. \nl
2) In $\bold V^{P_\alpha}$, from ${\underset\tilde {}\to \tau_\alpha}
[G_{Q_\alpha}]$ we can reconstruct $G_{Q_\alpha}$ and vice versa.  From
$\langle {\underset\tilde {}\to \tau_\gamma}:\gamma < \alpha \rangle
[G_{P_\alpha}]$ we can reconstruct $G_{P_\alpha}$ and vice versa.  So
$\bold V^{P_\alpha} = \bold V[\langle {\underset\tilde {}\to \tau_\beta}:
\beta < \alpha \rangle]$. \nl
3) If $\mu$ is any cardinal, and $\underset\tilde {}\to X$ is a
$P_{\alpha^*}$-name of a subset of $\mu$, \ub{then} there is a set $A \subseteq
\alpha^*$ such that $|A| \le \mu$ and $\Vdash_{P_{\alpha^*}}
``\underset\tilde {}\to X \in \bold V[\langle
{\underset\tilde {}\to \tau_\gamma}:\gamma \in A \rangle]"$.  Moreover, for
each $\zeta < \mu$ there is in $\bold V$ a Borel function ${\Cal B}_\zeta(x_0,
\dotsc,x_n,\ldots)_{n < \omega}$ with domain and range the set
$\{$true,false$\}$ and $\gamma_\ell \in A,\xi_{\zeta,\ell} < \mu_{\gamma_\zeta}$
for $\ell < \omega$ i.e. $\langle {\Cal
B}_\zeta,\gamma_{\zeta,\ell},\xi_{\zeta,\ell}:\zeta < \mu,\ell <
\omega \rangle \in \bold V$) such that

$$
\Vdash_{P_{\alpha^*}} ``\zeta \in \underset\tilde {}\to X \text{ iff true }
= {\Cal B}_\zeta(\ldots,``\text{truth value of } \xi_{\zeta,\ell} \in
{\underset\tilde {}\to \tau_{\gamma_{\zeta,\ell}}}[G_{Q_{\gamma_{\zeta,\ell}}}]",\ldots)".
$$
\endproclaim
\bigskip

\demo{Proof}  1) Easy. \nl
2), 3) By induction on $\alpha$. \nl
4) Let $\chi^*$ be such that $\{\bar Q,\lambda\} \in {\Cal H}(\chi^*)$, and
let $\zeta < \mu$; let $M$ be an elementary submodel of $({\Cal H}(\chi^*),
\in,<^*_{\chi^*})$ to which $\{\bar Q,\lambda,\kappa,\mu,
\underset\tilde {}\to X,\zeta\}$ belongs and $\mu \subseteq M$, so
$\Vdash_{P_{\alpha^*}} ``M[{\underset\tilde {}\to G_{P_{\alpha^*}}}] \cap
{\Cal H}(\chi^*) = M"$.  Hence by \scite{3.4}(3) (i.e. as
$\bold V^{P_\alpha} = \bold V[\langle {\underset\tilde {}\to \tau_\beta}:
\beta < \alpha \rangle])$ we have $M[{\underset\tilde {}\to G_{P_{\alpha^*}}}]
= M[\langle \tau_i:i \in \alpha \cap M \rangle]$ and the conclusion should
be clear.    \hfill$\square_{\scite{3.4}}$\margincite{3.4}
\enddemo
\bigskip

\proclaim{\stag{3.4A} Claim}   For $A \subseteq \alpha^*$, every real in $\bold V[\langle
{\underset\tilde {}\to \tau_\gamma}:\gamma \in A \rangle]$ of even
subset $X$ of $\mu$ (some $\mu$) has the form mentioned in
\scite{3.4}(4) with $\gamma_{\zeta,\ell} \in A$.
\endproclaim
\bigskip

\demo{Proof}  This does not follow by \scite{3.4}(4) as e.g. maybe
$\neg P'_A \lessdot P_\alpha$.  Let ${\Cal L}_{\lambda^+,\omega}$
denote the propositional logic, allowing conjunctions and
disjunctions of size $\le \lambda$.  We know that if $X \subseteq
\mu,X \in \bold V[ \langle \tau_\gamma:\gamma \in A \rangle]$ where
${\underset\tilde {}\to \tau_\gamma} \subseteq \mu_\gamma$ (and
$\langle \mu_\gamma:\gamma \in A \rangle \in \bold V)$, \ub{then} we can
find in $\bold V$ a sequence $\langle {\Cal B}_\zeta:\zeta < \mu
\rangle,\langle (\xi_{\zeta,i},\gamma_{\zeta,i}):\zeta < \mu,i <
\mu_\zeta \rangle,\xi_{\gamma,i} < \mu_{\gamma_{\zeta,i}},{\Cal
B}_\zeta = {\Cal B}_\zeta(\ldots,x_i,\ldots)_{i < \mu_\zeta} \in {\Cal
L}_{\lambda^+,\omega}$ for some $\lambda$ such that:
\mr
\item "{$(*)$}"  $X = \{\zeta < \mu:{\Cal B}_\zeta(\ldots,\text{truth
value}(\xi_{\zeta,i} \in \tau_{\gamma_{\zeta,i}}),\ldots)_{i < \mu_i} = \text{ truth}\}$.
\ermn
So if $p \in P_\alpha,p \Vdash ``\underset\tilde {}\to \tau \in \mu"$,
we can find $q,p \le q \in P_\alpha$ and $\langle({\Cal
B}_\zeta,\mu_\zeta,\xi_{\zeta,i},\gamma_{\zeta,i}):\zeta < \mu,i <
\mu_\zeta \rangle \in \bold V$ as above with $\gamma_{\zeta,i} \in A$
such that $q \Vdash_{P_\alpha} ``\underset\tilde {}\to \tau = \{\zeta
< \mu:{\Cal B}_\zeta(\ldots,\text{truth value}(\xi_{\zeta,i} \in
{\underset\tilde {}\to \tau_{\gamma,\zeta,i}}),\ldots)_{i <
\mu_\zeta}\}$.

Now as $P_\alpha$ satisfies the c.c.c., for each $\zeta$ separately we
can replace the conjunction and disjunction inside ${\Cal B}_\zeta$ by
countable ones, so we are done.  \hfill$\square_{\scite{3.5}}$\margincite{3.5} 
\enddemo
\bigskip

\proclaim{\stag{3.5} Claim}  1) If $\bar Q^1,\bar Q^2 \in {\Cal K}^s_\kappa$
and $f \in \text{ PAUT}(\bar Q^1,\bar Q^2)$ and ${\text{\rm dom\/}}(f) =
\ell g(\bar Q^1)$, \nl
${\text{\rm rang\/}}(f) = \ell g(\bar Q^2)$ \ub{then} $\hat f$
is an isomorphism from $(P'_{\ell g(\bar Q^1)})^{(\bar Q^1)}$ onto
$(P'_{\ell g(\bar Q^2)})^{(\bar Q^2)}$. \nl
2) If $\bar Q \in {\Cal K}^s_\kappa$ and $A \subseteq \ell g(\bar Q)$ is
strongly $\bar Q$-closed \ub{then} $P'_A \lessdot P'_{\ell g(\bar Q)}$, in fact
if $q \in P'_{\ell g(\bar Q)},q \restriction A \le p \in P'_A$ \ub{then} $p \cup
(q \restriction (\ell g(\bar Q \backslash A)) \in P'_{\ell g(\bar Q)}$ is
a lub of $p,q$ and there are unique $(\bar Q',f)$ such that $\bar Q' \in
{\Cal K}^s_\kappa,f \in \text{ PAUT}(\bar Q',\bar Q),f$ is order preserving,
dom$(f) = \ell g(\bar Q')$, rang$(f) = A$. \nl
3) If $A^0 \subseteq \alpha,A^{n+1} = A^n \cup \{A_\alpha:\alpha \in A^n\}$
\ub{then} $\dbcu_{n < \omega} A^n$ is strongly $\bar Q$-closed and $\bigcup
\{\alpha+1:\alpha \in A^0\} = \bigcup\{\alpha +1:\alpha \in \dbcu_{n < \omega}
A^n\}$. \nl
4) PAUT$(\bar Q)$ is closed under composition and inverse.  Similarly if
$f_\ell \in {\text{\rm PAUT\/}}(\bar Q^\ell,\bar Q^{\ell +1})$ for 
$\ell =1,2$, then $f_2 \circ f_1 \in {\text{\rm PAUT\/}}
(\bar Q^1,\bar Q^3),f^{-1}_i \in {\text{\rm PAUT\/}}(\bar Q^2,\bar
Q^1)$. \nl
5) If $\bar Q \in {\Cal K}^s_\kappa$ and $A \subseteq \ell g(\bar Q)$
is strongly $\bar Q$-closed and $\bar p \in P'_{\ell g(\bar Q)}$ then
$p \restriction A \in P'_A$ and $P_{\ell g(\bar Q)} \models ``(p
\restriction A) \le p"$ and $q \in P'_A \and q \le_{P_{\ell g(\bar
Q)}} p \Rightarrow q \le_{P_{\ell g(\bar Q)}} p \restriction A$.
\endproclaim
\bigskip

\demo{Proof}  Straightforward.
\enddemo
\bigskip

\definition{\stag{3.6} Definition}  We say $\bar {\Cal F}$ witnesses $A$ for
$\bar Q,\kappa$ and $A'$ (or for $(\bar Q,\kappa,A')$, we may omit $A'$ if
it is the strong $\bar Q$-closure of $A$) \ub{if}:
\mr
\item "{$(a)$}"  $\bar Q \in {\Cal K}^s_\kappa$
\sn
\item "{$(b)$}"  $A \subseteq \ell g(\bar Q),A \subseteq A' \subseteq
\ell g(\bar Q),A$ is $\bar Q$-closed, $A'$ is strongly $\bar Q$-closed
\sn
\item "{$(c)$}"  $\bar{\Cal F} = \langle {\Cal F}_\gamma:\gamma < \omega_1 \rangle,
{\Cal F}_\gamma$ decreasing with $\gamma$ (this just for notational simplicity)
\sn
\item "{$(d)$}"  ${\Cal F}_\gamma \subseteq \text{ PAUT}(\bar Q)$ for $\gamma <
\omega_1$
\sn
\item "{$(e)$}"  if $f \in {\Cal F}_\gamma$ then
{\roster
\itemitem{ $(\alpha)$ }   rang$(f) \subseteq A$,
\sn
\itemitem{ $(\beta)$ }  $|\text{dom}(f)| < \kappa$
\sn
\itemitem{ $(\gamma)$ }  $(\forall \alpha \in \text{ dom}(f))(\exists
\beta \in \text{ dom}(f))[f(\beta) = \beta \and \alpha \in \text{
scl}_{\bar Q}\{\beta\}]$ and
\sn
\itemitem{ $(\delta)$ }  dom$(f) \subseteq A'$
\endroster}
\sn
\item "{$(f)$}"  if $f_1 \in {\Cal F}_{\gamma_1},\gamma_2 < \gamma_1 < \omega_1$
and $C \subseteq A' (\subseteq \ell g(\bar Q))$ has cardinality $< \kappa$
then for some $f_2 \in {\Cal F}_{\gamma_2}$ we have: $f_1 \subseteq f_2$ and
$C \subseteq \text{ dom}(f_2)$ and $C \cap A \subseteq \text{ rang}(f_2)$
\sn
\item "{$(g)$}"  if $C \subseteq A$ is such that $|C| < \kappa$ then
id$_C \in \dbca_{\gamma < \omega_1} {\Cal F}_\gamma$.
\endroster
\enddefinition
\bigskip

\demo{\stag{3.7} Observation}  Let $\bar Q \in {\Cal K}^s_\kappa$. \nl
1) If $\beta^* \le \ell g(\bar Q)$ then $\bar Q \restriction \beta^* \in
{\Cal K}^s_\kappa$, and for $A \subseteq \beta^*$ we have $(P'_A)^{\bar Q
\restriction \beta^*} = (P'_A)^{\bar Q}$. \nl
2) If $\bar{\Cal F}$ witnesses $A$ for $\bar Q,\kappa$ in $A'$ and $A'' \subseteq
A' \subseteq \ell g(\bar Q),A \subseteq A'',A''$ is strongly $\bar Q$-closed
\ub{then} witnesses $A$ for $\bar Q,\kappa$ and $A''$.
\nl
3) If $\bar{\Cal F}$ witnesses $A$ for $\bar Q,\kappa$ in $A'$ and $\beta^* \le
\ell g(\bar Q)$ and $A''$ is the strong $\bar Q$-closure of $A \cap \beta^*$
then $\langle \{f \restriction \cup\{\text{scl}_{\bar
Q}\{\beta\}:\beta \in \text{ dom}(f) \cap \beta^*\}):f \in {\Cal F}_{1
+ \gamma}\}:\gamma < \omega_1 \rangle$ witnesses $A \cap \beta^*$ for $\bar Q,\kappa$ in
$A''$ (why $1 + \gamma$? to preserve also
``transitive closure of $\gamma_1 \in A_{\gamma_2}$"; see
\scite{3.10}(1) below). 
\enddemo
\bigskip

\demo{Proof}  Straightforward.
\enddemo
\bigskip

\proclaim{\stag{3.8} Claim}  Let $\bar Q \in {\Cal K}^s_\kappa$. \nl
1) If $\bar{\Cal F}$ witnesses $A$ for $\bar Q,\kappa$ inside $A' =
\text{ scl}_{\bar Q}(A)$, \ub{then}
\mr
\item "{$(a)$}"  $P'_A \lessdot P'_{\ell g(\bar Q)}$
\sn
\item "{$(b)$}"  $p,q \in P'_A$ are compatible in $P'_{\ell g(\bar Q)}$ \ub{iff}
they are compatible in $P'_A$.  \nl
Also if $p,q \in P'_A,\gamma \in \text{ dom}(p) \cap \text{ dom}(q)$ \ub{then}
$$
q \restriction \gamma \Vdash_{P'_\gamma} ``p(\gamma)
\le_{\underset\tilde {}\to Q_\gamma} q(\gamma)" \Rightarrow q \restriction
\gamma \Vdash_{P'_{\gamma \cap A}} ``p(\gamma)
\le_{\underset\tilde {}\to Q_\gamma} q(\gamma)"
$$
\sn
\item "{$(c)$}"  ${\Cal I} \subseteq P'_A$ is 
predense in $P'_{\ell g(\bar Q)}$ \ub{iff} it is predense in 
$P'_A$ (\wilog \, ${\Cal I}$ is countable)
\sn
\item "{$(d)$}"  there are \footnote{in fact, this does not depend on
$A$ having a witness (but not necessarily the ``moreover".}
unique $f,\bar Q'$ such that $\bar Q' \in
{\Cal K}^s_\kappa,f \in \text{ PAUT}(\bar Q',\bar Q),f$ order preserving,
dom$(f) = \ell g(\bar Q')(= \text{ otp}(A))$ and rang$(f) = A$,
moreover $\hat f$ is an isomorphism from $(P'_{\ell g(\bar
Q)'})^{{\bar Q}'}$ onto $(P'_A)^{\bar Q}$.
\ermn
2) Moreover
\mr
\item "{$(e)$}"  for every $p,q \in P'_{A'}$ there is $\gamma < \omega_1$
such that: 
{\roster
\itemitem{ $(i)$ }  if $f \in F_\gamma$ and supp$(p) \cup \text{ supp}(q)
\subseteq \text{ dom}(f)$ \ub{then} $P'_{A'} \models ``p \le q" \Leftrightarrow
P'_{\ell g(\bar Q')} \models ``\hat f(p) \le \hat f(q)"$
\sn
\itemitem{ $(ii)$ }  if $f \in F_\gamma$ and supp$(p) \cup \text{ supp}(q)
\subseteq \text{ rang}(f)$ \ub{then}
$$
P'_{A'} \models ``\hat f^{-1}(p) \le \hat f^{-1}(q)" \Leftrightarrow
P'_{A'} \models ``p \le q"
$$
\endroster}
\item "{$(f)$}"  for every $p,q \in P'_{A'}$ there is $\gamma < \omega_1$
such that:
{\roster
\itemitem{ $(i)$ }  if $f \in {\Cal F}_\gamma$ and supp$(p) \cup \text{ supp}(q)
\subseteq \text{ dom}(f)$ \ub{then} $p,q$ are compatible in $P'_{\ell g(\bar Q)}$
iff $\hat f(p),\hat f(q)$ are compatible in $P'_{\ell g(\bar Q')}$
\sn
\itemitem{ $(ii)$ }  if $f \in {\Cal F}_\gamma$ and supp$(p) \cup \text{ supp}(q)
\subseteq \text{ rang}(f)$ (so $p,q \in P'_A$) \ub{then} 
$\hat f^{-1}(p),\hat f^{-1}(q)$ are
compatible in $P'_{A'}$ iff $p,q$ are compatible in $P'_{A'}$
\endroster}
\item "{$(g)$}"  for every countable ${\Cal I} \subseteq P'_{A'}$, there is
$\gamma < \omega_1$ such that:
{\roster
\itemitem{ $(i)$ }  if $f \in {\Cal F}_\gamma$ and $\dbcu_{p \in {\Cal I}}
\text{ supp}(p) \subseteq \text{ dom}(f)$ \ub{then} ${\Cal I}$ is predense in
$P'_{A'}$ iff $\hat f({\Cal I})$ is predense in $P'_{A'}$
\sn
\itemitem{ $(ii)$ }  if $f \in {\Cal F}_\gamma$ and $\dbcu_{p \in {\Cal I}}
\text{ supp}(p) \subseteq \text{ rang}(f)$ (so ${\Cal I} \subseteq
P'_A$) \ub{then} $\{\hat f^{-1}(p):p \in
{\Cal I}\}$ is predense in $P'_{A'}$ iff ${\Cal I}$ is predense in $P'_{A'}$
\endroster}
\item "{$(h)$}"  for every Borel function ${\Cal B} = {\Cal B}(\ldots$,
truth value $(\xi_n \in {\underset\tilde {}\to \tau_{\alpha_n}}),\ldots)
_{n < \omega}$ from ${}^\omega 2$ to $\bold V$ such that $\alpha_n \in A'$
and $p \in P'_{A'}$ there is $\gamma < \omega_1$ such that:
{\roster
\itemitem{ $(i)$ }  if $f \in {\Cal F}_\gamma$ and supp$(p) \cup \{\alpha_n:n <
\omega\} \subseteq \text{ dom}(f)$ and $x \in \bold V$ \ub{then} \nl
$p \Vdash_{P'_{\ell g(\bar Q)}} ``{\Cal B}(\ldots,\text{truth value}
\, (\xi_n \in {\underset\tilde {}\to \tau_{\alpha_n}}),\ldots)_{n < \omega}
= x" \text{ iff}$ \nl
$\hat f(p) \Vdash_{P'_{\ell g(\bar Q)}} ``{\Cal B}(\ldots,\text{truth value}
\, (\xi_n \in {\underset\tilde {}\to \tau_{f(\alpha_n)}}),\ldots)_{n < \omega}
= x"$.
\sn
\itemitem{ $(ii)$ }  if $f \in {\Cal F}_\gamma$ and supp$(p) \cup \{\alpha_n:n <
\omega\} \subseteq \text{ rang}(f)$ and $x \in \bold V$ \ub{then} \nl
$p \Vdash_{P'_{\ell g(\bar Q)}} ``{\Cal B}(\ldots,\text{truth value}
\, (\xi_n \in {\underset\tilde {}\to \tau_{\alpha_n}}),\ldots)_{n < \omega}
= x" \text{ iff}$ \nl
$\hat f^{-1}(p) \Vdash_{P'_{\ell g(\bar Q)}} ``{\Cal B}(\ldots,
\text{truth value} \, (\xi_n \in 
{\underset\tilde {}\to \tau_{f^{-1}(\alpha_n)}}),\ldots)_{n < \omega}= x"$. 
\endroster}
\endroster
\endproclaim
\bigskip

\demo{Proof}  We prove ((1) + (2) together) by induction on $\alpha^*_0 =
\alpha^* = \ell g(\bar Q)$. \nl
Arriving to $\alpha^*$, we note various implications
\mr
\item "{$(*)_1$}"  for $\alpha < \alpha^*$, replacing $A,A'$ and $\bar{\Cal F}$
by $A \cap \alpha,A''$ (any strong $\bar Q$-closed set such that $A \cap
\alpha \subseteq A'' \subseteq A' \cap \alpha)$ and $\bar{\Cal F}
\upharpoonleft A'' =: \langle \{f \restriction B:B \subseteq \text{
dom}(f) \cap A''$ and $(\forall \alpha \in B)(\exists \beta \in
B)[f(\beta) = \beta \and \alpha \in \text{ scl}_{\bar Q}(\{\beta\}):f
\in {\Cal F}_i\}:i < \omega_1 \rangle$ respectively, Claim \scite{3.8}
holds
\sn
\item "{$(*)_2$}"  clause (d) (for $\alpha^*$) follows (from the induction
hypothesis)
\sn
\item "{$(*)_3$}"  if clause (e) holds \ub{then} clause (b) holds.
\ermn
[Why?  The second phrase in clause (b) holds by clause (e)(ii) as $C
\subseteq A \and |C| \le \aleph_0 \and i < \omega_1 \Rightarrow \text{
id}_C \in {\Cal F}_i$.  
If $p,q$ are compatible in $P'_A$ then they have a common upper bound
there, which ``works" in $P'_{\ell g(\bar Q)}$, too, by the previous
sentence.  Assume $p,q$ have a
common upper bound $r$ in $P'_{\ell g(\bar Q)}$, so by \scite{3.5}(5) also
$r' =: r \restriction A' \in P'_{A'}$ is a common upper bound of $p$ and $q$.

Let $\gamma_1,\gamma_2 < \omega_1$ be as guaranteed in clause (e) for
$p \le r',q \le r'$ respectively and let $\gamma = \text{ max}\{\gamma_1,
\gamma_2\}$ and let $f = \text{ id}_{\text{supp}(p) \cup \text{supp}(q)}$, now
$f \in F_{\gamma +1}$ by clause (g), Definition \scite{3.6} and there is
$f',f \subseteq f' \in {\Cal F}_\gamma$ such that supp$(r') \subseteq \text{ dom}
(f')$.

So $\hat f'(p) = \hat f(p) = p,\hat f'(q) = \hat f(q) = q,\hat f'(r') \in P'_{A}$,
and by clause (e) we have $p = \hat f'(p) \le \hat f'(r') \in P'_{A},
q = \hat f'(q) \le \hat f'(r') \in P'_{A}$ and we are done.]  
\mr
\item "{$(*)_4$}"  if clauses (e),(f) hold then clause (c) holds.
\ermn
[Why?  If ${\Cal I}$ is not predense in $P'_A$ then there is $q \in P'_A$
incompatible in $P'_A$ with every $p \in {\Cal I}$ hence by clause (b) which
holds by $(*)_3,q$ is incompatible with $p$ in $P'_{\ell g(\bar Q)}$, hence
${\Cal I}$ is not predense in $P'_{\ell g(\bar Q)}$.  Next assume ${\Cal I}$
is predense in $P'_A$, let ${\Cal J}$ be a maximal antichain of $P'_A$ of
elements above some member of ${\Cal I}$, so by clause (b) (which holds by
$(*)_3$) ${\Cal J}$ is an antichain in $P'_{\ell g(\bar Q)}$ hence is
countable.

Let $q$ be any member of $P'_{A'}$, let $\gamma = \text{ sup}\{\gamma(p,q):
p \in {\Cal I}\},\gamma(p,q)$ as in clause (f).

Now letting $C = \dbcu_{p \in {\Cal J}} \text{ supp}(p)$, it is a countable
subset of $A$ hence $f_0 = \text{ id}_C \in F_{\gamma +1}$, so there is
$f,f_0 \subseteq f \in F_\gamma$ such that supp$(q) \subseteq \text{ dom}(f)$.
Now for $p \in {\Cal J},\gamma \ge \gamma(p,q)$, hence we have that 
\sn
$(\forall p \in {\Cal J})(p,q$ are compatible in $P_{A'}$ iff $\hat f(p),
\hat f(q)$ are compatible in $P'_{\ell g(\bar Q)})$. 
\sn
As $\hat f(q) \in P'_A$ and ${\Cal J}$ is a maximal antichain of $P'_A$ for
some $p \in {\Cal J},\hat f(q),p$ are compatible in $P'_A$, hence in
$P'_{A'}$.  But $f(p) = p$ so $q,p$ are compatible in $P'_{\ell g(\bar Q)}$,
but by the choice of ${\Cal J}$ there is $p' \in {\Cal I},p' \le p$, so $q$
is compatible with some member of ${\Cal J}$, i.e. with $p'$.  As $q \in
P'_{A'}$ was arbitrary, this proves that ${\Cal I}$ is predense in $P'_{A'}$,
completing the second implication in the proof of $(*)_4$.]
\mr
\item "{$(*)_5$}"  if clauses (e),(f) hold then clause (a) holds.
\ermn
[Why?  By $(*)_3 + (*)_4$.]
\mr
\item "{$(*)_6$}"  if clauses (e),(f),(g) hold for $\alpha^*$ \ub{then}
clause (h) holds.
\ermn
[Why?  First, it is enough to deal with the case the range of the Borel
function is $\{0,1\}$.  Now we prove the assertion by induction on the depth of the Borel
function.
\bn
\ub{Case A}:  ${\Cal B}$ is atomic, i.e. ${\Cal B}$ is the truth value of
$\xi \in {\underset\tilde {}\to \tau_\gamma}$ for some $\xi < \mu_\gamma$.

Clearly $\gamma < \alpha^*$ so we can apply the induction hypothesis
to $\gamma$.  So $p \Vdash ``{\Cal B}(\xi \in
{\underset\tilde {}\to \tau_\gamma}) = i"$ where $p \in P'_{A'}$ is
equivalent to $p \restriction \gamma \Vdash ``p(\gamma) = \text{
truth}"$.
Now $p(\gamma)$ has the form ${\Cal B}'(\ldots,\text{truth
value}(\xi_n \in {\underset\tilde {}\to \tau_{\beta_n}}),\ldots)_{n <
\omega}$ where $\beta_n \in A_\gamma,\xi_n < \mu_{\beta_n},{\Cal B}'
\in \bold V$ a Borel function, so the statement is equivalent to $p
\restriction \gamma \Vdash ``{\Cal B}'(\ldots,\text{truth value}
(\xi_n \in {\underset\tilde {}\to \tau_{\beta_n}}),\ldots)_{n < \omega}
= \text{ truth}$.  As $A'$ is strongly $\bar Q$-closed clearly
$\beta_n \in A_\gamma \subseteq A'$ so we can apply the induction
hypothesis on $\gamma$ using clause (h) there.
\bn
\ub{Case B}:  ${\Cal B} = \neg{\Cal B}'$ (i.e. $1 - {\Cal B}$).

By the way we phrase the statement it follows from the statement on
${\Cal B}'$.
\bn
\ub{Case C}:  ${\Cal B} = \dsize \bigwedge_{n < \omega} {\Cal B}_n$.

Let ${\Cal I}$ be a maximal subset of

$$
\align
\bigl\{ q \in P'_{A'}:&(i) \quad q \text{ forces } {\Cal B}_n
\text{ for every } n \text{ \ub{or} for some } n,q \text{ forces } \neg
{\Cal B}_n \\
  &(ii) \quad p \le q \text{ or } q,p \text{ are incompatible} \bigr\}
\endalign
$$
\mn
which is an antichain in $P_{\ell g(\bar Q)}$.  Let $\gamma_{\Cal I} <
\omega_1$ be as guaranteed by clause (g) (of \scite{3.8}), let for $q \in
{\Cal I}:\gamma(p,q) < \omega_1$ be as guaranteed by clause (f) 
(of \scite{3.8}) if $p,q$ are incompatible and as guaranteed by clause (e)
(of \scite{3.8}) if $p \le q$.

For each $q \in {\Cal I}$ let $\gamma_n(q)$ be the $\gamma$ guaranteed for
$q,{\Cal B}_n$ for both $x=0$ and $x=1$.  
For $q_1 \ne q_2$ from ${\Cal I}$ (so incompatible) let
$\gamma(q_1,q_2) < \omega_1$ be as guaranteed in clause (f).
Let $\gamma^* = \sup(\{\gamma_{\Cal I}\} 
\cup \{\gamma(p,q):q \in {\Cal I}\} \cup \{\gamma_n(q):q
\in {\Cal I},n < \omega\} \cup \{\gamma(q_1,q_2):q_1 \ne q_2$ from
${\Cal I}\}) +1$.

Suppose $f_1 \in F_{\gamma^*}$, supp$(p) \cup \{\alpha_n:n < \omega\}
\subseteq \text{ dom}(f_1)$, where the $\alpha_n$'s are from the
statement of clause (h).  We can find $f_2,f_1 \subseteq f_2 \in
F_{\gamma^*-1}$ and $\dbcu_{q \in {\Cal I}} \text{ dom}(q) \subseteq
\text{ dom}(f_2)$.

Let ${\Cal I}_0 = \{q:p \le q \in {\Cal I}\},{\Cal I}_1 = {\Cal I} \backslash
{\Cal I}_0$.  Now $q_1 \ne q_2 \in {\Cal I} \Rightarrow
\gamma(q_1,q_2) \le \gamma^* -1$ and 
$\gamma_{\Cal I} \le \gamma^* - 1$ hence ${\Cal I}' =
\{\hat f_2(q):q \in {\Cal I}\}$ is a maximal antichain.  Also $\gamma(p,q)
\le \gamma_1$ hence

$$
q \in {\Cal I}'_0 =: \{\hat f_2(q);q \in {\Cal I}_0\} \Rightarrow \hat f_2
(p) \le q
$$

$$
q \in {\Cal I}'_1 = \{\hat f_1(q):q \in {\Cal I}_1\} \Rightarrow \hat f_2
(p),q \text{ incompatible}.
$$
\mn
Also for $q \in {\Cal I}_0$,

$$
q \Vdash ``{\Cal B}_n = 0" \Leftrightarrow f_2(q) \Vdash ``f_2({\Cal B}_n) = 0"
$$

$$
q \Vdash ``{\Cal B}_n = 1" \Leftrightarrow f_2(q) \Vdash ``f_2({\Cal B}_n) = 1.
$$
\mn
The rest should be clear.]
\bn
So together it is enough to prove clauses (e), (f), (g).  
\mr
\item "{$(*)_7$}"  clause (e) holds.
\ermn
If Dom$(p) = \emptyset$ this is trivial so suppose not, and if dom$(p)
\nsubseteq \text{ dom}(q)$ the equivalence is trivial so assume
$\emptyset \ne \text{ dom}(p) \subseteq \text{ dom}(q)$.  Let
$\alpha^\otimes = \text{ max}(\text{ dom}(q))$, now if $\alpha^\otimes
+ 1 < \alpha^*$ we can use the induction hypothesis on $\gamma^* +1$,
so assume $\alpha^* = \alpha^\otimes +1$ (we can also discard the case
$\gamma^* \notin \text{ Dom}(p)$ if we like). \nl
Now $P'_{A'} \models ``p \le q"$ iff $(\alpha) + (\beta)$ where
\mr
\item "{$(\alpha)$}"  $P'_{A' \cap \gamma^*} \models ``p \restriction
\alpha^\otimes \le q \restriction \alpha^\otimes$
\sn
\item "{$(\beta)$}"  $q \Vdash ``p(\alpha^\otimes) \le q"$.
\ermn
Now for $(\alpha)$ get $\gamma_1$ by applying clause (e) to $\bar Q
\restriction \alpha^\otimes$, and for $(\beta)$ get $\gamma_2$ by
applying clause (h) to $\bar Q \restriction \alpha^\otimes$, so
$\gamma = \text{ max}\{\gamma_1,\gamma_2\}$ is as required.
\mr
\item "{$(*)_8$}"  clause (f) holds.
\ermn
As in the proof of clause (e), \wilog \, $\alpha^* = \alpha^* +1$ and
$\alpha^\otimes = \text{ max}(\text{dom}(p)) = \max(\text{dom}(q))$.
Let $\{r_n:n < \omega\} \subseteq P'_{A_{\alpha^\otimes}}$ be a
maximal antichain, such that each $r_n$ satisfies
\mr
\item "{$(\alpha)$}"  it forces a truth value say $\bold t_n$ to
$``p(\alpha^\otimes),q(\alpha^\otimes)$ are compatible in
${\underset\tilde {}\to Q_{\alpha^\otimes}}"$
\sn
\item "{$(\beta)$}"  $r_n \ge p$ or $r_n,p$ are incompatible
\sn
\item "{$(\gamma)$}"  $r_n \ge q$ or $r_n,q$ are incompatible.  
\ermn
By applying clause (h) to $\bar Q \restriction \alpha^\otimes$, \wilog
\, $\{r_n:n < \omega\} \subseteq P'_A$.

Now clearly $p \le q$ \ub{iff} $\dsize \bigvee_n [p \le r_n \and q \le
r_n \and \bold t_n = \text{ truth}]$, and we can apply the induction
hypothesis to each of those countably many statements.
\mr
\item "{$(*)_9$}"  clause (g) holds.
\ermn
Without loss of generality ${\Cal I}$ is countable.  Now if $w =
\cup\{\text{dom}(p):p \in {\Cal I}\}$ has no last element, clearly
${\Cal I}$ is predense iff $\alpha \in w \Rightarrow {\Cal
I}^{[\alpha]} = \{p \restriction \alpha:p \in {\Cal I}\}$ is predense;
so we can finish by the induction hypothesis.  So assume
$\alpha^\otimes$ is the last element in $w$.  Let $\{r_n:n < \omega\}
\subseteq P'_{A' \cap \alpha^\otimes}$ be a maximal antichain of it,
each $r_n$ forcing a truth value to ``$\{p(\alpha^\otimes):p
\restriction \alpha^\otimes \in {\underset\tilde {}\to
G_{P_{\alpha^\otimes}}}\}$ is a maximal antichain of 
$Q_\alpha$".  \hfill$\square_{\scite{3.8}}$\margincite{3.8}
\enddemo
\bn
This follows from the claim below.
\definition{\stag{3.9} Definition}  1) We say $\bar F$ is a $\kappa$-witness
for $(\bar Q^1,\bar Q^2)$ if
\mr
\item "{$(a)$}"  $\bar Q^1 \in {\Cal K}^s_\kappa$ and $\bar Q^2 \in
{\Cal K}^s_\kappa$ (so $\kappa > \aleph_0$)
\sn
\item "{$(b)$}"  $\bar F = \langle F_\gamma:\gamma < \omega_1 \rangle$ and
$F_\gamma$ is decreasing with $\gamma$ (for notational simplicity)
\sn
\item "{$(c)$}"  $F_\gamma \subseteq \text{ PAUT}(\bar Q^1,\bar Q^2)$ for
$\gamma < \omega_1$
\sn
\item "{$(d)$}"  if $f \in F_\gamma$ \ub{then} dom$(f) \subseteq \ell g
(\bar Q^1)$ and rang$(f) \subseteq \ell g(\bar Q^2)$
\sn
\item "{$(e)$}"  if $f_1 \in F_{\gamma_1},\gamma_2 < \gamma_1 < \omega_1,
C_1 \subseteq \ell g(\bar Q^1),C_2 \subseteq \ell g(\bar Q^2),|C_1| <
\kappa,|C_2| < \kappa$ \ub{then} for some $f_2 \in F_{\gamma_2}$ we have
$$
f_1 \subseteq f_2,C_1 \subseteq \text{ dom}(f_2),C_2 \subseteq \text{ rang}
(f_2).
$$
\ermn
2) We say $\bar F$ is an explicit $\kappa$-witness for $(\bar Q^1,\bar Q^2)$
if (a) - (d) above hold and
\mr
\item "{$(e)'$}"  if $f_1 \in F_{\gamma_1},\gamma_2 < \gamma_1 < \omega_1,
C_1 \subseteq \ell g(\bar Q^1) \cap \text{ scl}_{(\bar Q^1)}(\text{dom}(f_1)),
C_2 \subseteq \ell g(\bar Q^2) \cap \text{ scl}_{(\bar Q^2)}(\text{dom}(f_2)),
|C_1| < \kappa,|C_2| < \kappa$, \ub{then} for some $f_2 \in F_{\gamma_2}$ we
have $f_1 \subseteq F_2,C_1 \subseteq \text{ dom}(f_2),C_2 \subseteq
\text{ rang}(f_2),\text{dom}(f_2) \subseteq \text{ scl}_{\bar Q^1}(\text{dom}
(f_1)),\text{rang}(f_2) \subseteq \text{ scl}_{\bar Q^2}(\text{rang}(f_1))$
\sn
\item "{$(f)$}"  for $f \in F_\gamma,g \subseteq f \Rightarrow g \in F_\gamma$
\sn
\item "{$(g)$}"  if $f \in F_\gamma$ and $\alpha,\beta \in \text{ dom}(f)$, \ub{then}
$$
\alpha \in \text{ scl}_{\bar Q^1}(\{\beta\}) \Leftrightarrow f(\alpha) \in
\text{ scl}_{\bar Q^2}(\{f(\beta)\})
$$
\sn
\item "{$(h)$}"    if $C_1 \subseteq \ell g(\bar Q^1),C_2 \subseteq
\ell g(\bar Q^2)$ are countable and $\gamma < \omega_1$, \ub{then} for
some $f \in F_\gamma$ we have $C_1 \subseteq \text{ dom}(f),C_2
\subseteq \text{ rang}(f)$.
\endroster
\enddefinition
\bigskip

\proclaim{\stag{3.10} Claim}  Assume $\bar F$ is a $\kappa$-witness for
$(\bar Q^1,\bar Q^2)$. \nl
1) If $\gamma \in (0,\omega_1)$ and $f \in F_\gamma$ and $\alpha,\beta <
\ell g(\bar Q^1)$, \ub{then}

$$
\alpha \in \text{ scl}_{\bar Q^1}(\{\beta\}) \Leftrightarrow f(\alpha) \in
\text{ scl}_{\bar Q^2}(\{f(\beta)\}).
$$
\mn
2) Let $\bar F^{-1} = \langle F^{-1}_\gamma:\gamma < \omega_1 \rangle,
F^{-1}_\gamma =  \{f^{-1}:f \in F_\gamma\}$.  
\ub{Then} $\bar F^{-1}$is a $\kappa$-witness for $(\bar Q^2,\bar Q^1)$. \nl
3) If $A_\ell \subseteq \ell g(\bar Q^\ell)$ and: $A_\ell =
\text{ scl}_{\bar Q^\ell}(A_\ell)$ for $\ell =1,2$ and we let $F'_\gamma =
\{f \cap X:f \in F_{1 + \gamma},X \subseteq A_1 \times A_2\}$, \ub{then}
$\langle F'_\gamma:\gamma < \omega_1 \rangle$ is a $\kappa$-witness for
$(\bar Q^1 \restriction A_1,\bar Q^2 \restriction A_2)$ and an explicit
$\kappa$-witness for it, note that by renaming $A_\ell$ can become an
ordinal. \nl
4) If $p \in (P^1_{\ell g(\bar Q^1)})',q \in (P^1_{\ell g(\bar Q^1)})'$ and
${\Cal I} \subseteq (P^1_{\ell g(\bar Q^1)})'$ is countable, \ub{then} for
some $\gamma < \omega_1$ we have: for every $f \in F_\gamma$ we have:
\mr
\item "{$(\alpha)$}"  if supp$(p) \cup \text{ supp}(q) \subseteq
\text{ dom}(f)$ \ub{then} $(P^1_{\ell g(\bar Q^1)})' \models p \le q$ iff
$(P^2_{\ell g(\bar Q^2)})' \models \hat f(p) \le \hat f(q)$
\sn
\item "{$(\beta)$}"  if supp$(p) \cup \text{ supp}(q) \subseteq
\text{ dom}(f)$ \ub{then} $p,q$ are compatible in 
$(P^1_{\ell g(\bar Q^1)})'$ iff $\hat f(p),\hat f(q)$ are incompatible
in $(P^2_{\ell g(\bar Q^1)})'$
\sn
\item "{$(\gamma)$}"  if $\dbcu_{r \in {\Cal I}}$ supp$(r) \subseteq
\text{ dom}(f)$ \ub{then} ${\Cal I}$ is predense in 
$(P^1_{\ell g(\bar Q^1)})'$ iff $\hat f({\Cal I}) = \{\hat f(r):r \in
{\Cal I}\}$ is predense in $(P^2_{\ell g(\bar Q^2)})'$.
\ermn
5) For every Borel function ${\Cal B} = {\Cal B}(\ldots,\bold t_n,\ldots)
_{n \in \omega}$ (for $\bold t_n$ a truth value, values in $\bold V$) and
$\langle (\xi_n,\alpha_n):n < \omega \rangle,\alpha_n < \ell g(\bar Q^1)$,
and $p \in (P^1_{\ell g(\bar Q^1)})'$ there is $\gamma < \omega_1$ such
that: if $f \in F_\gamma,\{\alpha_n:n < \omega\} \cup \text{ supp}(p) \subseteq
\text{ dom}(f)$ and $x \in \bold V$ \ub{then}

$$
p \Vdash_{P^1_{\ell g(\bar Q^1)}} ``{\Cal B}(\ldots,\text{truth value}
(\xi_n \in {\underset\tilde {}\to \tau_{\alpha_n}}),\ldots)_{n < \omega}
= x" \text{ iff}
$$

$$
\hat f(p) \Vdash_{P^2_{\ell g(\bar Q^2)}} ``{\Cal B}(\ldots,\text{truth value}
(\xi_n \in {\underset\tilde {}\to \tau_{f(\alpha_n)}}),\ldots)_{n < \omega}
= x".
$$
\endproclaim
\bigskip

\demo{Proof}  1) Easy. \nl
2) Trivial. \nl
3) Clear. \nl
4), 5) It suffices to prove \scite{3.11} below.
\enddemo
\bigskip

\proclaim{\stag{3.11} Claim}  If $\bar F$ is an explicit $\kappa$-witness for
$(\bar Q^1,\bar Q^2)$ then 4), 5) of \scite{3.10} holds.
\endproclaim 
\bigskip

\demo{Proof}  We prove this by induction on $\alpha^* = \text{ max}\{\ell g
(\bar Q^1),\ell g(\bar Q^2)\}$, and for a fixed $\alpha^*$ by induction on
$\beta^* = \text{ min}\{\ell g(\bar Q^1),\ell g(\bar Q^2)\}$.

As above it is enough to prove part (4). \nl
This is done by cases.  Without loss of generality $A'$ is the strong
$\bar Q$-closure of $A$.
\enddemo
\bn
\ub{Case 1}:  $\ell g(\bar Q^1) = 0$.

Trivial.
\bn
\ub{Case 2}:  cf$(\ell g(\bar Q^1)) > \aleph_0$.

So let $p,q,{\Cal I}$ be given and choose $\alpha^1 < \ell g(\bar Q^1)$ such
that $p,q \in (P^1_{\alpha^1})',{\Cal I} \subseteq (P^1_{\alpha^1})'$.  Let
$\bar Q^0 = \bar Q^1 \restriction \alpha^1,F'_\gamma = F_\gamma \cap
(\alpha^1 \times \ell g(\bar Q^2))$, and apply the induction hypothesis to
$\bar Q^0,\bar Q^2,\langle F'_\gamma:\gamma < \omega_1 \rangle$ and to
$p,q,{\Cal I}$ and get $\gamma < \omega_1$, and we shall prove that it works
for $\bar Q^1,\bar Q^2,\bar F,p,q,{\Cal I}$.  So let $f \in F_\gamma$ be such
that supp$(p) \cup \text{ supp}(q) \subseteq \text{ dom}(f)$ (needed if we
are dealing with clauses $(\alpha)$ or $(\beta)$ (of \scite{3.10}(4)) and
such that $\dbcu_{r \in cI} \text{ supp}(r) \subseteq \text{ dom}(f)$ (needed
if we are dealing with clause $(\gamma)$ of \scite{3.10}(4)).  Now $f' =:
f \restriction \alpha^1 \in F'_\gamma$ so:
\mr
\item "{$(*)_1$}"  $(P^1_{\alpha^1})' \models ``p \le q" \Leftrightarrow 
(P^2_{\ell g(\bar Q^2)})' \models ``\hat f'(p) \le \hat f'(q)"$
\sn
\item "{$(*)_2$}"  $p,q$ are compatible in $P^1_{\alpha^1}$ iff
$\hat f'(p),\hat f'(q)$ are compatible in $(P^2_{\ell g(\bar Q^2)})'$
\sn
\item "{$(*)_3$}"  ${\Cal I}$ is predense in $(P^1_{\alpha^1})'$ iff
$\{\hat f'(r):r \in {\Cal I}\}$ is predense in
$(P^2_{\ell g(\bar Q^2)})'$.
\ermn
As in all three cases we can replace $(P^1_{\alpha^1})'$ by
$(P^1_{\ell g(\bar Q^1)})'$ and $\hat f'$ by $\hat f$ we are done.
\bn
\ub{Case 3}:  cf$(\ell g(\bar Q^1)) = \aleph_0$.

Concerning clauses $(\alpha),(\beta)$ (of \scite{3.10}(4)) the proof is
just as in case 2, so we deal with clause $(\gamma)$.  Let ${\Cal I}
\subseteq (P^1_{\ell g(\bar Q^1)})'$ be countable.  We choose $\alpha^n <
\ell g(\bar Q^1)$ such that $\alpha^n < \alpha^{n+1}$ and $\dbcu_{n < \omega}
\alpha^n = \ell g(\bar Q^1)$ and let ${\Cal I}_n = \{r \restriction
\alpha^n:r \in {\Cal I}_n\}$ and $F^n_\gamma = F_\gamma \cap (\alpha^n \times
\ell g(\bar Q^2)),\bar F^n = \langle F^n_\gamma:\gamma < \omega_1 \rangle$.
For each $n$ apply induction hypothesis on $\bar Q^1 \restriction \alpha^n,
\bar Q^2,\bar F^n,{\Cal I}_n$ and get $\gamma_n < \omega_1$, and let
$\gamma^* = (\dbcu_{n < \omega} \gamma_n) +1$.  So let $f \in F_{\gamma^*}$
be such that $\dbcu_{r \in {\Cal I}} \text{ supp}(r) \subseteq \text{ dom}
(f)$. \nl
First assume ${\Cal I}$ is predense in $(P^1_{\ell g(\bar Q^1)})'$.

So there is $q \in (P^1_{(\bar Q^1)})'$ incompatible with every $r \in 
{\Cal I}$ hence for some $n,q \in (P^1_{\alpha_n})'$, so $q$ is incompatible
with every $r \in {\Cal I}_n$ in $(P^1_{\alpha_n})'$, hence ${\Cal I}_n$ is
not predense in $(P^1_{\alpha_n})'$.  As $\gamma^* > \gamma_n$, necessarily
$\{\hat f(r):r \in {\Cal I}_n\}$ is not predense in 
$(P^2_{\ell g(\bar Q^2)})'$, which means that some $q' \in
(P^1_{\ell g(\bar Q^2)})'$ is incompatible with $\hat f(r)$ for every
$r \in {\Cal I}_n$.  Now trivially $r \in {\Cal I} \Rightarrow \hat f(r
\restriction \alpha_n) \le_{(P^2_{\ell g(\bar Q^2)})'},\hat f(r)$ (look at
the definition and note $r \restriction \alpha_n \le r$ and increasing
$\gamma^*$), clearly $q' \in (P^2_{\ell g(\bar Q^2)})'$ is incompatible with
every memer of $\hat f({\Cal I}) =: \{\hat f(r):r \in {\Cal I}\}$, so
$\hat f({\Cal I})$ is not predense in $(P^2_{\ell g(\bar Q^2)})'$, as
required. \nl
Second assume ${\Cal I}$ is predense in $(P^1_{\ell g(\bar Q^1)})'$.  So for
each $n,{\Cal I}_n = \{r \restriction \alpha^n:r \in {\Cal I}\}$ is predense
in $(P^1_{\alpha_n})'$.

So by the induction hypothesis for each $n,\hat f({\Cal I}_n) = \{\hat f(r):
r \in {\Cal I}_n\}$ is a predense subset of $(P^2_{\ell g(\bar Q^2)})'$.

Let $A^1_n$ be the strong $\bar Q$-closure of $\dbcu_{p \in {\Cal I}_n}
\text{ supp}(p),A^1_\omega = \dbcu_{n < \omega} A^1_n$.

Let $A^2_n$ be the strong $\bar Q^2$-closure of $\dbcu_{p \in \hat f
({\Cal I}_n)} \text{ supp}(p),A^2_\omega = \dbcu_{n < \omega} A^2_n$, so if
$k < m \le \omega$ then $A^1_k \subseteq A^1_m$ and $A^2_k \subseteq A^2_m$,
hence $(P^1_{A^1_k})' \lessdot (P^1_{A^1_m})',(P^2_{A^2_k} \lessdot
P^2_{A^2_m})'$.

Remembering clause (g) of Definition \scite{3.9}(2), clearly
\mr
\item "{$\oplus$}"  if $\alpha \in \dbcu_{p \in {\Cal I}} \text{ supp}(p)
\subseteq \text{ dom}(f_1) \text{ then } \dsize \bigwedge_n(\forall n)
[\alpha \in A^1_n \leftrightarrow f_1(\alpha) \in A^2_n]$
\ermn
hence for $p \in {\Cal I},\hat f(p \restriction \alpha^n) = \hat f(p
\restriction A^1_n) = \hat f(p) \restriction A^2_n$.  As $p \in {\Cal I}
\rightarrow \dsize \bigvee_n(p \in {\Cal I}_n)$ (as dom$(p)$ is bounded in
$\ell g(\bar Q^1))$ clearly $\hat f({\Cal I}) = \dbcu_{n < \omega}
\hat f({\Cal I}_n)$.  Hence $q \in P'_{A^2_\omega}$ is incompatible in
$(P^2_{\ell g(\bar Q^2)})'$ with every $p \in \hat f({\Cal I})$ iff $q$
is incompatible with every $p \in \hat f({\Cal I}_{n(q)})$ in
$(P^2_{\ell g(\bar Q^2)})'$ where $n(q) = \text{ min}\{n:q \in 
P^2_{A^2_n}\}$ which is well defined as $\langle A^2_m:m \le \omega \rangle$
is increasing continuous.  As each $\hat f({\Cal I}_{n(q)})$ is predense in
$(P^1_{\ell g(\bar Q^2)})'$ (see above using the induction hypothesis),
there is no such $q$.  So $\hat f({\Cal I})$ is predense in
$P^2_{A^2_\omega}$, but $A^2_\omega$ is strongly $\bar Q^2$-closed (being
the union of such sets) hence (by \scite{3.7}) $(P^2_{A^2_\omega})' \lessdot
(P^2_{\ell g(\bar Q^2)})'$, so $\hat f({\Cal I})$ is predense, so we have
finished.
\bn
\ub{Case 4}: $\ell g(\bar Q^1) = \gamma^* + 1$.

By Definition \scite{3.1} all is translated to the $\gamma^*$ case
(including instances of \scite{3.10}(5)), so we can apply the induction
hypothesis.  \hfill$\square_{\scite{3.11}}$\margincite{3.11}
\bigskip

\proclaim{\stag{3.12} Claim}  If $\bar F$ witnesses $A$ for $\bar Q$ inside
$A'$, then
\mr
\item "{$(a)$}"  $P'_A \underset \bar {}\to \lessdot
P'_{A'} \lessdot P_{\ell g(\bar Q)}$
\sn
\item "{$(b)$}"  $p,q \in P'_A$ are compatible in $P'_{\ell g(\bar Q)}$ iff
they are compatible in $P'_A$
\sn
\item "{$(c)$}"  ${\Cal I} \subseteq P'_A$ is predense in 
$P'_{\ell g(\bar Q)}$iff it is predense in $P'_A$ (\wilog \, 
${\Cal I}$ is countable)
\sn
\item "{$(d)$}"  there are unique $f,\bar Q'$ such that $\bar Q' \in
{\Cal K}^s_\kappa,f \in \text{ PAUT}(\bar Q',\bar Q),f$ order preserving
dom$(f) = \ell g(\bar Q')(= \text{ otp}(A))$ and rang$(f) = A$; the 
essential part is: if $p,q \in P'_A,\gamma \in \text{ dom}(p) \cap
\text{ dom}(q)$ then
$$
q \restriction \gamma \Vdash_{P'_\gamma} p(\gamma)
\le_{\underset\tilde {}\to Q_\gamma} q(\gamma) \Rightarrow q \restriction
\gamma \Vdash_{P'_{\gamma \cap A}} p(\gamma) 
\le_{\underset\tilde {}\to Q_\gamma} q(\gamma).
$$
\endroster
\endproclaim
\bigskip

\remark{Concluding Remark}  We describe below an application of the method.
\endremark
\bigskip

\proclaim{\stag{3.13} Claim}  Assume for simplicity $\bold V \models GCH$.
For some c.c.c. forcing notion $\bold P$ of cardinality $\aleph_3$ we have:
in $\bold V^{\bold P}$
\mr
\item "{$(a)$}"  $2^{\aleph_0} = \aleph_3$
\sn
\item "{$(b)$}"  add (meagre) $= \aleph_1$
\sn
\item "{$(c)$}"  cov(meagre) $= \aleph_2$, moreover unif(meagre)
$= \aleph_2$
\sn
\item "{$(d)$}"  ${\frak d} = \aleph_3$.
\endroster
\endproclaim
\bigskip

\remark{\stag{3.14} Remark}  We can use other three cardinals, and use very
little cardinal arithmetic assumption.
\endremark
\bigskip

\demo{Proof}  Let $\bar Q = \langle P_i,{\underset\tilde {}\to Q_j},a_j,
{\underset\tilde {}\to \eta_i}:i \le \omega_3 + \omega_2,j < \omega_3 +
\omega_2 \rangle$ be FS iteration, with
\mr
\item "{$(a)$}"  ${\underset\tilde {}\to \eta_i}$ the generic real of
$Q_j$
\sn
\item "{$(b)$}"  if $j < \omega_3$ then ${\underset\tilde {}\to Q_j}$ is
Cohen say $({}^{\omega >}\omega,\triangleleft)$ so
${\underset\tilde {}\to \eta_j}$ is undominated
\sn
\item "{$(c)$}"  if $j = \omega_3 + \zeta,\zeta < \omega_2$ then $Q_j$ is
Random$^{\bold V[\langle \eta_i:i \in a_j \rangle]}$
\sn
\item "{$(d)$}"  if $b \subseteq a \in [\omega_3 + \omega_2]^{\le \aleph_1}$
then for arbitrarily large $i \in (\omega_3,\omega_3 + \omega_2)$ we have
$a \cap a_i = b$.
\ermn
Let $P = P_{\omega_3 + \omega_2}$.  Trivially $|P| = \aleph_3$.  So in
$\bold V$ we know:
\bn
\ub{Clause $(a)$}: $2^{\aleph_0} = \aleph_3$.

As $|P| \le \aleph_3$ the $\le$ inequality holds, the other inequality, $\ge$,
holds as $\langle {\underset\tilde {}\to \eta_i}:i < \omega_3 \rangle$ is a
sequence of $\aleph_3$ distinct reals.
\bn
\ub{Clause $(b)$}:  add(meagre) $= \aleph_1$.

As in \cite{Sh:592} we know that ${\frak b}$ is not increased by $\bold P$
so ${\frak b} = \aleph_1$ but add(meagre) $\le {\frak b}$.
\bn
\ub{Clause $(c)$}:  cov(meagre) $= \aleph_2$.

If for $i < \omega_1$ we have $A_i \subseteq {}^\omega 2$ are Borel sets in
$\bold V^{\bold P}$ then for some $\alpha < \omega_2, \langle A_i:i <
\omega_1 \rangle \in \bold V^{{\bold P}_{\omega_3 +i}}$, the forcing
$\bold P_{\omega_3 + i + \omega}/\bold P_{\omega_3 +i}$ add a Cohen real over
$\bold V^{{\bold P}_{\omega_3 +i}}$ (as we are using FS iteration).

So cov(meagre) $> \aleph_1$.

On the other hand $\langle {\underset\tilde {}\to \eta_{\omega_3 +i}}:i <
\omega_2 \rangle$ is a non-null set of reals hence unif(null) $\le \aleph_2$,
but cov(meagre) $\le$ unif(null) so cov(meagre) $\le \aleph_2$.  So together
cov(meagre) $= \aleph_2$.
\bn
\ub{Clause $(d)$}:  ${\frak d} = \aleph_3$.

Now ${\frak d} \le 2^{\aleph_0} = \aleph_3$, on the other hand, as in
\cite{Sh:592} for no $A \in [\omega_3]^{\aleph_2}$ from $\bold V$, is there
$<^*$-bound to $\{{\underset\tilde {}\to \eta_i}:i \in A\}$.
\enddemo

\newpage
    
REFERENCES.  
\bibliographystyle{lit-plain}
\bibliography{lista,listb,listx,listf,liste}

\def\germ{\frak} \def\scr{\cal} \ifx\documentclass\undefinedcs
  \def\bf{\fam\bffam\tenbf}\def\rm{\fam0\tenrm}\fi 
  \def\defaultdefine#1#2{\expandafter\ifx\csname#1\endcsname\relax
  \expandafter\def\csname#1\endcsname{#2}\fi} \defaultdefine{Bbb}{\bf}
  \defaultdefine{frak}{\bf} \defaultdefine{mathbb}{\bf}
  \defaultdefine{mathcal}{\cal}
  \defaultdefine{beth}{BETH}\defaultdefine{cal}{\bf} \def\bbfI{{\Bbb I}}
  \def\mbox{\hbox} \def\text{\hbox} \def\om{\omega} \def\Cal#1{{\bf #1}}
  \def\pcf{pcf} \defaultdefine{cf}{cf} \defaultdefine{reals}{{\Bbb R}}
  \defaultdefine{real}{{\Bbb R}} \def\restriction{{|}} \def\club{CLUB}
  \def\w{\omega} \def\exist{\exists} \def\se{{\germ se}} \def\bb{{\bf b}}
  \def\equivalence{\equiv} \let\lt< \let\gt> \def\cite#1{[#1]}
  \def\implies{\Rightarrow}
\begin{thebibliography}{GiSh 582}
\makeatletter \renewcommand{\@biblabel}[1]{[#1]} \makeatother
\def\eprintfootnotetext{References of the form {\tt math.XX/$\cdots$}
 refer to the {\tt xxx.lanl.gov} archive}
\ifx\documentstyle\undefinedcontrolsequence
   \def\anyfootnote{\footnote{*}}
   \else\def\anyfootnote{\footnote}\fi
\def\eprintfn{\ifEprint\anyfootnote{\eprintfootnotetext}\fi\Eprintfalse }
\newif\ifEprint  \Eprinttrue

\bibitem[Fe94]{Fe94}David Fremlin.
\newblock {Problem list}.
\newblock Circulated notes; available from {\tt
  http://www.essex.ac.uk/maths/staff/fremlin/measur.htm}.

\bibitem[GiSh 582]{GiSh:582}Moti Gitik and Saharon Shelah.
\newblock {More on real-valued measurable cardinals and forcing with ideals}.
\newblock {\em {Israel Journal of Mathematics}}, {\bf submitted}.
 {\tt math.LO/9507208}\eprintfn\

\bibitem[J]{J}Thomas Jech.
\newblock {\em {Set theory}}.
\newblock Academic Press, New York, 1978.

\bibitem[Ko]{Ko}Peter Komjath.
\newblock On second-category sets.
\newblock {\em Proc. Amer. Math. Soc.}, {\bf 107}:653--654, 1989.

\bibitem[Sh:b]{Sh:b}Saharon Shelah.
\newblock {\em {Proper forcing}}, volume 940 of {\em {Lecture Notes in
  Mathematics}}.
\newblock {Springer-Verlag, Berlin-New York, xxix+496 pp}, 1982.

\bibitem[Sh:f]{Sh:f}Saharon Shelah.
\newblock {\em {Proper and improper forcing}}.
\newblock {Perspectives in Mathematical Logic}. {Springer}, 1998.

\bibitem[Sh 592]{Sh:592}Saharon Shelah.
\newblock {Covering of the null ideal may have countable cofinality}.
\newblock {\em {Fundamenta Mathematicae}}, {\bf 166}:109--136, 2000.
 {\tt math.LO/9810181}\eprintfn\

\end{thebibliography}

\shlhetal

\enddocument

\bye